\documentclass{gtart}

\def\ifplaintex{\expandafter\ifx\csname documentclass\endcsname\relax}

\ifplaintex 
\else
\fi

\expandafter\ifx\csname beginpicture\endcsname\relax
\expandafter\ifx\csname documentclass\endcsname\relax
\input pictex \else
\input prepictex \input pictex \input postpictex \fi\fi

\def\gt{{\mathsurround=0pt\it $\cal G\mskip-2mu$eometry \&\ 
$\cal T\!\!$opology}}        

\def\gtp{{\mathsurround=0pt\it $\cal G\mskip-2mu$eometry \&\ 
$\cal T\!\!$opology $\cal P\!$ublications}}  

\def\volumenumber#1{\def\thevolumenumber{#1}}
\def\papernumber#1{\def\thepapernumber{#1}}
\def\volumeyear#1{\def\thevolumeyear{#1}}

\def\pagenumbers#1#2{\def\startpage{#1}\def\finishpage{#2}}
\def\published#1{\def\publishdate{#1}}
\def\proposed#1{\def\theproposer{#1}}
\def\seconded#1{\def\theseconders{#1}}
\def\received#1{\def\receiveddate{#1}}

\def\accepted#1{\def\accepteddate{#1}}

\def\asciiaddress#1{\def\theasciiaddress{#1}}
\def\asciiemail#1{\def\theasciiemail{#1}}
\long\def\asciiabstract#1{\long\def\theasciiabstract{#1}}

\def\shortauthors#1{\def\theshortauthors{#1}}

\let\\\par
\let\thevolumenumber\relax\let\thepapernumber\relax
\let\thevolumeyear\relax\let\thesamplenumber\relax\let\startpage\relax
\let\finishpage\relax\let\publishdate\relax\let\receiveddate\relax
\let\reviseddate\relax\let\accepteddate\relax\let\theasciititle\relax
\let\theasciiauthors\relax\let\theasciiaddress\relax
\let\theasciiabstract\relax
\let\theasciiemail\relax\let\theshortauthors\relax\let\theshorttitle\relax

\long\def\maketitlep{   

\count0=\startpage

\gt\hfill      
\beginpicture
\setcoordinatesystem units <0.33truein, 0.33truein> point at 2.2 0.9
\setplotsymbol ({$\cal G$})
\plotsymbolspacing=9truept
\circulararc 315 degrees from 0 1 center at 0 0
\setplotsymbol ({$\cal T$})
\circulararc 315 degrees from 1 -1 center at 1 0
\endpicture
\break
{\small\ifx\thesamplenumber\relax 
Volume \else Sample
\fi\thevolumenumber\ (\thevolumeyear)
\startpage--\finishpage\nl
Published: \publishdate}
\vglue 0.5truein plus 0.4fil minus 0.1truein

{\parskip=0pt\leftskip 0pt plus 1fil\def\\{\par\smallskip}{\ifplaintex\large
\else\Large\fi\bf\thetitle}\par\medskip}   

\vglue 0pt plus 0.1fil 

{\parskip=0pt\leftskip 0pt plus 1fil\def\\{\par}{\sc\theauthors}
\par\medskip}

\vglue 0pt plus 0.1fil 

{\small\parskip=0pt\let\newline\\
{\leftskip 0pt plus 1fil\def\\{\par}{\sl\theaddress}\par}
\expandafter\ifx\theemail\relax    
\relax\else\vglue 5pt plus 0.02fil minus 2pt\def\\{\stdspace{\rm 
and}\stdspace} 
\cl{Email:\stdspace\tt\theemail}\fi
\ifx\theurl\relax                  
\relax\else\vglue 5pt plus 0.02fil minus 2pt\def\\{\stdspace{\rm 
and}\stdspace}
\cl{URL:\stdspace\tt\theurl}\fi\par}

\vglue 7pt plus 0.3fil minus 3pt

{\bf Abstract}
\vglue 5pt plus 0.1fil minus 2pt

\theabstract

\vglue 7pt plus 0.3fil minus 3pt

{\bf AMS Classification numbers}\quad Primary:\quad \theprimaryclass

Secondary:\quad \thesecondaryclass

\vglue 5pt plus 0.3fil minus 2pt

{\bf Keywords}\quad \thekeywords

\vglue 10pt plus 0.5fil minus 5pt

{\small  Proposed: \theproposer\hfill Received: \receiveddate\nl
Seconded: \theseconders\hfill 
\ifx\reviseddate\relax                         
Accepted: \accepteddate                        
\else
Revised: \reviseddate                          
\fi}
\eject
}       

\let\maketitlepage\maketitlep
\let\maketitle\maketitlepage

\font\phead=cmsl9 scaled 950
\font=cmsl9 scaled 1050
\font\pnum=cmbx10 scaled 913
\font\lnum=cmbx10 
\font\pfoot=cmsl9 scaled 950
\font=cmsl9 scaled 1050
\ifplaintex
\headline{\vbox to 0pt{\vskip -4.5mm\line{\small\phead\ifnum
\count0=\startpage ISSN 1364-0380 (on line)
1465-3060 (printed) \hfill {\pnum\folio}\else\ifodd\count0\def\\{ }%
\ifx\theshorttitle\relax\thetitle\else\theshorttitle\fi\hfill{\pnum\folio}
\else\def\\{ and }{\pnum\folio}\hfill\ifx\theshortauthors\relax\theauthors
\else\theshortauthors\fi\fi\fi}\vss}}
\footline{\vbox to 0pt{\vglue 0mm\line{\small\pfoot\ifnum\count0=\startpage
\copyright\ \gtp\hfill\else
\gt, Volume \thevolumenumber\ (\thevolumeyear)\hfill\fi}\vss
}}
\else
\makeatletter
\def\@oddhead{{\small\lhead\ifnum\count0=\startpage ISSN 1364-0380 (on line)
1465-3060 (printed) \hfill {\lnum\number\count0}\else\ifodd\count0
\def\\{ }\ifx\theshorttitle\relax \thetitle \else\theshorttitle\fi\hfill
{\lnum\number\count0}\else\def\\{ and }{\lnum\number\count0}
\hfill\ifx\theshortauthors\relax 
\theauthors\else\theshortauthors\fi\fi\fi}}\def\@evenhead{\@oddhead}
\def\@oddfoot{\small\lfoot\ifnum\count0=\startpage\copyright\ \gtp\hfill\else
\gt, Volume \thevolumenumber\ (\thevolumeyear)\hfill\fi}
\def\@evenfoot{\@oddfoot}
\makeatother
\fi

\newwrite\gtoutfile
\long\gdef\makeheadfile{  
{\def\\{, }\def\s{ }
\immediate\openout\gtoutfile head.xxx
\immediate\write\gtoutfile{To: math@arxiv.org}
\immediate\write\gtoutfile{Subject: put or rep NNNNN:pppp}
\immediate\write\gtoutfile{--text follows this line--}
\immediate\write\gtoutfile{Proxy-for: \ifx\theasciiauthors\relax
\theauthors\else\theasciiauthors\fi\s<\ifx\theasciiemail\relax\theemail\else\theasciiemail\fi>}
\immediate\write\gtoutfile{\noexpand\\}
\immediate\write\gtoutfile{Authors: \ifx\theasciiauthors\relax
\theauthors\else\theasciiauthors\fi}
\immediate\write\gtoutfile{Title: \ifx\theasciititle\relax
\thetitle\else\theasciititle\fi}
\immediate\write\gtoutfile{Subj-class: GT or SG or MG etc}
\immediate\write\gtoutfile{MSC-class: \theprimaryclass\ifx\thesecondaryclass\relax\else, \thesecondaryclass\fi}
\immediate\write\gtoutfile{Journal-ref: Geom. Topol. \thevolumenumber
(\thevolumeyear) \startpage-\finishpage}
\immediate\write\gtoutfile{Comments: Published by Geometry and Topology at}
\immediate\write\gtoutfile{\s\s http://www.maths.warwick.ac.uk/gt/GTVol\thevolumenumber/paper\thepapernumber.abs.html}
\immediate\write\gtoutfile{\noexpand\\}
\immediate\write\gtoutfile{}
\ifx\theasciiabstract\relax
\immediate\write\gtoutfile{\theabstract}\else
\immediate\write\gtoutfile{\theasciiabstract}\fi
\immediate\write\gtoutfile{}
\immediate\write\gtoutfile{\noexpand\\}
\immediate\write\gtoutfile{}
\immediate\closeout\gtoutfile}}  

\def\maketitlepage{\maketitlep\makeheadfile}
\let\maketitle\maketitlepage

\def\ifplaintex{\expandafter\ifx\csname documentclass\endcsname\relax}

\ifplaintex 
\else
\fi

\expandafter\ifx\csname beginpicture\endcsname\relax
\expandafter\ifx\csname documentclass\endcsname\relax
\input pictex \else
\input prepictex \input pictex \input postpictex \fi\fi

\def\gt{{\mathsurround=0pt\it $\cal G\mskip-2mu$eometry \&\ 
$\cal T\!\!$opology}}        

\def\gtp{{\mathsurround=0pt\it $\cal G\mskip-2mu$eometry \&\ 
$\cal T\!\!$opology $\cal P\!$ublications}}  

\def\volumenumber#1{\def\thevolumenumber{#1}}
\def\papernumber#1{\def\thepapernumber{#1}}
\def\volumeyear#1{\def\thevolumeyear{#1}}

\def\pagenumbers#1#2{\def\startpage{#1}\def\finishpage{#2}}
\def\published#1{\def\publishdate{#1}}
\def\proposed#1{\def\theproposer{#1}}
\def\seconded#1{\def\theseconders{#1}}
\def\received#1{\def\receiveddate{#1}}

\def\accepted#1{\def\accepteddate{#1}}

\def\asciiaddress#1{\def\theasciiaddress{#1}}
\def\asciiemail#1{\def\theasciiemail{#1}}
\long\def\asciiabstract#1{\long\def\theasciiabstract{#1}}

\def\shortauthors#1{\def\theshortauthors{#1}}

\let\\\par
\let\thevolumenumber\relax\let\thepapernumber\relax
\let\thevolumeyear\relax\let\thesamplenumber\relax\let\startpage\relax
\let\finishpage\relax\let\publishdate\relax\let\receiveddate\relax
\let\reviseddate\relax\let\accepteddate\relax\let\theasciititle\relax
\let\theasciiauthors\relax\let\theasciiaddress\relax
\let\theasciiabstract\relax
\let\theasciiemail\relax\let\theshortauthors\relax\let\theshorttitle\relax

\long\def\maketitlep{   

\count0=\startpage

\gt\hfill      
\beginpicture
\setcoordinatesystem units <0.33truein, 0.33truein> point at 2.2 0.9
\setplotsymbol ({$\cal G$})
\plotsymbolspacing=9truept
\circulararc 315 degrees from 0 1 center at 0 0
\setplotsymbol ({$\cal T$})
\circulararc 315 degrees from 1 -1 center at 1 0
\endpicture
\break
{\small\ifx\thesamplenumber\relax 
Volume \else Sample
\fi\thevolumenumber\ (\thevolumeyear)
\startpage--\finishpage\nl
Published: \publishdate}
\vglue 0.5truein plus 0.4fil minus 0.1truein

{\parskip=0pt\leftskip 0pt plus 1fil\def\\{\par\smallskip}{\ifplaintex\large
\else\Large\fi\bf\thetitle}\par\medskip}   

\vglue 0pt plus 0.1fil 

{\parskip=0pt\leftskip 0pt plus 1fil\def\\{\par}{\sc\theauthors}
\par\medskip}

\vglue 0pt plus 0.1fil 

{\small\parskip=0pt\let\newline\\
{\leftskip 0pt plus 1fil\def\\{\par}{\sl\theaddress}\par}
\expandafter\ifx\theemail\relax    
\relax\else\vglue 5pt plus 0.02fil minus 2pt\def\\{\stdspace{\rm 
and}\stdspace} 
\cl{Email:\stdspace\tt\theemail}\fi
\ifx\theurl\relax                  
\relax\else\vglue 5pt plus 0.02fil minus 2pt\def\\{\stdspace{\rm 
and}\stdspace}
\cl{URL:\stdspace\tt\theurl}\fi\par}

\vglue 7pt plus 0.3fil minus 3pt

{\bf Abstract}
\vglue 5pt plus 0.1fil minus 2pt

\theabstract

\vglue 7pt plus 0.3fil minus 3pt

{\bf AMS Classification numbers}\quad Primary:\quad \theprimaryclass

Secondary:\quad \thesecondaryclass

\vglue 5pt plus 0.3fil minus 2pt

{\bf Keywords}\quad \thekeywords

\vglue 10pt plus 0.5fil minus 5pt

{\small  Proposed: \theproposer\hfill Received: \receiveddate\nl
Seconded: \theseconders\hfill 
\ifx\reviseddate\relax                         
Accepted: \accepteddate                        
\else
Revised: \reviseddate                          
\fi}
\eject
}       

\let\maketitlepage\maketitlep
\let\maketitle\maketitlepage

\font\phead=cmsl9 scaled 950
\font=cmsl9 scaled 1050
\font\pnum=cmbx10 scaled 913
\font\lnum=cmbx10 
\font\pfoot=cmsl9 scaled 950
\font=cmsl9 scaled 1050
\ifplaintex
\headline{\vbox to 0pt{\vskip -4.5mm\line{\small\phead\ifnum
\count0=\startpage ISSN 1364-0380 (on line)
1465-3060 (printed) \hfill {\pnum\folio}\else\ifodd\count0\def\\{ }%
\ifx\theshorttitle\relax\thetitle\else\theshorttitle\fi\hfill{\pnum\folio}
\else\def\\{ and }{\pnum\folio}\hfill\ifx\theshortauthors\relax\theauthors
\else\theshortauthors\fi\fi\fi}\vss}}
\footline{\vbox to 0pt{\vglue 0mm\line{\small\pfoot\ifnum\count0=\startpage
\copyright\ \gtp\hfill\else
\gt, Volume \thevolumenumber\ (\thevolumeyear)\hfill\fi}\vss
}}
\else
\makeatletter
\def\@oddhead{{\small\lhead\ifnum\count0=\startpage ISSN 1364-0380 (on line)
1465-3060 (printed) \hfill {\lnum\number\count0}\else\ifodd\count0
\def\\{ }\ifx\theshorttitle\relax \thetitle \else\theshorttitle\fi\hfill
{\lnum\number\count0}\else\def\\{ and }{\lnum\number\count0}
\hfill\ifx\theshortauthors\relax 
\theauthors\else\theshortauthors\fi\fi\fi}}\def\@evenhead{\@oddhead}
\def\@oddfoot{\small\lfoot\ifnum\count0=\startpage\copyright\ \gtp\hfill\else
\gt, Volume \thevolumenumber\ (\thevolumeyear)\hfill\fi}
\def\@evenfoot{\@oddfoot}
\makeatother
\fi

\newwrite\gtoutfile
\long\gdef\makeheadfile{  
{\def\\{, }\def\s{ }
\immediate\openout\gtoutfile head.xxx
\immediate\write\gtoutfile{To: math@arxiv.org}
\immediate\write\gtoutfile{Subject: put or rep NNNNN:pppp}
\immediate\write\gtoutfile{--text follows this line--}
\immediate\write\gtoutfile{Proxy-for: \ifx\theasciiauthors\relax
\theauthors\else\theasciiauthors\fi\s<\ifx\theasciiemail\relax\theemail\else\theasciiemail\fi>}
\immediate\write\gtoutfile{\noexpand\\}
\immediate\write\gtoutfile{Authors: \ifx\theasciiauthors\relax
\theauthors\else\theasciiauthors\fi}
\immediate\write\gtoutfile{Title: \ifx\theasciititle\relax
\thetitle\else\theasciititle\fi}
\immediate\write\gtoutfile{Subj-class: GT or SG or MG etc}
\immediate\write\gtoutfile{MSC-class: \theprimaryclass\ifx\thesecondaryclass\relax\else, \thesecondaryclass\fi}
\immediate\write\gtoutfile{Journal-ref: Geom. Topol. \thevolumenumber
(\thevolumeyear) \startpage-\finishpage}
\immediate\write\gtoutfile{Comments: Published by Geometry and Topology at}
\immediate\write\gtoutfile{\s\s http://www.maths.warwick.ac.uk/gt/GTVol\thevolumenumber/paper\thepapernumber.abs.html}
\immediate\write\gtoutfile{\noexpand\\}
\immediate\write\gtoutfile{}
\ifx\theasciiabstract\relax
\immediate\write\gtoutfile{\theabstract}\else
\immediate\write\gtoutfile{\theasciiabstract}\fi
\immediate\write\gtoutfile{}
\immediate\write\gtoutfile{\noexpand\\}
\immediate\write\gtoutfile{}
\immediate\closeout\gtoutfile}}  

\def\maketitlepage{\maketitlep\makeheadfile}
\let\maketitle\maketitlepage

\def\ifplaintex{\expandafter\ifx\csname documentclass\endcsname\relax}

\ifplaintex 
\else
\fi

\expandafter\ifx\csname beginpicture\endcsname\relax
\expandafter\ifx\csname documentclass\endcsname\relax
\input pictex \else
\input prepictex \input pictex \input postpictex \fi\fi

\def\gt{{\mathsurround=0pt\it $\cal G\mskip-2mu$eometry \&\ 
$\cal T\!\!$opology}}        

\def\gtp{{\mathsurround=0pt\it $\cal G\mskip-2mu$eometry \&\ 
$\cal T\!\!$opology $\cal P\!$ublications}}  

\def\volumenumber#1{\def\thevolumenumber{#1}}
\def\papernumber#1{\def\thepapernumber{#1}}
\def\volumeyear#1{\def\thevolumeyear{#1}}

\def\pagenumbers#1#2{\def\startpage{#1}\def\finishpage{#2}}
\def\published#1{\def\publishdate{#1}}
\def\proposed#1{\def\theproposer{#1}}
\def\seconded#1{\def\theseconders{#1}}
\def\received#1{\def\receiveddate{#1}}

\def\accepted#1{\def\accepteddate{#1}}

\def\asciiaddress#1{\def\theasciiaddress{#1}}
\def\asciiemail#1{\def\theasciiemail{#1}}
\long\def\asciiabstract#1{\long\def\theasciiabstract{#1}}

\def\shortauthors#1{\def\theshortauthors{#1}}

\let\\\par
\let\thevolumenumber\relax\let\thepapernumber\relax
\let\thevolumeyear\relax\let\thesamplenumber\relax\let\startpage\relax
\let\finishpage\relax\let\publishdate\relax\let\receiveddate\relax
\let\reviseddate\relax\let\accepteddate\relax\let\theasciititle\relax
\let\theasciiauthors\relax\let\theasciiaddress\relax
\let\theasciiabstract\relax
\let\theasciiemail\relax\let\theshortauthors\relax\let\theshorttitle\relax

\long\def\maketitlep{   

\count0=\startpage

\gt\hfill      
\beginpicture
\setcoordinatesystem units <0.33truein, 0.33truein> point at 2.2 0.9
\setplotsymbol ({$\cal G$})
\plotsymbolspacing=9truept
\circulararc 315 degrees from 0 1 center at 0 0
\setplotsymbol ({$\cal T$})
\circulararc 315 degrees from 1 -1 center at 1 0
\endpicture
\break
{\small\ifx\thesamplenumber\relax 
Volume \else Sample
\fi\thevolumenumber\ (\thevolumeyear)
\startpage--\finishpage\nl
Published: \publishdate}
\vglue 0.5truein plus 0.4fil minus 0.1truein

{\parskip=0pt\leftskip 0pt plus 1fil\def\\{\par\smallskip}{\ifplaintex\large
\else\Large\fi\bf\thetitle}\par\medskip}   

\vglue 0pt plus 0.1fil 

{\parskip=0pt\leftskip 0pt plus 1fil\def\\{\par}{\sc\theauthors}
\par\medskip}

\vglue 0pt plus 0.1fil 

{\small\parskip=0pt\let\newline\\
{\leftskip 0pt plus 1fil\def\\{\par}{\sl\theaddress}\par}
\expandafter\ifx\theemail\relax    
\relax\else\vglue 5pt plus 0.02fil minus 2pt\def\\{\stdspace{\rm 
and}\stdspace} 
\cl{Email:\stdspace\tt\theemail}\fi
\ifx\theurl\relax                  
\relax\else\vglue 5pt plus 0.02fil minus 2pt\def\\{\stdspace{\rm 
and}\stdspace}
\cl{URL:\stdspace\tt\theurl}\fi\par}

\vglue 7pt plus 0.3fil minus 3pt

{\bf Abstract}
\vglue 5pt plus 0.1fil minus 2pt

\theabstract

\vglue 7pt plus 0.3fil minus 3pt

{\bf AMS Classification numbers}\quad Primary:\quad \theprimaryclass

Secondary:\quad \thesecondaryclass

\vglue 5pt plus 0.3fil minus 2pt

{\bf Keywords}\quad \thekeywords

\vglue 10pt plus 0.5fil minus 5pt

{\small  Proposed: \theproposer\hfill Received: \receiveddate\nl
Seconded: \theseconders\hfill 
\ifx\reviseddate\relax                         
Accepted: \accepteddate                        
\else
Revised: \reviseddate                          
\fi}
\eject
}       

\let\maketitlepage\maketitlep
\let\maketitle\maketitlepage

\font\phead=cmsl9 scaled 950
\font=cmsl9 scaled 1050
\font\pnum=cmbx10 scaled 913
\font\lnum=cmbx10 
\font\pfoot=cmsl9 scaled 950
\font=cmsl9 scaled 1050
\ifplaintex
\headline{\vbox to 0pt{\vskip -4.5mm\line{\small\phead\ifnum
\count0=\startpage ISSN 1364-0380 (on line)
1465-3060 (printed) \hfill {\pnum\folio}\else\ifodd\count0\def\\{ }%
\ifx\theshorttitle\relax\thetitle\else\theshorttitle\fi\hfill{\pnum\folio}
\else\def\\{ and }{\pnum\folio}\hfill\ifx\theshortauthors\relax\theauthors
\else\theshortauthors\fi\fi\fi}\vss}}
\footline{\vbox to 0pt{\vglue 0mm\line{\small\pfoot\ifnum\count0=\startpage
\copyright\ \gtp\hfill\else
\gt, Volume \thevolumenumber\ (\thevolumeyear)\hfill\fi}\vss
}}
\else
\makeatletter
\def\@oddhead{{\small\lhead\ifnum\count0=\startpage ISSN 1364-0380 (on line)
1465-3060 (printed) \hfill {\lnum\number\count0}\else\ifodd\count0
\def\\{ }\ifx\theshorttitle\relax \thetitle \else\theshorttitle\fi\hfill
{\lnum\number\count0}\else\def\\{ and }{\lnum\number\count0}
\hfill\ifx\theshortauthors\relax 
\theauthors\else\theshortauthors\fi\fi\fi}}\def\@evenhead{\@oddhead}
\def\@oddfoot{\small\lfoot\ifnum\count0=\startpage\copyright\ \gtp\hfill\else
\gt, Volume \thevolumenumber\ (\thevolumeyear)\hfill\fi}
\def\@evenfoot{\@oddfoot}
\makeatother
\fi

\newwrite\gtoutfile
\long\gdef\makeheadfile{  
{\def\\{, }\def\s{ }
\immediate\openout\gtoutfile head.xxx
\immediate\write\gtoutfile{To: math@arxiv.org}
\immediate\write\gtoutfile{Subject: put or rep NNNNN:pppp}
\immediate\write\gtoutfile{--text follows this line--}
\immediate\write\gtoutfile{Proxy-for: \ifx\theasciiauthors\relax
\theauthors\else\theasciiauthors\fi\s<\ifx\theasciiemail\relax\theemail\else\theasciiemail\fi>}
\immediate\write\gtoutfile{\noexpand\\}
\immediate\write\gtoutfile{Authors: \ifx\theasciiauthors\relax
\theauthors\else\theasciiauthors\fi}
\immediate\write\gtoutfile{Title: \ifx\theasciititle\relax
\thetitle\else\theasciititle\fi}
\immediate\write\gtoutfile{Subj-class: GT or SG or MG etc}
\immediate\write\gtoutfile{MSC-class: \theprimaryclass\ifx\thesecondaryclass\relax\else, \thesecondaryclass\fi}
\immediate\write\gtoutfile{Journal-ref: Geom. Topol. \thevolumenumber
(\thevolumeyear) \startpage-\finishpage}
\immediate\write\gtoutfile{Comments: Published by Geometry and Topology at}
\immediate\write\gtoutfile{\s\s http://www.maths.warwick.ac.uk/gt/GTVol\thevolumenumber/paper\thepapernumber.abs.html}
\immediate\write\gtoutfile{\noexpand\\}
\immediate\write\gtoutfile{}
\ifx\theasciiabstract\relax
\immediate\write\gtoutfile{\theabstract}\else
\immediate\write\gtoutfile{\theasciiabstract}\fi
\immediate\write\gtoutfile{}
\immediate\write\gtoutfile{\noexpand\\}
\immediate\write\gtoutfile{}
\immediate\closeout\gtoutfile}}  

\def\maketitlepage{\maketitlep\makeheadfile}
\let\maketitle\maketitlepage

\def\ifplaintex{\expandafter\ifx\csname documentclass\endcsname\relax}

\ifplaintex 
\else
\fi

\expandafter\ifx\csname beginpicture\endcsname\relax
\expandafter\ifx\csname documentclass\endcsname\relax
\input pictex \else
\input prepictex \input pictex \input postpictex \fi\fi

\def\gt{{\mathsurround=0pt\it $\cal G\mskip-2mu$eometry \&\ 
$\cal T\!\!$opology}}        

\def\gtp{{\mathsurround=0pt\it $\cal G\mskip-2mu$eometry \&\ 
$\cal T\!\!$opology $\cal P\!$ublications}}  

\def\volumenumber#1{\def\thevolumenumber{#1}}
\def\papernumber#1{\def\thepapernumber{#1}}
\def\volumeyear#1{\def\thevolumeyear{#1}}

\def\pagenumbers#1#2{\def\startpage{#1}\def\finishpage{#2}}
\def\published#1{\def\publishdate{#1}}
\def\proposed#1{\def\theproposer{#1}}
\def\seconded#1{\def\theseconders{#1}}
\def\received#1{\def\receiveddate{#1}}

\def\accepted#1{\def\accepteddate{#1}}

\def\asciiaddress#1{\def\theasciiaddress{#1}}
\def\asciiemail#1{\def\theasciiemail{#1}}
\long\def\asciiabstract#1{\long\def\theasciiabstract{#1}}

\def\shortauthors#1{\def\theshortauthors{#1}}

\let\\\par
\let\thevolumenumber\relax\let\thepapernumber\relax
\let\thevolumeyear\relax\let\thesamplenumber\relax\let\startpage\relax
\let\finishpage\relax\let\publishdate\relax\let\receiveddate\relax
\let\reviseddate\relax\let\accepteddate\relax\let\theasciititle\relax
\let\theasciiauthors\relax\let\theasciiaddress\relax
\let\theasciiabstract\relax
\let\theasciiemail\relax\let\theshortauthors\relax\let\theshorttitle\relax

\long\def\maketitlep{   

\count0=\startpage

\gt\hfill      
\beginpicture
\setcoordinatesystem units <0.33truein, 0.33truein> point at 2.2 0.9
\setplotsymbol ({$\cal G$})
\plotsymbolspacing=9truept
\circulararc 315 degrees from 0 1 center at 0 0
\setplotsymbol ({$\cal T$})
\circulararc 315 degrees from 1 -1 center at 1 0
\endpicture
\break
{\small\ifx\thesamplenumber\relax 
Volume \else Sample
\fi\thevolumenumber\ (\thevolumeyear)
\startpage--\finishpage\nl
Published: \publishdate}
\vglue 0.5truein plus 0.4fil minus 0.1truein

{\parskip=0pt\leftskip 0pt plus 1fil\def\\{\par\smallskip}{\ifplaintex\large
\else\Large\fi\bf\thetitle}\par\medskip}   

\vglue 0pt plus 0.1fil 

{\parskip=0pt\leftskip 0pt plus 1fil\def\\{\par}{\sc\theauthors}
\par\medskip}

\vglue 0pt plus 0.1fil 

{\small\parskip=0pt\let\newline\\
{\leftskip 0pt plus 1fil\def\\{\par}{\sl\theaddress}\par}
\expandafter\ifx\theemail\relax    
\relax\else\vglue 5pt plus 0.02fil minus 2pt\def\\{\stdspace{\rm 
and}\stdspace} 
\cl{Email:\stdspace\tt\theemail}\fi
\ifx\theurl\relax                  
\relax\else\vglue 5pt plus 0.02fil minus 2pt\def\\{\stdspace{\rm 
and}\stdspace}
\cl{URL:\stdspace\tt\theurl}\fi\par}

\vglue 7pt plus 0.3fil minus 3pt

{\bf Abstract}
\vglue 5pt plus 0.1fil minus 2pt

\theabstract

\vglue 7pt plus 0.3fil minus 3pt

{\bf AMS Classification numbers}\quad Primary:\quad \theprimaryclass

Secondary:\quad \thesecondaryclass

\vglue 5pt plus 0.3fil minus 2pt

{\bf Keywords}\quad \thekeywords

\vglue 10pt plus 0.5fil minus 5pt

{\small  Proposed: \theproposer\hfill Received: \receiveddate\nl
Seconded: \theseconders\hfill 
\ifx\reviseddate\relax                         
Accepted: \accepteddate                        
\else
Revised: \reviseddate                          
\fi}
\eject
}       

\let\maketitlepage\maketitlep
\let\maketitle\maketitlepage

\font\phead=cmsl9 scaled 950
\font=cmsl9 scaled 1050
\font\pnum=cmbx10 scaled 913
\font\lnum=cmbx10 
\font\pfoot=cmsl9 scaled 950
\font=cmsl9 scaled 1050
\ifplaintex
\headline{\vbox to 0pt{\vskip -4.5mm\line{\small\phead\ifnum
\count0=\startpage ISSN 1364-0380 (on line)
1465-3060 (printed) \hfill {\pnum\folio}\else\ifodd\count0\def\\{ }%
\ifx\theshorttitle\relax\thetitle\else\theshorttitle\fi\hfill{\pnum\folio}
\else\def\\{ and }{\pnum\folio}\hfill\ifx\theshortauthors\relax\theauthors
\else\theshortauthors\fi\fi\fi}\vss}}
\footline{\vbox to 0pt{\vglue 0mm\line{\small\pfoot\ifnum\count0=\startpage
\copyright\ \gtp\hfill\else
\gt, Volume \thevolumenumber\ (\thevolumeyear)\hfill\fi}\vss
}}
\else
\makeatletter
\def\@oddhead{{\small\lhead\ifnum\count0=\startpage ISSN 1364-0380 (on line)
1465-3060 (printed) \hfill {\lnum\number\count0}\else\ifodd\count0
\def\\{ }\ifx\theshorttitle\relax \thetitle \else\theshorttitle\fi\hfill
{\lnum\number\count0}\else\def\\{ and }{\lnum\number\count0}
\hfill\ifx\theshortauthors\relax 
\theauthors\else\theshortauthors\fi\fi\fi}}\def\@evenhead{\@oddhead}
\def\@oddfoot{\small\lfoot\ifnum\count0=\startpage\copyright\ \gtp\hfill\else
\gt, Volume \thevolumenumber\ (\thevolumeyear)\hfill\fi}
\def\@evenfoot{\@oddfoot}
\makeatother
\fi

\newwrite\gtoutfile
\long\gdef\makeheadfile{  
{\def\\{, }\def\s{ }
\immediate\openout\gtoutfile head.xxx
\immediate\write\gtoutfile{To: math@arxiv.org}
\immediate\write\gtoutfile{Subject: put or rep NNNNN:pppp}
\immediate\write\gtoutfile{--text follows this line--}
\immediate\write\gtoutfile{Proxy-for: \ifx\theasciiauthors\relax
\theauthors\else\theasciiauthors\fi\s<\ifx\theasciiemail\relax\theemail\else\theasciiemail\fi>}
\immediate\write\gtoutfile{\noexpand\\}
\immediate\write\gtoutfile{Authors: \ifx\theasciiauthors\relax
\theauthors\else\theasciiauthors\fi}
\immediate\write\gtoutfile{Title: \ifx\theasciititle\relax
\thetitle\else\theasciititle\fi}
\immediate\write\gtoutfile{Subj-class: GT or SG or MG etc}
\immediate\write\gtoutfile{MSC-class: \theprimaryclass\ifx\thesecondaryclass\relax\else, \thesecondaryclass\fi}
\immediate\write\gtoutfile{Journal-ref: Geom. Topol. \thevolumenumber
(\thevolumeyear) \startpage-\finishpage}
\immediate\write\gtoutfile{Comments: Published by Geometry and Topology at}
\immediate\write\gtoutfile{\s\s http://www.maths.warwick.ac.uk/gt/GTVol\thevolumenumber/paper\thepapernumber.abs.html}
\immediate\write\gtoutfile{\noexpand\\}
\immediate\write\gtoutfile{}
\ifx\theasciiabstract\relax
\immediate\write\gtoutfile{\theabstract}\else
\immediate\write\gtoutfile{\theasciiabstract}\fi
\immediate\write\gtoutfile{}
\immediate\write\gtoutfile{\noexpand\\}
\immediate\write\gtoutfile{}
\immediate\closeout\gtoutfile}}  

\def\maketitlepage{\maketitlep\makeheadfile}
\let\maketitle\maketitlepage

\volumenumber{5}\papernumber{6}\volumeyear{2001}
\pagenumbers{143}{226}
\accepted{7 March 2001}
\proposed{Tomasz Mrowka}
\seconded{Ronald Fintushel, Ronald Stern}
\received{20 September 1999}
\published{21 March 2001}

\usepackage{amscd, amsmath, amssymb}
\usepackage{graphicx}

\numberwithin{equation}{section}

\begin{document}

\def\CC{{\mathbb C}}
\def\PP{{\mathbb CP}}
\def\QQ{{\mathbb Q}}
\def\RR{{\mathbb R}}
\def\ZZ{{\mathbb Z}}

\def\cA{{\mathcal A}}
\def\cB{{\mathcal B}}
\def\cG{{\mathcal G}}
\def\cH{{\mathcal H}}
\def\cL{{\mathcal L}}
\def\fM{{\mathfrak M}}
\def\fR{{\mathfrak R}}

\def\al{\alpha}
\def\be{\beta}
\def\ga{\gamma}
\def\de{\delta}
\def\ep{\varepsilon}
\def\ka{\kappa}
\def\la{\lambda}
\def\om{\omega}
\def\si{\sigma}
\def\th{\theta}
\def\Th{\Theta}
\def\Ga{\Gamma}
\def\La{\Lambda}
\def\Om{\Omega}
\def\Si{\Sigma}

\def\tA{\widetilde A}
\def\tB{\widetilde B}
\def\tC{\widetilde C}
\def\tit{\tilde t}
\def\tJ{\tilde J}
\def\tV{\widetilde V}
\def\tP{\widetilde P}
\def\tR{\widetilde R}
\def\hR{\widehat R}
\def\hC{\widehat C}
\def\tm{\tilde m}
\def\tn{\tilde n}
\def\tmu{\tilde\mu}
\def\tla{\tilde\la}
\def\tal{\tilde\alpha}
\def\tbe{\tilde\beta}
\def\tOm{\widetilde \Omega}

\def\hA{\widehat A}

\def\lk{{\rm lk}}
\def\lto{\longrightarrow}
\newcommand{\codim}{\operatorname{codim}}
\newcommand{\im}{\operatorname{im}}
\newcommand{\hol}{\operatorname{hol}}
\newcommand{\sspan}{\operatorname{span}}
\newcommand{\Image}{\operatorname{Im}}
\newcommand{\Mas}{\operatorname{Mas}}
\newcommand{\Lag}{\operatorname{Lag}}
\newcommand{\Id}{\operatorname{Id}}
\newcommand{\proj}{\operatorname{proj}}
\newcommand{\limit}{\operatornamewithlimits{limit}}
\newcommand{\vol}{\operatorname{vol}}
\newcommand{\id}{\operatorname{id}}
\def\ccspan{\sspan_{{\mathbb C}^2}}
\def\contract{\lrcorner}
\def\Hom{\hbox{Hom}}

\def\mapright#1{\smash{ \mathop{\longrightarrow}\limits^{#1}}}

\newtheorem{theorem}{Theorem}[section]
\newtheorem{lemma}[theorem]{Lemma}
\newtheorem{corollary}[theorem]{Corollary}
\newtheorem{proposition}[theorem]{Proposition}

\theoremstyle{definition}
\newtheorem{definition}[theorem]{Definition}

\theoremstyle{remark} \newtheorem{remark}[theorem]{Remark}
\newtheorem{example}[theorem]{Example}
\newtheorem{claim}[theorem]{Claim}

\title[Gauge Theoretic Invariants]{Gauge Theoretic Invariants of\\Dehn 
Surgeries on Knots}
\authors{Hans U Boden\\Christopher M Herald\\Paul A Kirk\\Eric P Klassen}
\shortauthors{Boden, Herald, Kirk, and Klassen}
\address{
McMaster University, Hamilton, Ontario L8S 4K1, Canada \\   
University of Nevada, Reno, Nevada 89557, USA \\
Indiana University, Bloomington, Indiana 47405, USA  \\ 
Florida State University, Tallahassee, Florida 32306, USA  \\ \medskip
\rm Email addresses: \stdspace \tt boden@math.mcmaster.ca\ \ herald@unr.edu \\
pkirk@indiana.edu\ \ klassen@zeno.math.fsu.edu}

\asciiaddress{McMaster University, Hamilton, Ontario L8S 4K1, Canada\\   
University of Nevada, Reno, Nevada 89557, USA\\
Indiana University, Bloomington, Indiana 47405, USA \\ 
Florida State University, Tallahassee, Florida 32306, USA}

\asciiemail{boden@math.mcmaster.ca, herald@unr.edu,
pkirk@indiana.edu, klassen@zeno.math.fsu.edu}

\asciiabstract{New methods for computing a variety of gauge theoretic
invariants for homology 3-spheres are developed. These invariants
include the Chern-Simons invariants, the spectral flow of the odd
signature operator, and the rho invariants of irreducible SU(2)
representations. These quantities are calculated for flat SU(2)
connections on homology 3-spheres obtained by 1/k Dehn surgery on
(2,q) torus knots.  The methods are then applied to compute the SU(3)
gauge theoretic Casson invariant (introduced in [H U Boden and C M
Herald, The SU(3) Casson invariant for integral homology 3--spheres,
J. Diff. Geom.  50 (1998) 147-206]) for Dehn surgeries on (2,q)
torus knots for q=3,5,7 and 9. }

\begin{abstract}
New methods for computing a variety of gauge theoretic invariants for homology 
3--spheres   are developed. These
invariants include the Chern--Simons invariants, the spectral flow of 
the odd signature operator, and the rho invariants of irreducible $SU(2)$
representations. These 
quantities are calculated for flat $SU(2)$
connections on homology 3--spheres
 obtained by $1/k$ Dehn surgery on $(2,q)$ torus knots. 
The methods are then applied to compute the $SU(3)$ gauge theoretic Casson
invariant (introduced in \cite{bh}) for  Dehn surgeries 
on $(2,q)$ torus knots for $q=3,5,7$ and $9$. 
\end{abstract}
\primaryclass{57M27}
\secondaryclass{ 53D12, 58J28, 58J30}
\keywords{Homology 3--sphere, gauge theory,  
3--manifold invariants, spectral flow, Maslov index }
\date{July 28, 1999} \maketitlepage
 
\section{Introduction}
The goal of this article is to develop new methods for computing
a variety of gauge theoretic
invariants for 3--manifolds obtained by Dehn surgery on knots.
These invariants include the Chern--Simons invariants,
the spectral flow of the odd signature operator,
 and
the rho invariants of irreducible  $SU(2)$ representations.  The rho
invariants and
spectral flow considered here  are different from
the ones usually   studied in $SU(2)$ gauge theory in that they
 do not come from the adjoint representation on $su(2)$
but rather from   the canonical representation  on $\CC^2$.
Their values   are necessary to compute the $SU(3)$ Casson
invariant $\la_{SU(3)}$  defined in \cite{bh}. The
methods developed here are used together with results from \cite{boden}
to calculate $\la_{SU(3)}$  for a number of examples.

Gathering data on
    the $SU(3)$ Casson invariant is important for several reasons. First, in a
    broad sense it is unclear whether $SU(n)$ gauge theory for
    $n>2$ contains more information than can be obtained by studying
    only $SU(2)$ gauge theory. Second, as more and more combinatorially
    defined 3--manifold invariants  have recently emerged, the task of
    interpreting these new invariants in geometrically meaningful
    ways has become ever more important.  In particular, one would like to
    know whether or not $\la_{SU(3)}$ is of
    finite type. Our calculations here show that
   $\la_{SU(3)}$ not a finite type invariant
    (see Theorem \ref{notfinitetype}).

The behavior of the finite type invariants under Dehn surgery is well
    understood (in some sense it is built into their definition),
    but their relationship to  the fundamental group is not so clear.
    For example, it is unknown whether the invariants vanish on homotopy
spheres. The situation with the  $SU(3)$ Casson invariant is the
    complete opposite. It is obvious from the definition
    that $\la_{SU(3)}$ vanishes
    on  
    homotopy spheres, but its behavior under Dehn
    surgery is  subtle and not well understood.

In order to better explain the results in this paper,
we briefly recall the definition of the 
 $SU(3)$ Casson invariant $\la_{SU(3)}(X)$ for integral homology
spheres
$X$. It is given as
the sum of two terms. The first is
a signed count of the conjugacy classes of  irreducible  $SU(3)$
representations, and the second, which is called the correction
term, involves only  conjugacy classes of
irreducible  $SU(2)$ representations.
 
 To understand the need for a correction term, recall Walker's
 extension of the Casson invariant to rational homology spheres
 \cite{walker}.  Casson's invariant for integral homology spheres
 counts (with sign) the number of irreducible $SU(2)$ representations
 of $\pi_{1}X$ modulo conjugation. Prior to the count, a perturbation
 may be required to achieve transversality, but the assumption that
 $H_{1}(X;\ZZ)=0$ guarantees that the end result is independent of the
 choice of perturbation. The problem for rational homology spheres
 is that the signed count of
    irreducible $SU(2)$ representations   depends in a subtle way
    on the perturbation. To compensate,
    Walker  defined a correction term  using
  integral symplectic invariants of the reducible (ie, abelian)
    representations. This correction term can alternatively be viewed as
   a sum of  differences between the Maslov index and a nonintegral term
   \cite{CLM-bull} or as a sum of $U(1)$ rho invariants
    \cite{mrowka-walker}. 

In \cite{bh}, the objects of study are $\ZZ$--homology spheres, but the
    representations are taken in $SU(3)$.   As in the $SU(2)$ case
    there are no nontrivial abelian representations, but inside the
    $SU(3)$ representation variety there are those that reduce to
    $SU(2).$ This means that simply counting (with sign) the irreducible
    $SU(3)$ representations will not in general yield a well-defined
    invariant, and in \cite{bh} is a definition for the appropriate
    correction term involving a difference of the spectral flow
    and Chern--Simons invariants of  the reducible flat connections.
    In the simplest case, when the $SU(2)$ moduli space is  regular
    as a subspace of the $SU(3)$ moduli space, this quantity can be
    interpreted in terms of the rho invariants of Atiyah, Patodi
    and Singer \cite{APS} for flat $SU(2)$ connections (see Theorem
    \ref{2pcorrect}, for instance).

Neither the spectral flow nor the Chern--Simons invariants are
    gauge invariant, and as a result they are typically only computed
    up to some  indeterminacy. Our goal of calculating $\la_{SU(3)}$
    prevents us from  working modulo gauge, and this   technical point
    complicates the present work. In overcoming this obstacle, we
    establish  a Dehn surgery type formula (Theorem \ref{rhoinvts}) 
    for the  rho invariants in $\RR$ 
    (as opposed to the much simpler $\RR/\ZZ$--valued
    invariants).

The main results of this article are formulas which express the 
$\CC^2$--spectral flow (Theorem \ref{nicerformula}), 
the Chern--Simons invariants (Theorem \ref{csthm}),  and
the  rho invariants (Theorem \ref{rhoinvts}) for 3--manifolds $X$ 
obtained by Dehn surgery on a knot  in terms of simple
invariants of the curves in $\RR^2$ parameterizing the
$SU(2)$ representation variety of the knot complement.  The primary tools include
a splitting theorem for the $\CC^2$--spectral flow adapted for our purposes (Theorem
\ref{zorro}) and a detailed analysis of the spectral flow on a solid torus
(Section \ref{sec-csandrho}). 
These results are then applied to Dehn surgeries on torus knots,
culminating in the formulas of Theorem \ref{+2qfull}, Theorem \ref{-2qfull}, 
Table \ref{table3}, and Table \ref{table4}
giving the $\CC^2$--spectral flow,  the Chern--Simons invariants, the rho
invariants, and the $SU(3)$ Casson invariants for homology spheres
obtained by surgery on a $(2,q)$ torus
knot.

Theorem \ref{rhoinvts} can also be viewed as a small step in the program
    of extending the results of \cite{fl}. There, the rho invariant
    is shown to be a homotopy invariant up to path components of
    the representation space. More precisely, the difference in rho
    invariants of homotopy equivalent closed manifolds is a locally
    constant function on the representation space of their fundamental
    group. Our method of computing rho invariants differs from others
    in the literature in that it is a  cut--and--paste technique  rather
than one which relies on  flat bordisms or factoring representations
 through finite groups.

Previous surgery formulas for  computing spectral flow   require that
    the dimension of the cohomology of the boundary manifold be constant
    along the path of connections (see, eg~\cite{JDGpaper}). This
    restriction had to be eliminated in the present work since we
    need to compute the spectral flow starting at the
    trivial connection, where this assumption fails to
    hold. Our success in treating this issue promises
    to have other important applications to cut--and--paste methods
    for computing spectral flow.

The methods used in this article are delicate and draw on a number of
areas. The tools we use include the seminal work of
Atiyah--Patodi--Singer on the eta invariant and the index theorem
for manifolds with boundary \cite{APS}, analysis of
$SU(2)$ representation spaces of knot groups following
\cite{klassen-thesis}, the infinite dimensional symplectic analysis of
spectral flow from
\cite{Nicolaescu}, and the  analysis of the moduli of stable parabolic
bundles over Riemann surfaces from \cite{boden}.  We have attempted to
give an exposition which presents the material in bite-sized pieces, with
the goal of computing the gauge theoretic invariants in terms of a few
easily computed numerical invariants associated to $SU(2)$ representation
spaces of knot groups.

{\bf Acknowledgements}\qua The authors would like to thank
Stavros Garoufalidis for  his strengthening of
Theorem \ref{notfinitetype} and Ed Miller
for pointing out a mistake in an earlier version of
Proposition \ref{princ}.
HUB and PAK were partially supported by grants
from the National Science Foundation (DMS-9971578 and DMS-9971020).
CMH was partially supported by a Research Grant
from Swar\-thmore College.
HUB would also like to thank the
Mathematics Department at
Indiana University for the invitation to visit during
the Fall Semester of 1998.


\section{Preliminaries}

\subsection{Symplectic linear algebra} We define symplectic vector
    spaces and Lagrangian subspaces in the complex setting.

\begin{definition} \label{defn-herm-symp-vs} Suppose 
$(V, \langle \, \cdot \, , \cdot \, \rangle)$ 
    is a finite-dimensional complex vector space with positive
    definite Hermitian inner product.

\begin{enumerate}    \item[(i)] A {\it symplectic structure} is defined to be a
    skew-Hermitian nondegenerate form $\om\co V\times V\to \CC$ such
    that the signature of $i\omega$ is zero.
    Namely, $\om(x,y) = -\overline{\om(y,x)}$ for all $x,y \in V$ and
    $0=\om(x,\cdot \, ) \in V^* \Leftrightarrow x=0.$

  \item[(ii)] An {\it almost complex structure} is an isometry $J\co V\to
    V$ with $J^2=-\Id $ so that the signature of $iJ$ is zero.

  \item[(iii)] $J$ and $\om$ are {\it compatible} if $
    \om(x,y)=\langle x,Jy\rangle$ and $\om(Jx,Jy)=\om(x,y).$
    \item[(iv)] A subspace $L\subset V$ is {\it Lagrangian} if
    $\om(x,y)=0$ for all $x,y\in L$ and $\dim L =\frac{1}{2}\dim
    V.$  \end{enumerate} We shall refer to 
    $(V, \langle \, \cdot \, , \cdot \, \rangle, J,
    \om)$ as a {\it Hermitian symplectic space with compatible
    almost complex structure}.  \end{definition}

We use the same language for the complex Hilbert spaces
    $L^2(\Om^{0+1+2}_\Si\otimes \CC^2)$ of differential forms on a
    Riemannian surface $\Si$ with values in $\CC^2.$ The definitions in
    the infinite-dimensional setting are given below.

A Hermitian symplectic space can be obtained by complexifying a real
    symplectic space and extending the real inner product to a
    Hermitian inner product. The symplectic spaces we consider
    will essentially be of this form, except that we will usually
    tensor with $\CC^2$ instead of $\CC$.

In our main application (calculating $\CC^2$--spectral flow), the
    Hermitian symplectic spaces we consider are of the form
    $U\otimes_\RR\CC^2$ for a real symplectic vector space $U$. In
    most cases $U=H^{0+1+2}(\Si;\RR)$ with the symplectic structure
    given by the cup product. Furthermore, many of the Lagrangians
    we will encounter are of a special form; they are ``induced''
    from certain Lagrangians in $U\otimes_\RR\CC$. For the rest of
    this subsection we investigate certain algebraic properties of
    this special situation.

Suppose, then, that $(U, (\, \cdot \,  , \cdot \, ),J,\om)$ is 
a real symplectic vector
    space with compatible almost complex structure. Construct the
    complex symplectic vector space $$V= U \otimes_\RR \CC$$ with
    compatible almost complex structure as follows. Define $\om$ on
    $V$ by setting $$\om( u_1\otimes z_1,u_2\otimes z_2 ) = z_1
    \bar{z}_2 \, \om(u_1,u_2).$$  Similarly,
    define a Hermitian inner product 
    $ \langle \, \cdot \, , \cdot \, \rangle$ and a
    compatible almost complex structure $J$ by setting $$\langle
    u_1\otimes z_1,u_2\otimes z_2 \rangle=z_1\bar{z}_2(u_1,u_2)
    \quad \text{and} \quad J(u\otimes z)=(Ju)\otimes z.$$   It is a simple
matter to verify that the conditions
    of Definition \ref{defn-herm-symp-vs} hold and from this it follows
    that $(V, \langle \, \cdot \, , \cdot \, \rangle, J, \om)$ is a Hermitian symplectic
    space with compatible almost complex structure. Furthermore, $V$
    admits an involution $V \to V$ given by conjugation: $
    \overline{u\otimes z} \mapsto u\otimes\bar{z}.$

Now consider $$W=U \otimes_\RR \CC^2=V \otimes_\CC \CC^2.$$ Extending
    $\om, J$  and $\langle \, \cdot \, , \cdot \, \rangle$ to $W$ in the natural way, it
    follows that $W$ is also a Hermitian symplectic space with
    compatible almost complex structure. 
    Given a  linearly independent subset $ \{u_1, \ldots, u_n\}$
     of $U,$ then it follows  that the set $\{ u_1\otimes
    e_1, u_1 \otimes e_2, \ldots, u_n \otimes e_1, u_n \otimes e_2
    \}$ is  linearly independent in $W$, where $\{e_1, e_2\}$
    denotes the standard basis for $\CC^2.$ In later sections, it will be
    convenient to adopt the following notation:
\begin{equation} \label{ccspan}
    \ccspan \{u_1, \ldots,u_n\} := \sspan \{ u_1\otimes e_1, u_1 \otimes e_2,
    \ldots, u_n
    \otimes e_1, u_n \otimes e_2 \}.
\end{equation}

\subsection{The signature operator on a 3--manifold with boundary}
    \label{ssecsign}

Next we introduce the two first order differential operators which will
    be used throughout this paper. These depend on Riemannian metrics
    and orientation. We adopt the sign conventions for the Hodge
    star operator and the formal adjoint of the de~Rham differential
    for a $p$--form on an oriented Riemannian $n$--manifold whereby
$$*^2=(-1)^{p(n-p)},\ \ d^*=(-1)^{n(p+1)+1}*d*.$$
The Hodge star operator
    is defined by the formula $a\wedge *b=(a,b)\ dvol,$ where 
    $(\, \cdot \, , \, \cdot\, )$ denotes the inner product on forms 
    induced by the Riemannian
    metric and $dvol$ denotes the volume form, which depends on a
    choice of orientation.  To distinguish the star operator on the
    3--manifold from the one on the 2--manifold, we denote the former
    by $\star$ and the latter by $*$.

Every principal $SU(2)$ bundle over a 2 or 3--dimensional manifold is
    trivial. For that reason we work only with trivial bundles $P=X
    \times SU(2)$ and thereby identify connections with $su(2)$--valued
    1--forms in the usual way. Given a 3--manifold $Y$ with nonempty boundary
    $\Si$, we choose compatible trivializations of the principal $SU(2)$
    bundle over $Y$ and its restriction to $\Si$.  We will generally use  
    upper case letters such as $A$ for connections on the
    3--manifold and lower case letters such as $a$ for connections on the boundary
    surface.

Given an $SU(2)$ connection $A\in\Om^1_X\otimes su(2)$ and an $SU(2)$
representation $V$, we associate to $A$ the
covariant derivative $$d_A\co \Om^p_X\otimes V\to
    \Om^{p+1}_X\otimes V,~~ d_{A}=d+A.$$ The two representations
    that arise in this paper are  the canonical
    representation of $SU(2)$ on $\CC^{2}$ and the adjoint
    representation of $SU(2)$ on its Lie algebra $su(2)$.

The first operator we consider is the {\it twisted de~Rham operator
}
    $S_a$\glossary{$S_a$}  on the closed oriented  Riemannian
    2--manifold $\Si$.

\begin{definition} For an $SU(2)$ connection $a\in \Om^1_\Si\otimes
    su(2)$,  define the {\it twisted de~Rham operator } $S_a$ to be
    the elliptic first order differential operator\begin{eqnarray*}
    &S_a\co \Om^{0+1+2}_\Si\otimes \CC^2 \lto \Om^{0+1+2}_\Si \otimes\CC^2
    &\\ &S_a (\al,\be,\ga)=(*d_a\be,-*d_a\al-d_a*\ga,d_a*\be).&
\end{eqnarray*} This operator is self-adjoint with respect to
    the $L^2$ inner product on $\Om^{0+1+2}_\Si\otimes \CC^2$ given
    by the formula $$\langle (\al_1,\be_1,\ga_1),(\al_2,\be_2,\ga_2)
    \rangle= \int_\Si (\al_1\wedge*\al_2+\be_1\wedge*\be_2 +\ga_1\wedge
    *\ga_2),$$ where the notation for the Hermitian inner product
    in the fiber $\CC^2$ has been suppressed.  \end{definition}

It is convenient to introduce the almost complex structure
    $$J\co \Om^{0+1+2}_\Si\otimes \CC^2 \lto \Om^{0+1+2}_\Si\otimes \CC^2$$
    defined by 
    $$J(\al,\be, \ga)= (-*\ga,*\be,*\al).$$
    Clearly
    $J^2=-\Id$ and $J$ is an isometry of $ L^2(\Om^{0+1+2}_\Si\otimes
    \CC^2)$. To avoid confusion later, we point out that changing
    the orientation of $\Si$ does not affect the $L^2$ inner product but
    does change the sign of $J$.

With this almost complex structure, the Hilbert space
    $L^2(\Om^{0+1+2}_\Si\otimes \CC^2)$ becomes an infinite-dimensional
    Hermitian symplectic space, with symplectic form defined
    by $\om(x,y)=\langle x,Jy\rangle$. Recall (see, eg,
    \cite{Nicolaescu, Kirk-klassen-memoirs}) that a closed subspace
    $\La\subset L^2(\Om^{0+1+2}_\Si\otimes \CC^2)$ is called a
    {\it Lagrangian} if $\La$ is orthogonal  to $J\La$ and $\La+J\La=
    L^2(\Om^{0+1+2}_\Si\otimes \CC^2)$ (equivalently
$J\Lambda=\Lambda^\perp$). More generally a closed
    subspace $V$ is called {\it isotropic} if $V$ is orthogonal to
    $JV$.


The other operator we consider is the {\it odd signature operator}
    $D_A$\glossary{$D_A$}  on a compact, oriented, Riemannian 3--manifold
    $Y,$ with or without boundary.
    \begin{definition} For an $SU(2)$
    connection $A\in \Om^1_Y\otimes su(2)$, define the {\it odd
    signature operator} $D_A$ on $Y$ twisted by $A$ to be the formally
    self-adjoint first order differential operator\begin{eqnarray*}
    &D_A\co  \Om^{0+1}_Y \otimes \CC^2\lto \Om^{0+1}_Y\otimes \CC^2&
    \\ &D_A(\si,\tau)= (d_A^*\tau, d_A\si + \star d_A\tau).
\end{eqnarray*}
 \end{definition}

We wish to relate the operators $D_A$ and $S_a$ in the case when $Y$ has
    boundary $\Si$ and $a=A|_{\Si}$. The easiest way to avoid confusion
arising from
    orientation conventions is to first work on the cylinder
    $[-1,1]\times\Si$. So assume that $\Si $ is an oriented closed
    surface with Riemannian metric and that $[-1,1]\times\Si$ is
    given the product metric and the product orientation ${\mathcal
    O}_{[-1,1]\times\Si}=\{ du,{\mathcal O}_\Si\}$. Thus $\partial
    ([-1,1]\times\Si)= (\{1\} \times \Si) \cup -( \{-1\} \times
    \Si)$ using the outward normal first convention.

Assume further that $a\in \Om^{1}_{\Si}\otimes su(2)$ and
$A=\pi^{*}a\in \Om^{1}_{[-1,1]\times \Si}\otimes su(2)$,
the pullback of $a$ by the projection
    $$\pi \co [-1,1]\times \Si \to \Si.$$
    In other words,
    $$ d_A = d_a+du\wedge
    \tfrac{\partial}{\partial u} ,$$ where $u$ denotes
    the $[-1,1]$ coordinate.

Denote by $\tOm^{0+1+2}_{[-1,1]\times\Si}$ the
    space of forms on the cylinder with no $du$ component 
    and define $$\Phi\co  \Om^{0+1}_{[-1,1]\times\Si} \otimes \CC^2 \lto
    \tOm^{0+1+2}_{[-1,1]\times\Si} \otimes \CC^2$$
    $$\Phi(\si,\tau)= (i_u^*(\si)
    ,i_u^*(\tau), * i_u^*(\tau \contract \tfrac{\partial}{\partial
    u})),$$ where
    $i_{u} \co  \Si \hookrightarrow [-1,1]\times\Si$ is the
    inclusion at $u$ and $\contract$ denotes contraction.  
The following lemma is well known and follows from a straightforward
    computation.

\begin{lemma} {\sl $\Phi \circ D_A = J\circ (S_a+{\scriptstyle
    \tfrac{\partial}{\partial u}})\circ \Phi$.  }\end{lemma}

The analysis on the cylinder carries over to a general 3--manifold with boundary
$\Si$ given
    an identification of the collar of the boundary with $I\times\Si$.
    In the terminology of Nicolaescu's article \cite{Nicolaescu},
     the generalized Dirac operator $D_A$ is {\it
    neck compatible}  and {\it cylindrical} near the boundary provided
    the connection is in cylindrical form in a collar.

We are interested in decompositions of closed, oriented 3--manifolds $X$
    into two pieces $Y\cup_{\Si}Z$. Eventually $\Si$ will be a torus
    and $Y$ will be a solid
    torus, but for the time being $Y$ and
    $Z$ can be any 3--manifolds with boundary $\Si$.  Fix an orientation
    preserving identification of a tubular neighborhood of $\Si$
    with $[-1,1]\times\Si$ so that $ \{-1\}\times\Si$ lies in the
    interior of $Y$ and $ \{1\}\times\Si$ lies in the interior of
    $Z$. We identify $\Si$ with $\{0\}\times \Si$.
    As oriented boundaries, $\Si=\partial Y=-\partial Z$
    using the outward normal first convention.

\begin{figure}[ht!]\small
 \centering
\includegraphics[scale=.8]{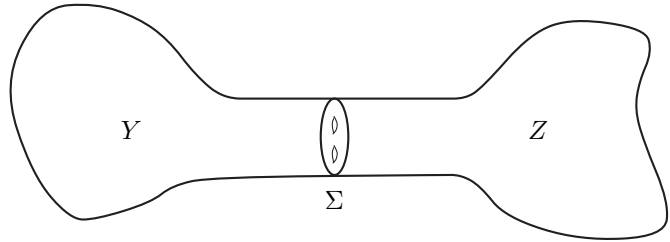}
\vskip-.7in $Y$\hskip2in $Z$\\
\vskip.2in$\Sigma$
 \caption{ The split 3--manifold $X$ }
 \label{dumbell}
 \end{figure}

To stretch the collar of $\Si$,  we introduce the notation
   \begin{eqnarray*} Y^R&=&Y\cup
    ([0,R]\times\Si) \\ Z^R&=&Z\cup ([-R,0]\times\Si)\end{eqnarray*}
    for all $R\ge 1$.  We
    also define $Y$ and $Z$ with infinite collars attached:
    \begin{eqnarray*} Y^\infty&=& Y\cup ([0,\infty)\times\Si) \\
    Z^\infty&=& Z\cup((-\infty,0]\times\Si).  \end{eqnarray*}

Notice that since $\Phi\circ D_A= J\circ (S_a+{\frac{\partial}{\partial
    u}})\circ \Phi$, the operator $D_A$ has natural extensions to $Y^R$,
    $Z^R$, $Y^\infty$, and $Z^\infty$.

\subsection{The spaces P$^+$ and P$^-$}\label{dfklkjdfa}

In this section we identify certain subspaces of the $L^2$ forms on
    $\Si$ associated to the operators $S_a$ and $D_A$. We first
    consider $L^2$ solutions to $D_A(\sigma,\tau)=0$ on $Y^\infty$
    and $Z^\infty$. Since $S_a$ is elliptic on the closed surface
    $\Si$, its spectrum is discrete and each eigenspace is a finite-dimensional
space of smooth forms.

Suppose $(\si,\tau) \in \Om^{0+1}_{Y^\infty} \otimes \CC^2$ is a
    solution to $D_A(\si,\tau)=0$ on $Y^\infty$. Following \cite{APS},
    write $\Phi(\si,\tau) = \sum c_\la(u) \phi_\la$ along
    $[0,\infty)\times\Si,$ where $\phi_\la \in \Om^{0+1+2}_\Si
    \otimes \CC^2 $ is an eigenvector of $S_a$ with eigenvalue
    $\la.$ Since $\Phi \circ D_A = J\circ (S_a + \tfrac{\partial}{\partial u})
    \circ \Phi$, it follows by hypothesis that\begin{eqnarray}\label{fourier}
    0 & = & ( S_a + \tfrac{\partial}{\partial u}) (\Phi(\si,\tau))
    \nonumber\\ & = & \sum_\la (\la c_\la + \tfrac{\partial
    c_\la}{\partial u}) \phi_\la  \end{eqnarray} hence $$c_\la(u) =
    e^{-\la u} b_\la$$ for some constants $b_\la$. Thus $(\si,\tau)\in L^2(
\Om^{0+1}_{Y^\infty} \otimes
    \CC^2)$ if and only if $c_\la(u) =0 $ for all $\la \leq 0.$

This implies that there is a one-to-one correspondence, given by
    restricting from $Y^\infty$ to $Y$,
    between the $L^2$ solutions to $D_A(\si,\tau) = 0$
    on $Y^\infty$ and the solutions to $D_{A}(\si,\tau)=0$ on $Y$
    whose restriction to the boundary $\Si$ lie in the {\it positive
    eigenspace} $P^+_a\subset L^2(\Om^{0+1+2}_\Si\otimes \CC^2)$ of
    $S_a$, defined by
    $$
     P^+_a=\sspan_{L^2}
    \{\phi_\la \mid \la>0\}. $$

Recalling that $\Si=\partial Y=-\partial Z$, we obtain a similar
one-to-one correspondence between  the space of $L^2$ solutions to $D_A(\si,\tau)
= 0 $
on   $Z^\infty$ and the space of  solutions to
    $D_{A}(\si,\tau)$ on $Z$ whose restriction to the boundary $\Si$ lie
    in the negative eigenspace $P^-_a\subset L^2(\Om^{0+1+2}_\Si\otimes
    \CC^2)$ of $S_a$, defined by
    $$
    P^-_a=\sspan_{L^2} \{\phi_\la \mid \la<0\}.$$

The spectrum of $S_a$ is symmetric and $J$ preserves the kernel of $S_a$
    since $S_a J=-J S_a$. In fact, $J$ restricts to an isometry $J\co P^+_a
    \lto P^-_a$. The eigenspace decomposition of $S_a$ determines the orthogonal
decomposition into closed
    subspaces\begin{equation}\label{decompofL2}
    L^2(\Om^{0+1+2}_\Si\otimes \CC^2) =P^+_a \oplus \ker S_a\oplus
    P^-_a.  \end{equation}

The spaces $P^\pm_a$ are isotropic subspaces and are Lagrangian if and
    only if $\ker S_a=0$. Since $\Si$ bounds the 3--manifold
    $Y$ and the operator $D_A$ is defined on $Y$, it is not hard to
    see that the signature of the restriction of $iJ$ to $\ker S_a$ is zero.
Hence
    $\ker S_a$ is a finite-dimensional  sub-symplectic space of
    $L^2(\Om^{0+1+2}_\Si\otimes \CC^2)$. The restrictions of the complex
structure
$J$
     and  
   the  inner product to $\ker S_a$ depend on the Riemannian metric, whereas the
    symplectic structure $\om(x,y)= \langle x,Jy\rangle$ depends on the
orientation
    but not on the metric.

An important observation is that if $L\subset \ker S_a$
    is any Lagrangian subspace, then $P^+_a \oplus L$ and $P^-_a
    \oplus L$ are infinite-dimensional Lagrangian subspaces of
    $L^2(\Om^{0+1+2}_\Si\otimes \CC^2)$.


If $a\in \Om^{1}_{\Si}\otimes su(2)$ is a {\it flat connection}, that is,
if the curvature 2--form $F_{a}=da+a\wedge
a$ is everywhere zero, then the kernel of $S_{a}$ consists of harmonic
forms, ie, $S_a(\alpha,\beta,\gamma)=0$ if and only if
$d_a\al=d_a\beta=d_a^*\beta=d_a^*\ga=0$.  
The Hodge and de~Rham theorems identify $\ker S_a$ with the
cohomology  group
$H^{0+1+2}(\Si;\CC^2_a)$, where $\CC^2_a$ denotes the local coefficient
    system determined by the holonomy representation of the flat connection $a$. 
Under this identification,  the induced symplectic structure on
$H^{0+1+2}(\Si;\CC^2_a)$ agrees with the direct sum of the
 symplectic structures on
    $H^{0+2}(\Si;\CC^2_a)$ and $H^1(\Si,\CC^2_a)$  given by the negative of
the cup product.  This is because the wedge products of differential forms
induces the cup product on de~Rham cohomology, and because of the 
 formula $$\om(x,y)= \langle x,Jy\rangle=-\int_\Si x\wedge y=-(x\cup
y)[\Si],$$
    where the   forms $x$ and $y$ are either
    both are 1--forms or  0-- and  2--forms, respectively.
    In this formula, we have suppressed the notation for the complex
    inner product on $\CC^2$ for the forms as well as in the cup product.
    Notice   that $H^0( \Si;\CC^2_a)$ and $H^2( \Si;\CC^2_a)$ are
    Lagrangian subspaces of $H^{0+2}( \Si;\CC^2_a)$.

\subsection{Limiting values of extended L$^2$ solutions and Cauchy data
    spaces}

Our next task is to identify the Lagrangian of { limiting values of
    extended $L^2$ solutions},  and
    its infinite-dimensional generalization, the Cauchy data
    spaces, in the case when $A$ is a flat connection
    in cylindrical form on a 3--manifold $Y$ with boundary $\Si$.

Atiyah, Patodi and Singer define the space of limiting values of
    extended $L^2$ solutions to $D_A \phi=0$ to be a certain
finite-dimensional
    Lagrangian subspace $$L_{Y,A} \subset\ker S_a,$$
    where $a$ denotes the restriction of $A$ to the boundary. We
    give a brief description of this subspace and refer to
    \cite{APS,JDGpaper} for further details.

First we define the Cauchy data spaces; these will be crucial in our
    later analysis. We follow \cite{Nicolaescu} closely; our terminology
    is derived from that article. In \cite{ Booss-Wojciechowski} it
    is shown that there is a well-defined, injective restriction
    map\begin{equation}\label{restrictionmap}
    r\co \ker\left(D_A\co L^2_{\frac{1}{2}}(\Om^{0+1}_Y\otimes \CC^2)\to
    L^2_{-{\frac{1}{2}}}(\Om^{0+1}_Y\otimes \CC^2) \right) \lto
    L^2(\Om^{0+1+2}_\Si\otimes \CC^2). \end{equation} Unique continuation
    for the operator $D_A$ (which holds for any generalized Dirac
    operator) implies that $r$ is injective.

\begin{definition} The image of $r$ is a closed, infinite-dimensional
    Lagrangian subspace of $L^2(\Om^{0+1+2}_\Si\otimes \CC^2)$. It
    is called the {\it Cauchy data space} of the operator $D_A$ on
    $Y$ and is denoted
    $$ \La_{Y,A}.$$
    Thus the Cauchy data space is the space of
    restrictions to the boundary  of solutions to $D_A(\si,\tau)=0$. It
    is shown in \cite{Nicolaescu} that $\La_{Y,A}$ varies continuously
    with the connection $A$.  \end{definition}

\begin{definition} The {\it limiting values of extended $L^2$ solutions}
    is defined as the symplectic reduction of $\La_{Y,A}$ with respect
    to the isotropic subspace $P^+_a$, the positive eigenspace of
    $S_a$. Precisely,
    $$  L_{Y,A}=
    \proj_{\ker S_a}\left(\La_{Y,A}\cap (P^+_a \oplus \ker S_a)\right)=
     \frac  { \La_{Y,A} \cap(P^+_a \oplus \ker S_a) }
    {\La_{Y,A} \cap P^+_a}  \subset\ker S_a .
   $$ \end{definition}

This terminology comes from \cite{APS}, where the restriction $r$ is
    used to identify the space of $L^2$ solutions of $D_A(\si,\tau)=0$
    on $Y^\infty$ with the subspace $\La_{Y,A}\cap P^+_a$, and the
    space of extended $L^2$ solutions with $\La_{Y,A}\cap
    (P^+_a\oplus\ker S_a)$. Thus $L_{Y,A}$ is the symplectic reduction
    of the extended $L^2$ solutions:\begin{equation}\label{extendeddd}
    L_{Y,A}= \frac{   \La_{Y,A} \cap(P^+_a \oplus \ker S_a) }  
    {\La_{Y,A} \cap P^+_a}   \cong{ \frac{ \hbox{Extended }\ L^2 \
    \hbox{solutions} } { L^2 \ \hbox{solutions} } } \end{equation}

We now recall a result of Nicolaescu on the ``adiabatic limit'' of the
    Cauchy data spaces \cite{Nicolaescu}. To avoid some technical
    issues, we make the assumption $\La_{Y,A}\cap P^{+}=0$;  in the
    terminology of \cite{Nicolaescu}, this means that {\it
    0 is a non-resonance level} for $D_A$ acting on $Y$.
    This assumption does not hold  in general, but it does hold in all the
cases
    considered in this article.

To set this up, replace $Y$ by $Y^R$ and extend $D_A$ to $Y^R$. This
    determines a continuous family $\La^R_{Y,A}=\La_{Y^R,A}$ of
    Lagrangian subspaces of $L^2(\Om^{0+1+2}_\Si\otimes\CC^2)$ by
    Lemma 3.2 of \cite{kirk-daniel}. The corresponding subspace
    $L_{Y,A}^R$ of limiting values of extended $L^2$ solutions  is
    independent of $R$.

Nicolaescu's theorem asserts that
as $R \to \infty$, $\La^R_{Y,A}$ limits to a certain
    Lagrangian.  Our assumption that 0 is a
    non-resonance level ensures that its limit is $L_{Y,A} \oplus
    P^-_a$.  Recall from Equation (\ref{decompofL2}) that
    $L^2(\Om^{0+1+2}_\Si\otimes\CC^2)$ is decomposed into the orthogonal
    sum of $P^+_a$, $P^-_a$, and $\ker S_a$. Notice also that the
    definition of $L_{Y,A}$ in Equation (\ref{extendeddd}) shows that it
    is independent of the collar length, ie, that $$\hbox{proj}_{\ker
    S_a} \left(\La^R_{Y,A} \cap (P^+_a\oplus \ker S_a)\right)$$ is
    independent of $R$. This follows easily from the eigenspace
    decomposition of $S_a$ in Equation (\ref{fourier}).

We now state Nicolaescu's adiabatic limit theorem \cite{Nicolaescu}, as
    sharpened in \cite{kirk-daniel}.

\begin{theorem}[Nicolaescu] \label{adiab} 
{\sl Assume that $\La_{Y,A}\cap
    P^+_a=0$ (equivalently assume that there are no $L^2$ solutions
    to $D_A (\si,\tau)=0$ on $Y^\infty$). Let $L_{Y,A}\subset \ker S_a$
    denote the limiting values of extended $L^2$ solutions. Then
    $$\lim_{R\to\infty} \La^R_{Y,A}=L_{Y,A}\oplus P^-_a,$$ with
    convergence in the gap topology on closed subspaces, and moreover
    the path of Lagrangians $$t\mapsto\begin{cases} \La^{1/(1-t)}_{Y,A}&
    t<1\\ L_{Y,A}\oplus P^-_a&t=1 \end{cases}$$ is continuous for
    $t\in[0,1]$ in the gap topology on closed subspaces.}
    \end{theorem}

Next we introduce some notation for the extended $L^2$ solutions. 
Although we use the terminology of extended $L^2$ solutions and limiting
values from \cite{APS}, it is more convenient for us to use the
characterization of these solutions in terms of forms on $Y$ with
$P^+_a\oplus \ker S_a$ boundary conditions.

\begin{definition} Let $\tV_A$ be the {\it space of  extended $L^2$
    solutions to $D_A (\si,\tau)=0$}. This is defined by setting
$$\tV_A=\{(\si,\tau)\in
    \Om_Y^{0+1}\otimes\CC^2 \mid D_A(\si,\tau)=0 \hbox{ and }
    r(\si,\tau)\in P^+_a\oplus \ker S_a\}.$$ Define also the {\it limiting
    value map} $ p\co \tV_A\lto \ker S_a$ by setting $p(\si,\tau) =
    \hbox{proj}_{\ker S_a} (r (\si,\tau))$ for $(\si,\tau)\in \tV_A$,
    where $r$ is the restriction map of Equation (\ref{restrictionmap}).
    Notice that $p(\tV_A)=L_{Y,A}$.  The choice of terminology is
explained by Equation (\ref{extendeddd}).\end{definition}

Let $\Th$ denote the trivial connection on $Y$ and $\th$ the trivial
    connection on $\Si = \partial Y.$ Let $\La_Y=\La_{Y,\Th}$ and
    $L_Y =L_{Y,\Th}$. The following theorem identifies $L_Y,$ the
    limiting values of extended $L^2$ solutions to $D_\Th(\si,\tau)=0$
    on $Y$. Since $\th$ is the trivial connection on
    $\Si$, $\ker S_\th$ can be identified with the (untwisted)
    cohomology $H^{0+1+2}(\Si;\CC^2)$.

\begin{theorem}\label{limvalattriv} 
{\sl Suppose $Y$ is a compact, oriented,
  connected 3--manifold with connected boundary $\Si$. Let $\Th$ be
the
    trivial connection on $Y$ and $\th$ the trivial connection on
    $\Si$. Identify $\ker S_\th$ with $H^{0+1+2}(\Si;\CC^2)$ using
    the Hodge theorem. Then the space of the limiting values of extended
    $L^2$ solutions decomposes as $$L_Y=H^0(\Si; \CC^2) \oplus
    \Image \left(H^1(Y; \CC^2)\to H^1(\Si;\CC^2) \right).$$}
\end{theorem}

\noindent {\bf Proof}\qua
   Proposition 4.2 of \cite{JDGpaper} says that if $D_\Th(\si,\tau)=0$ and
    $( \si,\tau)$ has boundary conditions in $P^+_\th \oplus \ker
    S_\th$ (ie, if $(\si,\tau)\in\tV$), then $d_\Th \si=0$, $d_\Th
    \tau=0$ and $d^*_\Th \tau=0$.  Regularity of solutions to this
elliptic boundary problem ensures that $\si$ and $\tau$ are smooth forms.

 If
$r(\si,\tau)=(\al,\beta,\gamma)$, then $\al\in\Omega^0_\Si\otimes \CC^2$
is a closed form whose cohomology class equals the restriction of the
cohomology class on $Y$ represented by $\si$. Similarly
$\beta \in\Omega^1_\Si\otimes \CC^2$  represents the restriction of
the cohomology class of $\tau$ to $\Si$.
Since the projection to harmonic forms does not change the
cohomology class of a closed form,
\begin{eqnarray*}
p(\tilde V) & \subset &\Image (H^{0+1}(Y;\CC^2) \to H^{0+1}(\Si;\CC^2)
)\oplus H^2(\Si,\CC^2) \\
&=&H^0(\Si;\CC^2)\oplus
\Image (H^{1}(Y;\CC^2) \to H^{1}(\Si;\CC^2))\oplus H^2(\Si;\CC^2).
\end{eqnarray*}

All of $H^0(\Si;\CC^2)$ is contained in $p(\tilde V)$,
since constant 0--forms on $\Si$ extend over $Y$, and if $\si$ is a
constant $0$ form on $Y$ then $(\si,0)\in \tilde{V}$ because its
restriction to the boundary lies in $\ker S_\theta$.  This implies that
$$p(\tilde V)\subset H^0(\Si;\CC^2)\oplus
\Image (H^{1}(Y;\CC^2) \to H^{1}(\Si;\CC^2)).$$

Since
$p(\tilde{V})$ is Lagrangian, it is a half 
dimensional subspace of $H^{0+1+2}(\Si,\CC^2)$.
Poincar\'e duality and the long exact
sequence of the pair $(Y,\Sigma)$ show that
$H^0(\Si;\CC^2)\oplus \Image (H^{1}(Y;\CC^2) \to
 H^{1}(\Si;\CC^2))$ has the same dimension, so they are
 equal. \endproof

Suppose $A$ is a flat connection on $Y$ with
 restriction $a = A|_\Si.$ Denote the  kernel of the limiting value map
 by $K_A=\ker( p\co \tV_{A}\lto
 \ker S_{a})$. By definition,
 $ K_A$   
 is the kernel of $D_{A}$
 on $Y$ with $P^{+}$ boundary conditions, but it can be characterized in
 several other useful ways. The eigenvalue
 expansion of  
 Equation
 (\ref{fourier}) implies that every form in $ K_A$
  extends to an exponentially decaying $L^2$ solution to
  $D_{A}(\sigma, \tau) = 0$ on $Y^\infty$.  Moreover, the
 restriction map $r$ of Equation (\ref{restrictionmap}) sends
$K_A$ injectively to $P^+_a$ by unique continuation,
and $r(K_A) =\La_{Y,A} \cap P^{+}_a$.
 For more details, see the fundamental
    articles  of Atiyah, Patodi, and Singer \cite{APS}
    and the book \cite{Booss-Wojciechowski}.

Suppose that $(\si,\tau)\in   K_A$.
    Then Proposition 4.2 of \cite{JDGpaper} implies that $d_A\si=0$,
    $d_A\tau=0$ and $d_A^*\tau=0$. Since $A$ is an $SU(2)$ connection,
    we have that $$d\langle \si,\si\rangle=\langle
    d_A\si,\si\rangle+\langle \si,d_A\si\rangle=0$$ pointwise. Thus
    the pointwise norm of $\si$ is constant. Since $\si$ extends to
    an $L^2$ form on $Y^\infty$, $\si=0$. Also $\tau$ is an $L^2$
    harmonic 1--form on $Y^\infty$. In \cite{APS} it is shown that
    if $A$ is a flat connection then
    the space of $L^2$ harmonic 1--forms on $Y^{\infty}$
    is isomorphic to $$\Image \left(
    H^1(Y,\Si;\CC^2_A) \to H^1(Y;\CC^2_A) \right),$$ the image of
    the relative cohomology in the absolute. Hence
     there is a short exact sequence $$0 \to \Image
    \left(H^1(Y,\Si;\CC^2_{A})\to H^1(Y;\CC^2_A)\right)\lto \tV_A\lto
    L_{Y,A} \to 0.$$ More generally, for any subspace $Q\subset
    \ker S_a$, restricting $p$ to $\tV_A\cap(P^+_a\oplus Q)$, one
    obtains the following very useful proposition.   

\begin{proposition}\label{exact} {\sl Suppose that $A$ is a flat connection
    on a 3--manifold $Y$ with boundary $\Si$. Let $a$ be the restriction
    of $A$ to $\Si$. If $Q\subset \ker S_a$ is any subspace (not
    necessarily Lagrangian), then there is a short exact sequence
    $$ 0 \to \Image \left(H^1(Y,\Si;\CC^2_{A})\to
    H^1(Y;\CC^2_A)\right)\lto \ker D_A(P^+_a \oplus Q)\mapright{p}
    L_{Y,A} \cap Q \to 0, $$ where $ \ker D_A(P^+_a \oplus Q)$ consists
    of solutions to $D_A(\si,\tau)=0$ whose restrictions to the boundary
    lie in $P^+_a \oplus Q$.

If $Q=0$, then this gives the isomorphisms $$\La_{Y,A}\cap P^+_a\cong
    K_A \cong \Image \left(H^1(Y,\Si;\CC^2_{A})\to
    H^1(Y;\CC^2_A)\right).$$  }\end{proposition}

\subsection{Spectral flow and Maslov index conventions}

If $D_t,\ t\in [0,1]$ is a 1--parameter family of self-adjoint operators
    with compact resolvents  and with $D_0$ and $D_1$ invertible, the spectral
    flow $SF(D_t)$\glossary{$SF(D_t)$}  is the algebraic
     number of eigenvalues crossing
    from negative to positive along the path. For precise definitions,
    see \cite{APS} and \cite{ CLM-splitting}.
  In case $D_0$ or $D_1$ is not invertible, we adopt the $(-\ep,-\ep)$
    convention to handle zero eigenvalues at the endpoints.

\begin{definition} Given a continuous 1--parameter family of self-adjoint
    operators with compact resolvents $D_t,\ t\in[0,1]$, choose $\ep>0$
    smaller than the modulus of the largest
    negative  eigenvalue of
    $D_0$ and $D_1$. Then the {\it spectral flow} $SF(D_t)$ is defined
    to be the algebraic intersection number in $[0,1] \times \RR$ of
    the track of the spectrum
$$\{(t,\la) \mid \ t\in [0,1], \ \la\in\hbox{Spec}(D_t) \}$$
    and the line segment from $(0,-\ep)$ to $(1,-\ep)$.
    The orientations are chosen so that if $D_t$ has spectrum $\{
    n+t \mid  n\in{\ZZ}\}$ then $SF(D_t)=1$.
\end{definition}

The proof of the following proposition is clear.

\begin{proposition} {\sl With the convention set above, the spectral flow is
    additive with respect to composition of paths of operators. It
    is an invariant of homotopy rel endpoints of paths of self-adjoint
    operators. If $\dim \ker D_t$ is constant, then $SF(D_t)
    =0.$  }  \end{proposition}

We will apply this definition to families of odd signature operators
    obtained from paths $A_t$ of $SU(2)$ connections. Suppose $A_t$
    is a path of $SU(2)$ connections on the closed 3--manifold $X$
    for $0 \leq t \leq 1$. We denote by
    $SF(A_t;X)$\glossary{$SF(A_t;X)$}  the spectral flow of the
    family   of odd signature operators $D_{A_t}$ on
    $\Om^{0+1}_X \otimes \CC^2.$ Since the space of all connections
    is contractible, the spectral flow $ SF(A_t;X)$ depends only on
    the endpoints $A_0$ and $A_1$ and we shall occasionally adopt
    the notation $SF(A_0,A_1;X)$  \glossary{ $SF(A_0,A_1;X)$} to
    emphasize this point.

\bigskip
We next introduce a compatible
    convention for the Maslov index \cite{Daniel}. A good reference for these ideas is
    Nicolaescu's article \cite{Nicolaescu}.  Let $H$ be a symplectic Hilbert
    space with compatible almost complex structure $J$.  A pair of
    Lagrangians $(L,M)$ in $H$ is called {\it Fredholm} if $L+M$ is
closed and both
    $\dim (L\cap M)$ and $\codim (L+M)$ are finite. We will say that
    two Lagrangians are {\it transverse} if they intersect trivially.

Consider a continuous path $(L_t,M_t)$ of Fredholm pairs of Lagrangians in
 $H$. Here, continuity is measured in
    the gap topology on closed subspaces. If
     $L_i$ is transverse to $M_i$
    for $i=0,1$, then the Maslov index $\Mas(L_t,M_t)$
    is  the number of times the
    two Lagrangians intersect, counted with sign and multiplicity.
    We choose the sign  so that if $(L,M)$ is a fixed
    Fredholm pair of Lagrangians such that $e^{sJ}L$ and $M$
    are transverse for all $0 \neq s \in [-\ep,\ep]$, then
    $\Mas( e^{ \ep(2t-1) J} L, M)=\dim(L \cap M)$.
    A precise definition is given in
    \cite{Nicolaescu} and more general properties of the Maslov
    index are detailed in \cite{CLM-Maslov, KL}.

    Extending the Maslov index to paths where the pairs at the
    endpoints
    are not transverse requires more care.  We  use  $e^{sJ}$,
    the  1--parameter group of symplectic transformations associated
    to $J$, to make them transverse.
    If $L$ and $M$ are any two Lagrangians, then $e^{sJ}L$ and $M$ are
    transverse for all small nonzero $s$. By \cite{kato}, the set
    of Fredholm pairs is open in the space of all pairs of Lagrangians.
    Hence, if $(L,M)$ is a Fredholm pair, then so is $(e^{sJ}L,M)$
    for all $s$ small.

\begin{definition}  Given a continuous 1--parameter family of Fredholm
    pairs of Lagrangians $(L_t,M_t),\ t\in[0,1],$ choose $\ep >0$
    small enough  that
    \begin{enumerate}
    \item[(i)] $e^{sJ}L_i$ is transverse to $M_i$  for
    $i=0,1$ and $0<s\leq\ep$, and
    \item[(ii)] $(e^{sJ}L_t,M_t)$ is
    a Fredholm pair for all $t\in[0,1]$
    and all $0\leq s\leq \ep$.
\end{enumerate}
Then define the  { Maslov index} of the pair $(L_t,M_t)$
    to be the Maslov index of $ (e^{\ep J}L_t, M_t)$.
\end{definition}

The proof of the following proposition is easy.

\begin{proposition}\label{masprop} {\sl With the conventions set above, the
    Maslov index is additive with respect to composition of paths.
    It is an invariant of homotopy rel endpoints of paths of Fredholm
    pairs of Lagrangians. Moreover, if  $\dim(L_t\cap M_t)$ is constant, then
    $\Mas(L_t,M_t)=0$. }   \end{proposition}

For 1--parameter families of Lagrangians $(L_t,M_t)$ 
which are transverse
except at one of the endpoints, the Maslov index  $\Mas(L_t,M_t)$
is often easy to compute.    

\begin{proposition}\label{princ}
{\sl Let $(L_t,M_t),\ t \in [0,1],$  be a continuous
1--parameter family of
Fredholm pairs of Lagrangians which are transverse for $t \neq 0.$
Suppose  $s\co \RR \to \RR$ is a smooth function
with $s(0)=0$ and $s'(0) \neq 0.$
Choose $\de >0$ so that $s(t)$ is strictly monotone on
$[0,\de]$ and $\ep>0$ with $\ep< | s(\de) |$ and $\ep < | s(-\de) | .$
Suppose further that,
for all $-\ep \leq r \leq \ep$ and all $0 \leq t \leq \de,$
the pair $(e^{rJ} L_t,M_t)$ satisfies
$$\dim(e^{rJ}L_t\cap M_t)=
\begin{cases} 
\dim(L_0 \cap M_0) &\text{if } r=s(t)\\
0&\text{otherwise.}\end{cases}$$

Then
$$\Mas(L_t,M_t)= 
\begin{cases} 
- \dim(L_0 \cap M_0)&\text{if } s'(0)>0\\ 
0&\text{if } s'(0)<0.\end{cases}$$}
\end{proposition}

\noindent {\bf Proof}\qua 
Write 
$$\Mas(L_t,M_t) = 
\Mas(L_t,M_t; \ 0 \leq t \leq \de) +
\Mas(L_t,M_t; \ \de \leq t \leq 1).$$
Since $L_t$ and $M_t$
are transverse for $t \in [\de,1]$,
it follows that $$ \Mas(L_t,M_t; \ \de \leq t \leq 1)=0.$$
The convention for dealing with non-transverse endpoints now applies to
show
that
$$\Mas(  L_t,M_t) = \Mas(  L_t,M_t; \ 0 \leq t \leq \de) = 
\Mas(  e^{\ep J} L_t,M_t; \ 0 \leq t \leq \de).$$
 
If $s'(0) <0,$ then $s(t)$ is monotone decreasing  on $[0,\de]$
and the hypotheses imply that $e^{\ep J} L_t$ and $M_t$ 
are transverse for $t \in [0,\de]$.
Hence $ \Mas(L_t,M_t) = 0$ as claimed.

On the other hand, if $s'(0)>0,$ then
we write  
\begin{eqnarray*}
\Mas(e^{\ep J} L_t,M_t; \ 0 \leq t \leq \de) &=&
\Mas( e^{\ep(1-2t) J} L_0,M_0; \ 0 \leq t \leq 1) \\
&&+ 
\Mas( e^{-\ep J} L_t,M_t; \ 0 \leq t \leq \de) \\
&&+ \Mas( e^{\ep(2t-1) J} L_\de,M_\de; \ 0 \leq t \leq 1).
\end{eqnarray*}
Since $s(t)$ is now monotone increasing on $[0,\de]$,
the hypotheses imply that $e^{-\ep J} L_t$ and $M_t$
are transverse for $t \in [0,\de]$.
Furthermore, by choosing $\ep$ smaller, if necessary,
we can assume that $ e^{\ep(2t-1) J} L_\de$ and $M_\de$
are transverse for all $t \in [0,1]$. 
Hence 
\begin{eqnarray*}
\Mas(L_t, M_t) &=& \Mas( e^{\ep(1-2t) J} L_0,M_0) \\
&=& -\Mas (e^{\ep(2t-1) J} L_0,M_0) = -\dim (L_0,M_0)
\end{eqnarray*}
by our sign convention.
\endproof
\noindent{\bf Remark}\qua  There is a similar result
for pairs $(L_t,M_t)$ which are transverse for $t \neq 1$.
If $s(t)$ is a smooth function satisfying the analogous
conditions, namely that $s(1)=0, \ s'(1) \neq 0$ and 
$$\dim (e^{rJ}L_t \cap M_t) =  
\begin{cases} 
\dim(L_1 \cap M_1) &\text{if } r=s(t)\\
0&\text{otherwise,}\end{cases}$$
then 
$$\Mas(L_t,M_t)= 
\begin{cases} 
\dim(L_1 \cap M_1)&\text{if } s'(1) < 0\\ 
0&\text{if } s'(1) > 0.\end{cases}$$
The details of the proof  are left to the reader.

\subsection{Nicolaescu's decomposition theorem for spectral flow}

The spectral flow and Maslov index are related by the following result
    of Nicolaescu, which holds in the more general context of neck
    compatible generalized Dirac operators. The following is the
    main theorem of \cite{Nicolaescu}, as extended in \cite{Daniel},
    stated in the context of the odd signature operator $D_A$ on a
    3--manifold.

\begin{theorem}\label{nic1} {\sl Suppose $X$ is a 3--manifold decomposed along
    a surface $\Si$ into two pieces $Y$ and $Z$, with $\Si$
    oriented so that $\Si=\partial Y = -\partial Z$.  Suppose $A_t$ is
    a continuous path of $SU(2)$ connections on $X$ in cylindrical
    form in a collar of $\Si$. Let $\La_Y(t) = \La_{Y,A_t}$ and
    $\La_Z(t) = \La_{Z,A_t}$ be the Cauchy data spaces associated
    to the restrictions of $D_{A_t}$ to $Y$ and $Z.$ Then $(\La_Y(t),
    \La_Z(t))$ is a Fredholm pair of Lagrangians and
    $$SF(A_t;X)=\Mas(\La_Y(t),\La_Z(t)).$$  }\end{theorem}

There is also a theorem for manifolds with boundary, see
    \cite{Nicolaescu2, Daniel2}. This requires the introduction of
    boundary conditions.  The following is not the most general
notion, but suffices for our exposition. See
\cite{Booss-Wojciechowski, KL} for a more detailed analysis of 
elliptic boundary conditions.

\begin{definition} Let $D_A$ be the odd signature operator twisted by a
    connection $A$ on a 3--manifold $Y$ with non-empty boundary $\Si$. A
    subspace $\tP\subset L^2(\Om^{0+1+2}_\Si\otimes\CC^2)$ is
    called a {\it self-adjoint Atiyah--Patodi--Singer (APS) boundary
    condition for $D_A$} if $ \tP$ is a Lagrangian subspace and if,
    in addition, $ \tP$ contains all the eigenvectors of the tangential
    operator $S_a$ which have sufficiently large positive eigenvalue as
    a finite codimensional subspace. In other words, there exists a
    positive number $q$ so that $$\{\phi_\la \mid
    S_a(\phi_\la)=\la\phi_\la\hbox{ and }\la>q\}\subset \tP$$ with
    finite codimension.  \end{definition}

\begin{lemma}\label{bcok}
{\sl Suppose that $X=Y\cup_\Sigma Z$ and $A_0$, $A_1$ are 
$SU(2)$ connections  in cylindrical form on the collar of $\Sigma$
as above. Let $\tP_0$  (resp.~$\tP_1$) be  a self-adjoint APS
boundary condition for $D_{A_0}$  (resp.~$D_{A_1}$) restricted to
$Y$.

Then 
$$(\La_{Y,A_0}, \La_{Z,A_1}),  \
(\La_{Y,A_0},\tP_1), \ (J\tP_0, \La_{Z,A_1}), \
\text{ and }(J\tP_0,\tP_1) 
$$
are Fredholm pairs.}
\end{lemma}
\noindent {\bf Proof}\qua Let $S_{a_0}$ and $S_{a_1}$ denote the tangential
operators of $D_{A_0}$ and $D_{A_1}$.  It is proved in
\cite{Booss-Wojciechowski} that 
\begin{enumerate}
\item The
$L^2$--orthogonal  projections to $\La_{Y,A_0}$ and $\La_{Y,A_1}$
are zeroth--order
pseudo-differential operators whose principal symbols  are just  the
projections onto the positive eigenspace of the principal symbols of
$S_{a_0}$ and $S_{a_1}$, respectively.   
\item If $Q_0$ and $Q_1$ denote the  
$L^2$--orthogonal  projections to the positive eigen\-spans
of   $S_{a_0}$ and $S_{a_1}$, respectively, then $Q_0$ and $Q_1$ are
zeroth--order pseudo-differential operators whose principal symbols
are also the projections onto the positive eigenspaces of the principal
symbols of
$S_{a_0}$ and $S_{a_1}$.  
\end{enumerate}
From the definition one sees that the difference 
$S_{a_0}-S_{a_1}$  is a zeroth order differential operator, and
in particular the  principal symbols  of $S_{a_0}$ and
$S_{a_1}$ coincide. Hence 
$$\sigma(Q_0)=\sigma(Q_1)=\sigma(\text{proj}_{\La_{Y,A_0}})=
\sigma(\text{proj}_{\La_{Y,A_1}}),$$
where $\sigma$ denotes the principal symbol.
Moreover, $Q_i$ and the projection to
$\tP_i$   differ by a finite-dimensional projection.   This
implies that the projections to
$\La_{Y,A_0}$, 
$\La_{Y,A_1}$, $\tP_0$, and $\tP_1$  are compact perturbations of
$Q_0$. The lemma follows from this and the fact that viewed from
the ``$Z$ side,'' the roles of the positive and negative spectral
projections are reversed.
\endproof

It follows from the results of \cite{APS} (see also
\cite{Booss-Wojciechowski})   that restricting the domain of
$D_A$ to 
$
    r^{-1}(\tP) \subset L^2_1(\Om^{0+1}_Y\otimes\CC^2)$ yields a
self-adjoint
    elliptic operator. Moreover, unique continuation for solutions
    to $D_A(\si,\tau)=0$ shows that the kernel of $D_A$ on $Y$ with
    APS boundary conditions $ \tP$ is mapped isomorphically by the
    restriction map $r$ to $\La_{Y,A}\cap \tP$.

A generalization of Theorem \ref{nic1}, which is also due to Nicolaescu
    (see \cite{Nicolaescu2} and \cite{Daniel}), states the
    following.

\begin{theorem}[Nicolaescu] \label{nic2} {\sl Suppose $Y$ is a 3--manifold
    with boundary $\Si.$ If $A_t$ is a path of connections on $Y$
    in cylindrical form near $\Si$ and $ \tP_t$ is a continuous
    family of self-adjoint APS boundary conditions, then the spectral
    flow $SF( A_t;Y; \tP_t)$\glossary{ $SF( A_t;Y; \tP_t)$} is well
    defined and   $$SF(A_t;Y; \tP_t)=\Mas(\La_Y(t), \tP_t).$$}
\end{theorem}

\section{Splitting the spectral flow for Dehn surgeries}
    \label{splitme}

In this paper, the spectral flow theorems described in the previous
section will be applied to homology 3--spheres $X$ obtained by Dehn
surgery on a knot, so $X$ is decomposed as
$X=Y\cup_{\Si}Z$ where $Y=D^{2}\times S^{1}$ and $\Si=\partial Y$ is
the 2--torus.  In our examples, $Z$ will be the complement of a knot in
$S^{3}$, but
the methods work just as well for knot complements in other homology
spheres.

This section is devoted to proving a splitting theorem for
$\CC^{2}$--spectral flow of the
odd signature operator for paths of $SU(2)$ connections with certain
properties.
In the end, the splitting theorem expresses the spectral flow  as a sum of
two terms, one involving $Z$ and the other
involving $Y$.

\subsection{Decomposing X  along a torus}\label{torus decomp}
We make the following assumptions, which
    will hold for the rest of this article.

\begin{enumerate}   \item The surface $\Si$ is the torus $$T=S^1\times
    S^1=\{(e^{ix},e^{iy})\},$$ oriented so that the  1--forms
    $dx$ and $dy$ are ordered as $\{dx,dy\}$ and with the product
    metric, where the unit circle $S^1 \subset \CC$ is given the standard
    metric. The torus $T$ contains the two curves $$\mu=\{(e^{ix},1)\}
    \quad \text{and} \quad \la=\{(1,e^{iy})\},$$ and $\pi_1(T)$ is
    the free abelian group generated by these two loops.

  \item The 3--manifold $Y$ is the solid torus $$Y =D^2\times S^1 = \{
    (re^{ix},e^{iy}) \mid 0\leq r \leq 1\},$$ oriented so that $dr dx
    dy$ is a positive multiple of the volume form when $r>0$.
    The fundamental group $\pi_1(Y)$ is infinite cyclic generated
    by $\la$ and the curve $\mu$ is trivial in $\pi_1Y$ since it
    bounds the disc $D^2\times\{1\}$.
 There is a product metric on $Y$  such that a collar neighborhood of the
 boundary may be isometrically identified with $[-1,0]\times T$ and
 $\partial Y=\{0\}\times T$.
 The form $dy$ is a globally defined 1--form
    on $Y$, whereas the form $dx$ is well-defined off the core circle
    of $Y$ (ie, the set where $r=0$).   

  \item The 3--manifold $Z$ is the complement of an open tubular
  neighborhood of a knot in a homology sphere.   Moreover, we assume that the
    identification of $T$ with $\partial Z$ takes the loop $\la$ to
    a null-homologous loop in $Z$.

\smallskip \noindent
    There is a  metric on $Z$  such that a collar neighborhood of the
 boundary may be isometrically identified with $[0,1]\times
T$.  As oriented manifolds,  $\partial Z=-\{0\}\times T$.  The form
$dx$ on $\partial Z$ extends to a closed 1--form on $Z$  generating
the first cohomology $H^{1}(Z;\RR)$ which we continue to denote  
$dx$.

  \item The closed 3--manifold $X=Y\cup_{T}Z$ is a homology sphere.  The metric
  on $X$ is compatible with those on $Z$ and $Y$ and  $T$ is identified
  with the set $\{0\}\times T$ in the neck.
\end{enumerate}

\subsection{Connections in normal form and the moduli space of T}
\label{compsofco} Flat connections on the torus play a
central role here, and in this subsection 
we describe a  2--parameter family
of flat connections on the torus and discuss its relation to the flat
moduli space.

For notational  convenience, we identify elements of $SU(2)$ with
 unit quaternions via
    $$\begin{pmatrix} \al  & \be \\ -\bar \be  & \bar \al
    \end{pmatrix}\leftrightarrow \al+\be  j$$ where $\al,\be \in
    \CC$ satisfy $|\al|^2+|\be|^2 =1$. The Lie algebra $su(2)$
    is then identified with the purely imaginary quaternions $$\begin{pmatrix}
    ix &y+iz\\-y+iz& -ix \end{pmatrix}\leftrightarrow xi+yj+zk $$ for
    $x,y,z\in\RR$.

With these notational conventions, the action of $su(2)$ on $\CC^2$ can
    be written in the form
 $$(ix+jy+kz)\cdot (v_1e_1+v_2 e_2)=(ixv_1+(y+iz)v_2)e_1
    - ((y-iz)v_1+ixv_2)e_2.$$ In particular,
\begin{equation}\label{action}
    xi\cdot (v_{1}e_1+v_2e_2)=ixv_{1}e_1 -ixv_2e_2.
\end{equation} This corresponds
    to the standard inclusion $U(1)\subset SU(2)$ sending
    $\al \in U(1)$ to ${\rm diag}(\al, \al^{-1}) \in SU(2)$.
    On the level of Lie algebras,  this is the inclusion
    $u(1)\subset su(2)$ sending $ix$ to ${\rm diag}(ix, -ix)$.

\begin{definition} \label{Anf}   For $ (m,n)\in \RR^{2}$, let
$ a_{m,n}= -m i dx - n i dy$ and define the {\em
connections in normal form on $T$}  to be the set
$$\cA_{\rm nf}(T) = \{ a_{m,n} \mid 
(m,n)\in \RR^{2}\}.$$
    An $SU(2)$ connection $A$ on  $Z$ or $Y$ is said to be in
    {\it normal form along
    the boundary} if it is in cylindrical form on the collar
    neighborhood of $T$ and
    its restriction to the boundary is in normal form.
    \end{definition}
Notice that if $a=a_{m,n},$ then $\hol_{a}(\mu)=e^{2\pi im}$ and
$\hol_{a}(\la)=e^{2\pi i n}$.  
The relevance of connections in normal form is made clear by the
following proposition, which follows from a standard gauge fixing
argument.  We
will call a connection {\em diagonal} if its connection
1--form takes values in the diagonal Lie subalgebra $u(1)\subset su(2)$.

\begin{proposition}
{\sl Any flat $SU(2)$ connection on $T$ is gauge equivalent to a
diagonal connection.  Moreover, any flat diagonal $SU(2)$
connection on $T$ is gauge equivalent via a gauge
 transformation $g\co T\to U(1)\subset SU(2)$ to a connection in normal
 form, and the normal form connection is unique if $g$ is required to
 be homotopic to the constant map $\id \co T\to \{ \id \}\subset
 U(1)$.}
 \end{proposition}

  We will introduce a special gauge group for the set of
  connections in normal form in Section \ref{relative gauge group},
  but for now note that any constant gauge transformation of the form $\cos(s)j
  + \sin(s)k$ acts on $\cA_{\rm nf}(T)$ by sending $a_{m,n}$ to
  $a_{-m,-n}$.  Alternatively, one can view this as interchanging the
  complex conjugate eigenvalues of the $SU(2)$ matrices in the
  holonomy representation.

 For any manifold $X$ and compact Lie group $G$, denote by   $\fR_{G}(X)$ the
    space of conjugacy classes of  representations $\rho\co  \pi_1 X \to G, $
ie,
    $$\fR_{G}(X)=\mbox{Hom}(\pi_1X,G)/\mbox{conjugation}, $$
 and denote by $\fM_{G}(X)$ the space of flat connections on
 principal $G$--bundles over $X$ modulo gauge transformations of those
 bundles.  In all cases considered here, $G=SU(n),\ n=2,3$ and $\dim
 X\leq 3$, so all $G$--bundles over $X$ are necessarily trivial.
   The association to each flat connection   its holonomy
 representation provides a homeomorphism
 $$\hol\co \fM_{G}(X) \stackrel{\cong}\lto \fR_{G}(X),$$ so we will use
 whichever interpretation is convenient.

  By identifying $\cA_{\rm nf}(T)$ with $\RR^{2}$,
  the moduli space $\fM_{SU(2)}(T)$ of flat connections (modulo the
  full gauge group) can be identified with the quotient of
  $\RR^{2}$ by the semidirect product of $\ZZ /2$
  with $\ZZ^2$, where $\ZZ/2$ acts by
  reflections through the origin and $\ZZ^{2}$ acts by
  translations.  The quotient map is a branched covering.  Indeed, setting 
$f(m,n)=[\hol_{a_{m,n}}\co \pi_{1}T\to SU(2)]$ for $(m,n) \in \RR^2$ defines
the  branched covering map
  \begin{equation}\label{pillowpara}
  f\co \RR^{2}\to \fR_{SU(2)}(T).
  \end{equation}

Since the connection 1--form of any $a\in \cA_{\rm nf}(T)$ takes values in
$u(1)\subset su(2)$, the
    twisted cohomology splits $$H^{0+1+2}(T;\CC^2_a) =
    H^{0+1+2}(T;\CC_{\hat{a}})\oplus H^{0+1+2}(T;\CC_{- \hat{a}}),$$ where $\pm
    \hat{a}$ are the $u(1)$ connections given by the reduction of
    the bundle.  Similarly, the de~Rham operator splits as
    \begin{equation}\label{Sasplit}
    S_a =
    S_{\hat{a}}\oplus S_{-\hat{a}},\end{equation} where   $S_{\pm \hat{a}} \co 
    \Om^{0+1+2}_T\otimes \CC \to  \Om^{0+1+2}_T\otimes \CC$  are
    the de~Rham operators associated to the $u(1)$ connections $ \pm
    \hat{a}$.

We leave the following cohomology calculations to the reader.
(See Equation (\ref{ccspan})
    for the definition of $\ccspan$.)

\begin{enumerate}

\item The flat connection $a_{m,n}\in \cA_{\rm nf}(T)$
 is gauge equivalent to the trivial connection if and only
    if $(m,n)\in\ZZ^2$. Moreover,\begin{equation}\label{cohomolofT}
    H^{0+1+2}(T; \CC^2_a)=
\begin{cases}0&\text{if $ (m,n)\not\in\ZZ^2$,}\\
\ccspan\{ 1,dx,dy,dxdy \} & \text{if $(m,n)=(0,0)$.}  \end{cases}
\end{equation}

  \item If $A$ is a flat $SU(2)$ connection on $Y$ in normal form
along the boundary (so $A|_{T}=a_{m,n} = -m i dx - n i dy$ with $m\in\ZZ$), then
$A$ is gauge equivalent to the trivial connection if and only
if $n\in \ZZ$.
    Moreover, \begin{equation}\label{cohomolofW} H^{0+1}(Y; \CC^2_A)=
    \begin{cases}0&\text{if $ n\not\in\ZZ$,}\\ \ccspan \{ 1, dy
    \} & \text{if $n=0$}  \end{cases}  \end{equation}

  \item   For the trivial connection $\Th$  on $Z$,  the
    coefficients are untwisted and $
    H^{0+1}(Z; \CC^2)=\sspan_{\CC^{2}}\{ 1,dx\}$.
       \end{enumerate}

In terms of the limiting values of extended $L^2$ solutions, these
    computations together with Theorem \ref{limvalattriv} give the
    following result.

\begin{proposition}\label{compoflim} {\sl The spaces
of limiting values of extended $L^{2}$ solutions for the trivial
connection on $Y$ and $Z$ are $L_Y=\ccspan \{ 1,\ dy
    \}$ and $L_Z=\ccspan\{1,dx\}$ respectively.}
    \end{proposition}

\subsection{Extending connections in normal form on T over Y}
    \label{famofconn}
The main technical difficulty in the present work has at its core the
    special nature of the trivial connection. We begin by specifying
    a 2--parameter family of connections on $Y$ near $\Theta$
    which extend the
    connections on normal form on $T$.  We will use these connections
    to build paths of connections on $X$ which start at the trivial
    connection and, at first, move away in a specified way that is
    independent of $Z$ and $Y$ except through the homological
    information in the identification of their boundaries (which
 determine our coordinates on $T$).

Choose once and for all
a smooth non-decreasing cutoff function  $q\co [0,1]\to [0,1]$  with
    $q(r)=0$ for $r$ near $0$ and $q(r)=1$ for $r$ near enough to $1$
    that $(re^{ix},e^{iy})$ lies in the collar neighborhood of $T$.

For each point $(m,n) \in \RR^2$, let $A_{m,n}$ be the connection in
    normal form on the solid torus $Y$ whose value at the point
    $(re^{ix},e^{iy})$ is\begin{equation} \label{eqnofbmn}
    A_{m,n}(re^{ix},e^{iy})= -q(r) m i dx - n i dy.\end{equation}
    This can be thought of as a $U(1)$ connection, or as an $SU(2)$
    connection using quaternionic notation. Notice that $A_{m,n}$
    is flat if and only if $m=0$, and in general is flat away from an 
annular
    region in the interior of $Y$.

\subsection{Paths of connections on X and adiabatic limits at
    ${\Th}$} \label{defnofAt}

Suppose $X$ is a homology 3--sphere decomposed as $X= Y \cup_T Z.$ For
    the rest of this section, we will suppose that $A_t, \ t \in
    [0,1]$ is a continuous path of $SU(2)$ connections on $X$ satisfying
    the following properties:
    \begin{enumerate}   
    \item $A_0 =
    \Th$, the trivial connection on $X$, and $A_1$ is a flat
    connection on $X$. 
    
    \item  The restriction of $A_{t}$ to the neck is a path of
    cylindrical normal form connections
    $$A_{t}|_{[-1,1]\times T}=a_{m_{t},n_{t}}$$ for some
    piecewise smooth   path
    $(m_{t},n_{t})$ in $\RR^{2}$ with $(m_{t},n_{t})\not\in \ZZ^{2}$
    for $0<t\leq 1$. 
    
    \item There exists a small number $\de >0$ such that, for
    $0<t\leq \de $,
    \begin{enumerate}
    \item $(m_{t},n_{t})=(t,0)$,
    \item $A_{t}|_{Z}=-tidx$ and $A_{t}|_{Y}=-q(r)tidx$, and
    \item $\Delta_{Z}(e^{i2\pi t})\neq 0,$ where $\Delta_{Z}$ denotes
    the Alexander polynomial of $Z$.
    \end{enumerate}
    \end{enumerate}

Most of the time we will assume that the restriction of $A_t$ to $Z$ is
    flat for all $t$, but this is not a necessary hypothesis in Theorem
    \ref{zorro}. This extra bit of generality can be useful in contexts
    when the space $\fR_{SU(2)}(Z)$ is not connected.

The significance of the condition involving the Alexander polynomial
is made clear by the following lemma and corollary.

\begin{lemma} {\sl If $A_t$ is a path of connections satisfying conditions
    1--3 above and if $\de >0$ is the constant in condition 3,
    then $H^1(Z,T;\CC^2_{A_t})  =0$ for $0 \leq t \leq \de $.}
    \end{lemma}

\noindent{\bf Sketch of Proof}\qua For $A_0= \Th,$ the trivial connection,
    this follows from the long exact sequence in cohomology of the
    pair $(Z,T)$ for $t=0$. Using the Fox calculus to identify the
    Alexander matrix with the differential on 1--cochains in the infinite
    cyclic cover of $Z$ proves the lemma for $0 < t \leq \de $.
    A very similar computation is carried out in \cite{klassen-thesis}.
    \endproof

\begin{corollary} \label{adiabofX} {\sl With the same hypotheses as above, the
$L^2$ kernel of
    $D_{A_t}$ on $Z^\infty$ is trivial for $0 \leq t \leq \de $.
    Equivalently, letting $\La_Z(t) = \La_{Z,A_t}$, then for $ 0
    \leq t \leq \de $, $$\La_{Z}(t) \cap P^-_{a_t}=0.$$
    Furthermore, letting $\La^R_Z(t) =
    \La_{Z^R,A_t}$,
    $$\lim_{R \to \infty} \La^R_{Z}(t) =
    \begin{cases} L_Z \oplus P^+_{\th} & \text{if $t=0$} \\ P^+_{a_t} &
    \text{if $0<t \leq \de $.}\end{cases}$$}
    \end{corollary}
    \noindent {\bf Proof}\qua The first claim follows
    immediately from Proposition
    \ref{exact} applied to $Z$ with $K=0$. (The orientation conventions,
    as described in Section \ref{dfklkjdfa} explain why $P^-$ is
    used instead of $P^+$.)  In the terminology of \cite{Nicolaescu}, this
    means that $0$ is
    a non-resonance level for $D_{A_t}$ for $ 0 \leq t \leq \de $.
   Applying Theorem \ref{adiab}, Theorem
    \ref{limvalattriv}, and Equation  (\ref{cohomolofT}) gives the
    second claim.
    \endproof

\subsection{Harmonic limits of positive and negative
eigenvectors} \label{harmlimits}

In this section, we investigate some limiting properties of  the
    eigenvectors of $S_a$ where $a$ ranges over a neighborhood of the
    trivial connection $\theta$
    in the space of connections in normal form on $T$.

Let  $s \in \RR$ be a fixed number. (Throughout this
subsection, $s$ is a
fixed angle. In
Theorem \ref{counterexample}, the value $s=0$
is used.)  Consider the path of connections
$$a_t =-t \cos(s) i dx - t \sin(s) i dy$$
for $0 \leq t \leq
    \de $.  Notice that  $a_t$ is a
    path of connections in normal form approaching the trivial connection $\th$
and the angle of approach  is $ s$.

The path of operators $S_{a_t}$ is an analytic (in $t$) path of elliptic
    self-adjoint operators. It follows from the results of analytic
    perturbation theory that $S_{a_t}$ has a spectral decomposition
    with analytically varying eigenvectors and eigenvalues (see
    \cite{kato, KKillinoisJ}). By Equation (\ref{cohomolofT}) we have $$
    \dim (\ker S_{a_t}) =\begin{cases} 8 & \text{if $t=0$} \\ 0 &
    \text{if $0<t\leq\de $}\end{cases}$$ and $$ \ker S_\th = \ccspan \{ 1,\
dx,\ dy, \ dx dy \}.$$

Since the spectrum of $S_{a_t}$ is symmetric, it follows that for $t>0$
    there are four linearly independent positive eigenvectors and four
negative eigenvectors
    of $S_{a_t}$ whose eigenvalues limit to $0$ as $t \to 0^+$, ie, the
    eigenvectors limit to (untwisted) $\CC^2$--valued harmonic
    forms. More precisely, there exist 4--dimensional subspaces $K^+_s$
    and $K^-_s$ of $\ker S_{\th}$ so that $$\lim_{t\to
    0^+}P^+_{a_t}=K^+_s \oplus P^+_{\th} \quad \text{and} \quad
    \lim_{t\to 0^+}P^-_{a_t}=K^-_s \oplus P^-_{\th} .$$  In particular,
    the paths of Lagrangians 
    $$t\mapsto\begin{cases} K^+_s \oplus
    P^+_{\th}&\hbox{ if } t=0 \\ P^+_{a_t}&\hbox{ if } 0 < t \leq
    1\\\end{cases} \quad \text{ and } \quad t\mapsto\begin{cases} K^-_s \oplus
    P^-_{\th}&\hbox{ if } t=0 \\ P^-_{a_t}&\hbox{ if } 0 < t \leq
    1\\\end{cases}$$ are continuous.

The finite-dimensional Lagrangian subspace  $K^+_s$ will be used to
    extend the boundary conditions $P^+_a$  to a continuous family
    of boundary conditions up to $\th$. Similarly, $K^-_s$ will be
    used to extend  the boundary conditions $P^-_a$. The next
    proposition gives a useful description of these spaces.

\begin{proposition}\label{hrighgh} {\sl Define
    the 1--form $\xi_s = -\cos( s) i dx - \sin( s) idy.$ Consider the
    family of connections on $T$ given by $a_t = t\xi_s$ for $t \in
    [0,\delta]$. If $K^+_s$ and $K^-_s$ are defined as above, then
    \begin{eqnarray*} 
    K^+_s &=& \sspan \{ (1 - *\xi_s)\otimes e_1, \
    (\xi_s - dx dy)\otimes e_1, \\
    && \qquad (1 + *\xi_s)\otimes e_2, \ (-\xi_s- dx dy)\otimes e_2 \}, \\ 
    K^-_s &=&\sspan \{ (1+*\xi_s)\otimes
    e_1, \ (\xi_s + dx dy)\otimes e_1,  \\
    && \qquad (1-*\xi_s)\otimes e_2, \ (-\xi_s+ dx dy)\otimes e_2\}.
    \end{eqnarray*}}
\end{proposition}

\noindent {\bf Proof}\qua  Recalling the way a diagonal connection acts on the two
factors of $\CC^{2}$ from Equation (\ref{action}), we can decompose
$K^{\pm}_{s}$ into $K^{\pm}_{s}=\hat K ^{\pm}_{s} \oplus \hat
K^{\pm}_{-s}$ where $\hat K^{\pm}_{s}$ is the space of harmonic
limits of the operator $S_{\hat a_{t}}$ in Equation (\ref{Sasplit}).

Now $$S_{\hat{a}_t}(\al,\be,\ga)=S_{\th}(\al,\be,\ga)+ t \Psi_s(\al,
    \be, \ga),$$ where $\Psi_s(\al,\be,\ga) = (*(\xi_s \be), -*(\xi_s \al
    ) - \xi_s(* \ga), \xi_s(*\be)).$
A direct computation shows  that $\Psi_s(1,-*\xi_s,0) = (1,-*\xi_s,0) $ and
 $\Psi_s(0,\xi_s,-dxdy) = (0,\xi_s,-dxdy)$. Since $-\hat{a}_t=-\xi_s$,
    it follows that $$  \{ (1 - *\xi_s)\otimes e_1, \ (\xi_s - dx
    dy)\otimes e_1,  (1 + *\xi_s)\otimes e_2, \ (-\xi_s - dx dy)\otimes
    e_2 \} \subset K^+_s.$$ The first formula then follows since
    both sides are 4--dimensional subspaces of $\ker S_\th$.

The result for $K^-_s$ can also be computed directly. Alternatively,
it can obtained from the result for
$K^{+}_{s}$ by applying $J,$ using the fact that
$S_a J = - J S_a$ and so
$ K^-_s = J K^+_s$.\endproof

Comparing these formulas for $K^+_s$ and $K^-_s$ with that for $L_Z$
    from Proposition \ref{compoflim} yields the following important
    corollary.

\begin{corollary}\label{trasnveree} {\sl For $s=\tfrac{\pi}{2}$ or
$\tfrac{3\pi}{2}$, $\dim K^{\pm}_{s}\cap L_Z=2$  and for $s=0$ or $\pi $,
$\dim K^{\pm}_{s}\cap L_Y=2$.  For values of $s$ other than those
specified, the intersections are trivial.}
\end{corollary}

Next, we present an example which, though  
peripheral to the main thrust of this article, shows that extreme
care must be taken when dealing with paths of adiabatic limits of Cauchy
data spaces.
For the sake of argument, suppose that we could
replace the path of
Cauchy data spaces with the path of the adiabatic limits of the Cauchy data
spaces.  This would reduce all the Maslov indices from the
infinite dimensional setting  to a finite dimensional one.  
This would lead to a major simplification in  computing the spectral flow;
for example,  
one would be able to prove Theorem \ref{zorro}
by just stretching the neck of $T$ and reducing to finite dimension.

The next theorem shows that this is not the case because,
as suggested by Nicolaescu in \cite{Nicolaescu}, 
there may exist paths of
Dirac operators on a manifold with boundary for which the corresponding
paths of adiabatic limits of the Cauchy data spaces are not continuous.
Corollary \ref{adiabofX} and Proposition \ref{hrighgh} 
provide a specific example of this phenomenon, confirming Nicolaescu's
prediction.

\begin{theorem} \label{counterexample}
{\sl Let $A_t$, $0\leq t\leq \de $ be the path of connections on $Z$
specified in  Section \ref{defnofAt}. The path of operators
    $D_{A_t}, t\in[0,\de ]$ is a continuous (even analytic)  path of
    formally self-adjoint operators for which the adiabatic limits
    of the Cauchy data spaces are not continuous in $t$ at $t=0$.}
    \end{theorem}
    \noindent {\bf Proof}\qua We use $\La^R_Z(t)$ to denote the
    Lagrangian $\La_{Z^R, A_t}$. Corollary \ref{adiabofX} shows that the
    adiabatic limit of the Cauchy data spaces $\La^R_Z(t)$ is
    $P^+_{a_t}$ when $0<t\leq \de $ and $L_Z\oplus P^+_{\th}$
    when $t=0$. Since $K^+_0$ is transverse to $L_Z$, the adiabatic
    limits are not continuous in $t$ at $t=0$, ie, $$\lim_{t\to 0^+}
    \left( \lim_{R\to \infty}\La^R_{Z}(t) \right) =\lim_{t\to 0^+}
    P^+_{a_t} =K^+_0\oplus P^+_{\th} \ne L_Z \oplus P^+_{\th}
    =\lim_{R\to \infty}\La^{R}_Z(0).\eqno{\qed}$$

\subsection{Splitting the spectral flow} \label{zorrostate}

We now state the main result of this section, a splitting
    formula for the spectral flow $SF(A_t;X)$ of the family $D_{A_t}$
    when $X$ is decomposed as $X=Y \cup_T Z$.  We will use
    the machinery developed in \cite{kirk-daniel}. The technique of
    that article is perfectly suited to the calculation needed here. In
    particular, Theorem \ref{zorro} expresses the spectral flow of
    the odd signature operator on $X$ from the trivial connection
    in terms of the spectral flow on $Y$ and $Z$ between {\it
    nontrivial} connections. This greatly reduces the complexity of
    the calculation of spectral flow on the pieces.

In order to keep the notation under control, we make the following
    definitions. Given a path $A_t$ of connections on $X$ satisfying
    conditions 1--3 of Subsection \ref{defnofAt}, define the three
    paths $\xi, \eta,$ and $\si$ in $\RR^2$ with the property that
    $\xi \cdot \eta = (m_t,n_t)$ (here $\cdot$ denotes the composition
    of paths):\begin{enumerate}   \item $\xi$ is the straight
    line from $(0,0)$ to $(\de ,0)$.  \item $\eta$ is the
    remainder of $(m_t,n_t)$, ie, it is the path from $(\de ,0)$ to
    $(m_1,n_1)$ given by $(m_t,n_t)$ for $\de  \leq t \leq 1.$
      \item $\si$ is the small quarter circle centered at the
    origin from $(\de ,0)$ to $(0,\de )$. Thus $\si_t=
    (\de  \cos({\frac{t\pi} {2}}),\de  \sin({\frac{t\pi}{2}})).$
    \end{enumerate}

 \begin{figure}[ht!]\small
\begin{center} 
\includegraphics[scale=.4]{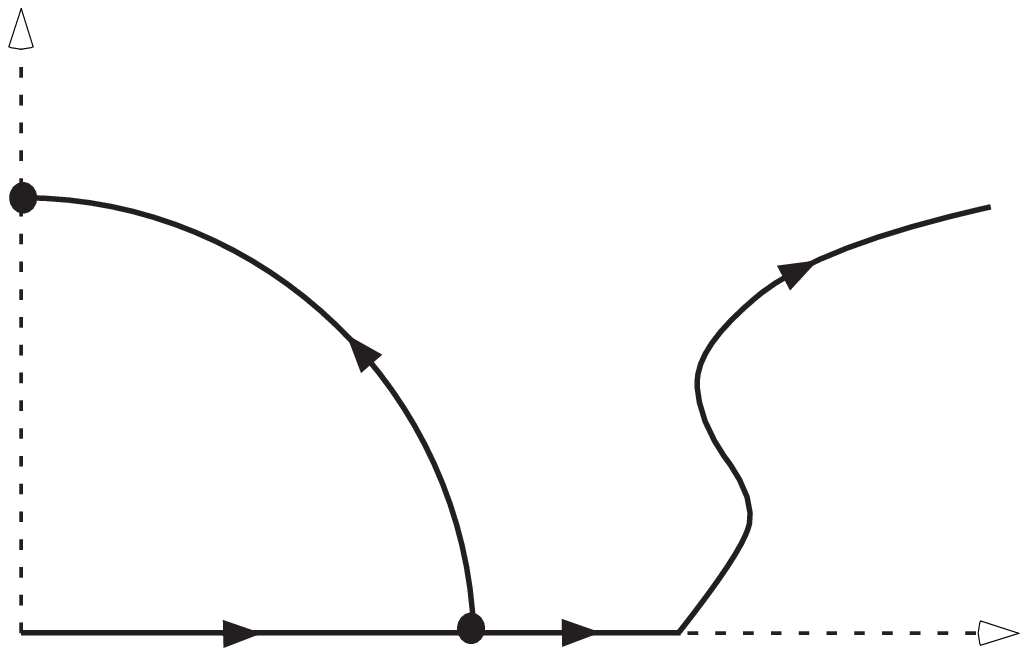}
\end{center} 
\vskip-2.25cm \hskip 3.6cm$(0,\delta)$ \hskip4cm$(m_t,n_t)$ 
 \vskip.2cm \hskip 6.1cm $\si$
\vskip .35cm \hskip 5.3cm $\xi$ \hskip 1cm $\eta$
 \vskip .2cm \hskip 6cm $(\de,0)$
 \vskip .1cm
\caption{The paths $\xi, \eta,$ and $\si$} 
\label{nearby} 
\end{figure}
 
We have paths of connections $A_\xi$ and $A_\eta$ on $X$ associated to
    $\xi$ and $\eta$. Here, $A_\xi$ is the path of connections on
    $X$ given by $A_t$ for $0 \leq t \leq \de $, and $A_\eta$
    is the path of connections on $X$ given by $A_t$ for $\de 
    \leq t \leq 1.$ In addition, using the construction of Subsection
    \ref{famofconn}, we can associate to $\si$ a path of connections
    $A_\si$ on $Y$ using the formula
$$A_{\si}(t) =    - q(r) \de  \cos(t) i dx - \de  \sin(t) i dy, 
\quad t \in [0,\tfrac{\pi}{2}].$$

\begin{theorem} \label{zorro}  {\sl Given a path $A_t$ of connections
    satisfying conditions 1--3 of Subsection \ref{defnofAt}, consider
    the paths $\xi, \eta,$ and $\si$ defined above and the associated paths of
    connections $A_\xi(t), A_\eta(t),$ and $ A_\si(t)$. Denote by $\bar{\si}
    \cdot \eta$ the path from $(0,\de )$ to $(m_1,n_1)$ which
    traces $\si$ backwards and then follows $\eta$, and denote by
    $A_{ \bar{\si} \cdot \eta}$ the corresponding path of
    connections on $Y$.
The spectral flow of $D_{A_t}$ on $X$ splits according to the
    decomposition $X=Y \cup_T Z$ as\begin{equation} \label{zorroed}
    SF(A_t;X) = SF(A_{\bar{\si} \cdot \eta}(t);Y;P^+) + SF(A_\eta(t); Z;
    P^-)-2.\end{equation}}\end{theorem}

The proof of Theorem \ref{zorro} is somewhat difficult and   has been
    relegated to the next subsection. The impatient reader is invited to
    skip ahead.

Section 4 contains a general computation of spectral flow
    on the solid torus. Regarding the other term,
    there are effective methods for computing
    the spectral flow on the knot complement
    when the
    restriction of $A_t$ to $Z$ is flat for all $t$  (see 
\cite{Fine-klassen-kirk, JDGpaper, Kirk-klassen-memoirs, 
kirk-klassen-massey prod, KKR}). For example, the main result of
    \cite{Fine-klassen-kirk} shows that after a homotopy of $A_\eta(t)$
    rel endpoints, one can assume that the paths $m_t$ and $n_t$
    are piecewise analytic. The results of \cite{Kirk-klassen-memoirs},
    combined with those  of \cite{kirk-klassen-massey prod}, can then be
    used to determine $SF(A_{\eta}(t);Z;P^-)$. The essential point
    is that the spectral flow along a path of flat connections on
    $Z$ is a homotopy invariant  calculable in terms of Massey products
    on the twisted cohomology of $Z$.

\subsection{Proof of Theorem \ref{zorro}} Applying Theorem \ref{nic1}
    shows that the spectral flow is given by the Maslov index, ie,
    that $$SF(A_t;X) = \Mas(\La_{Y}(t),\La_{Z}(t)).$$ Since the
    Maslov index is additive with respect to composition of paths
    and is invariant under homotopy rel endpoints, we prove
    (\ref{zorroed}) by decomposing $ \La_{Y}(t)$ and $\La_{Z}(t)$
    into 14 paths. That is, we define paths $M_i$ and $N_i$ of
    Lagrangians for $i =1,\ldots, 14$ so that $ \La_{Y}(t)$ and
    $\La_{Z}(t)$ are homotopic to the composite paths $M_{1} \cdots
    M_{14}$ and $N_1  \cdots   N_{14},$ respectively. We will then
    use the results of the previous section to identify $\Mas(M_i,N_i)$
    for $i=1,\ldots, 14$. The situation is not as difficult as it
    first appears, as most of the terms vanish. Nevertheless, introducing
    all the terms helps separate the contributions of $Y$ and $Z$
    to the spectral flow.

Let $a_\xi, a_\eta$ and $a_\si$ denote the paths of connections on $T$
    obtained by restricting $A_\xi, A_\eta$ and $A_\si.$ In order
    to define $M_i$ and $N_i,$ we need to choose a path $\cL_t$ of
  finite-dimensional Lagrangians in $\ker S_\th$ with the property
    that $\cL_0 = L_Z$ and $\cL_1 = K^+_0$.  A specific path $\cL_t$ will be given
    later, but it should be emphasized that the end
    result is independent of that particular choice.

We are ready to define the 14 paths $(M_i,N_i)$ of   pairs of
    infinite-dimensional Lagrangians. In each case Lemma 
\ref{bcok} shows these to be Fredholm pairs, so that their Maslov
indices are defined. 

\begin{enumerate}   \item[1.] Let $M_1 $ be the constant path at the
    Lagrangian $\La_{Y}(0)$ and $N_1 $ be the path which stretches
    $\La^R_Z$ to its adiabatic limit. Thus, using Corollary
    \ref{adiabofX}, we have $$N_1(t)=\begin{cases}
    \La^{1/(1-t)}_{Z}(0)&\text{if $0 \leq t<1$},\\ L_Z\oplus
P^+_{\th}&\text{if
    $t=1$}.\end{cases}$$  Theorem \ref{adiab} shows that $N_1$ is
continuous and      it follows from Lemma \ref{bcok} that $(M_1(t),
N_1(t)$ form a Fredholm pair for all $t$.

Since $\La_{Y}(0) \cap \La^R_{Z}(0) \cong H^{0+1}(X;\CC^2) $ is
independent of the length of the collar 
    $R$, it follows that $\dim (M_1(t)\cap N_1(t))=2 $ for $0\leq
    t<1$. At $t=1$, we have $$M_1(1)\cap N_1(1)=\La_{Y} \cap (L_Z\oplus
    P^+_{\th})=L_Y\cap L_Z$$ by Proposition \ref{exact}, since
    $H^1(Y,T;\CC^2)=0$. Since $\dim (L_Y\cap L_Z) =2$, it follows
    by Proposition \ref{masprop} that $\Mas(M_1,N_1)=0.$  

  \item[2.] Let $M_2 $ be the constant path at the Lagrangian
    $\La_{Y}(0)$. Let $N_2(t)=\cL_t \oplus P^+_{\th}.$ We claim 
    that $\Mas(M_2 ,N_2)=\Mas(L_Y, \cL_t).$

To see this, notice that $M_2$ is homotopic rel endpoints to the
composite of 3 paths, the first stretches $\La_Y(0)$ to its
adiabatic
    limit $P_\th^-\oplus L_Y$, the second is the constant path at
    $P_\th^-\oplus L_Y$, and the third is the reverse of the first,
    starting at the adiabatic limit $P_\th^-\oplus L_Y$ and returning to
    $\La_Y(0)$.
 
The path $N_2$ is homotopic rel endpoints to the composite of 3 paths,
    the first is constant at $\cL_0\oplus P_\th^+$, the second is
    $\cL_t\oplus P_\th^+$, and the third is constant at $\cL_1\oplus
    P_\th^+$.

Using homotopy invariance and additivity of the Maslov index,
    we can write $\Mas(M_2,N_2)$ as a sum of three terms. The first term is
    zero since $\La_Y^R(0)\cap (\cL_0\oplus P_\th^+)$ has dimension
    equal to  $\dim (L_Y\cap \cL_0)$ for all $R$ by Proposition
    \ref{exact}, and this also equals the dimension of
    $$ (\lim_{R\to\infty}\La^R_Y(0))\cap (\cL_0\oplus P_\th^+)=
    (P_\th^-\oplus L_Y)\cap (\cL_0\oplus P_\th^+).$$ Since the dimension
    of the intersections is constant, the Maslov index vanishes.
    Similarly the third term is zero. This leaves the second term,
    which equals $$\Mas(P_\th^-\oplus L_Y,\cL_t\oplus
    P_\th^+)=\Mas(L_Y,\cL_t).$$

  \item[3.] Let $M_3 $ be the path $\La_{Y}(t)$ for $0\leq t \leq
    \de  \ $ (this is the path of Lagrangians associated to $A_\xi$
    on $Y$). Let $$N_3(t)=\begin{cases} K^+_0 \oplus P^+_{\th}&\hbox{
    if } t=0 \\ P^+_{a_\xi(t)}&\hbox{ if } 0 < t \leq 1.\\\end{cases}$$
    That $N_3$ is continuous in $t$ was shown in the previous
    subsection.  

  \item[4.] Let $M_4$ be the path $\La_{Y,A_\si(t)}$ and $N_4$ the
    path $P^+_{a_\si(t)}. $\end{enumerate}

\begin{lemma} \label{MasIndComparison} $\Mas(M_3 \cdot M_4 ,N_3 \cdot
    N_4)= \Mas(L_Y,K^+_{t\pi/2})$.
    \end{lemma}

\noindent {\bf Proof}\qua Let $\zeta$ be the vertical line from $(0,0)$ to
    $(0,\de )$ and observe that the path $\xi \cdot \si$ is
    homotopic to $\zeta$.  Denote by
$A_\zeta(t)$ the associated path of
    flat connections on $Y$ with connection 1--form given by 
    $-t \de 
    \, idy$ (this is just the path $A_{0,t \de })$. 
    Then $M_3 \cdot M_4$ is homotopic rel endpoints to
    $M_3' \cdot M_4',$ where $M_3'$  is the constant path $\La_Y(0)$ and
    $M_4'(t) = \La_{Y,A_\zeta(t)}$. Similarly, $N_3 \cdot N_4$ is
    homotopic to $N_3' \cdot N_4'$, where $$N_3'(t) =
    K^+_{ {t\pi}/{2}} \oplus  P^+_{\th},$$ 
    $$N_4'(t) =
    \begin{cases} K^+_{\pi/2} \oplus P^+_{\th} & \text{for
    $t=0$} \\ P^+_{a_\zeta(t)} & \text{for $0<t\leq 1$,}\end{cases}$$
    and $a_\zeta(t)$ denotes the restriction of $A_\zeta(t)$ to
    $T.$

Decomposing $M_3'$ and $N_3'$ further into three paths  as in step 2
    (the proof that Mas$(M_2,N_2)=$Mas$(L_Y,\cL_t)$), we see that
    $\Mas(M_3',N_3')=\Mas(L_Y,K^+_{t\pi/2})$.

Next, Proposition \ref{exact} together with the cohomology computation
    of Equation (\ref{cohomolofW}) shows that $M_4'(0)\cap N_4'(0)$
    is isomorphic to $L_Y\cap K^+_{\pi/2}$, but Corollary \ref{trasnveree}
     shows that the latter intersection is zero.
    Another application of Proposition \ref{exact} together with
    Equation (\ref{cohomolofW}) shows that
    $M_4'(t)\cap N_4'(t)=0$ for positive $t$. Hence $M_4(t)$ and
    $N_4(t)$ are transverse for all $t$ so that Mas$(M_4',N_4')=0$.
    The proof now follows from additivity of the Maslov index under
    composition of paths.\endproof

\medskip\begin{enumerate}   \item[5.] Let $(M_5, N_5)$ be $(M_4,
    N_4)$ run backwards, so $M_5(t)= \La_{Y,A_{\bar{\si}}(t)}$ and
    $N_5(t)=P^+_{a_{\bar{\si}}(t)}$.  

  \item[6.] Let $M_6(t)=\La_{Y,A_\eta(t)}$ and $N_6(t) =
    P^+_{a_\eta(t)}$.\end{enumerate}

\medskip \noindent Theorem \ref{nic2} shows that
    \begin{equation}\label{sfonst} \Mas(M_5 \cdot M_6,N_5 \cdot
    N_6) = SF(A_ {\bar{\si} \cdot \eta}(t);Y;P^+),\end{equation}
    the advantage being that now both endpoints of $ A_ {\bar{\si}
    \cdot \eta}$ refer to {\it nontrivial} flat connections on $Y$.
    In the next section we will explicitly calculate this integer
    in terms of homotopy invariants of the path $\bar{\si} \cdot
    \eta$.

\medskip

\begin{enumerate}   \item[7.] Let $M_7 $ be the path obtained by
    stretching $\La^R_{Y}(1)$ to its adiabatic limit. Since $a_1,$
    the restriction of $A_1$ to $T^2,$ is a nontrivial flat connection,
    ${ \lim_{R\to \infty}} \La^R_{Y}(1) = P^-_{a_1}$. This follows
    from Corollary \ref{adiabofX} applied to $Y$, or directly by
    combining Theorem \ref{adiab}, Proposition \ref{exact} and
     Equation (\ref{cohomolofW}).

Let $N_7 $ be the constant path $P^+_{a_1}$. An argument similar to the
    one used in step 1 shows that $\Mas(M_7,N_7)=0.$  

  \item[8.] Let $M_8(t)=P^-_{a_\eta (1-t)}$ and $N_8(t) =P^+_{a_\eta
    (1-t)}$ (this is just $N_6$ run backwards). Observe that since
    $M_8(t)$ and $N_8(t)$ are transverse for all $t$, $\Mas(M_8,N_8)=0.$

  \item[9.] Let $$M_9(t)=\begin{cases} P^-_{a_\xi (1-t)}&\hbox{ if }
    0\leq t<1\\ K^-_0\oplus P^-_{\th}&\hbox{ if } t=1\end{cases}$$
    and $$N_9(t) =\begin{cases} P^+_{a_{\xi}(1-t)} & \text{for
    $t<1$} \\ K^+_0\oplus P^+_\th & \text{for $t=1$}.\end{cases}$$
    Now $N_9$ is just $N_3$ run backwards, and it is not difficult
    to see that $M_9(t)$ and $N_9(t)$ are transverse for all $t$,
    hence $\Mas(N_9,M_9)=0.$  

  \item[10.] Let $M_{10} $ be the constant path at $K^-_0 \oplus
    P^-_{\th}$ and let $N_{10} $ be $N_2$ run backwards, ie, $N_2(t) =
    \cL_{1-t} \oplus P^+_\th.$ Thus, $\Mas(M_{10},N_{10}) =
    \Mas(K^-_0,\cL_{1-t})$.  

  \item[11.] Let $M_{11} $ be the constant path at $K^-_0 \oplus
    P^-_{\th}$ and $N_{11} $ be $N_1$ run backwards, ie,
    $$N_{11}(t)=\begin{cases} L_Z\oplus P^+_{\th} &\text{if $t=0$}
    \\ \La^{1/t}_{Z}(0)&\text{if $t>0$.}\end{cases}$$  
     Propositions \ref{exact} and \ref{hrighgh} and
    Corollary \ref{adiabofX} show that $M_{11}(t)$ is transverse to
    $N_{11}(t)$ for all $t$, hence $\Mas(M_{11},N_{11})=0.$ 

  \item[12.] Let $M_{12} $ be $M_{9} $ run backwards, ie,
    $$M_{12}(t)=\begin{cases} K^-_0\oplus P^-_{\th}&\hbox{ if } t=0
    \\ P^-_{a_\xi(t)}&\hbox{ if } 0<t\leq 1.\end{cases}$$ Let
    $N_{12}(t)=\La_{Z,A_\xi(t)}$. Since the restriction of $A_\xi(t)$ to
    $Z$ is flat, Proposition \ref{exact} shows that $M_{12}(t)$ is
    transverse to $N_{12}(t)$ for all $t$. Hence
    $\Mas(M_{12},N_{12})=0.$  

  \item[13.] Let $M_{13} (t)= P^-_{a_\eta(t)}$ (ie, $M_8$ run
    backwards) and let $N_{13}(t)=\La_{Z,A_\eta(t)}$. Theorem \ref{nic2}
    then implies that $$\Mas(M_{13},N_{13})=SF(A_\eta(t);Z;P^-),$$
    the spectral flow on $Z$.  

  \item[14.] Let $M_{14}$ be $M_7$ run in reverse and $N_{14} $ the
    constant path at $\La_{Z,A_1}$. An argument like the one in
    step 1 (but simpler since $\ker S_{a_1}=0$) shows that
    $M_{14}(t)\cap N_{14}(t) \cong H^{0+1}(X;\CC^2_{A_1})$  for all
    $t$. This implies that $\Mas(M_{14},N_{14})=0.$

\end{enumerate}

\medskip \noindent We leave it to the reader to verify that the terminal
    points of $M_i$ and $N_i$ agree with the initial points of $M_{i+1}$
    and $N_{i+1}$ for $i=1, \ldots ,13,$ and that $M_1 \cdots  M_{14}$
    and $N_1   \cdots   N_{14}$ are homotopic rel endpoints to
    $\La_{Y}(t) $ and $\La_{Z}(t)$, respectively. Thus $$SF(A_t;M
    )=\Mas( M_1 \cdots  M_{14}, N_1   \cdots   N_{14}) = \sum_{i=1}^{14}
    \Mas(M_i,N_i).$$ The arguments above show that $\Mas(M_i,N_i) =
    0 $ for $i=1,7,8,9,11,12,$ and $14.$ Moreover, by Equation
    (\ref{sfonst}) and step 13, we see that\begin{eqnarray*} \Mas(
    M_5 \cdot  M_{6}, N_5 \cdot N_{6})  &=& SF(A_{\bar{\si} \cdot
    \eta}(t);Y;P^+), \hbox{ and}\\ \Mas(M_{13},N_{13}) &=&
    SF(A_\eta(t);Z;P^-).\end{eqnarray*}

To finish the proof of Theorem \ref{zorro}, it remains to show that
the sum of the remaining terms
 $$\Mas(M_2,N_{2})
      + \Mas(M_{3}\cdot M_4,N_{3}\cdot N_4)+
    \Mas(M_{10}, N_{10})$$
    equals $-2$.  By Step 2, Lemma \ref{MasIndComparison}, and Step 10,
    these summands equal $\Mas(L_Y,\cL_t)$, $\Mas(L_Y, K^+_{t\pi/2})$ and
    $\Mas(K^-_0,\cL_{1-t})$, respectively.

    Define  the path $\cL_{t}$
    to be
      \begin{eqnarray} \cL_t &=& \sspan \{ (1, \ (1-t) idx + t
    idy,0)\otimes e_1, \ (1, \ (t-1) idx - t idy,0)\otimes e_2,
    \nonumber \\  && \quad \quad \ (1-t, \ -i dx, -\ t dxdy) \otimes
    e_1, \ (1-t, \ i dx, -\ t dxdy) \otimes e_2\}. \label{defnofcL}
    \end{eqnarray}

\begin{lemma} \label{MasIndComp} {\sl For the path $\cL_{t}$ in Equation
(\ref{defnofcL}),
\begin{enumerate}
\item[\rm(i)]   $\Mas(L_Y,\cL_t) =0.$
    \item[\rm(ii)] $\Mas(L_Y,K^+_{t\pi/2}) = -2.$
\item[\rm(iii)] $\Mas(K^-_0, \cL_{1-t})=0.$\end{enumerate}}
\end{lemma}

\noindent
\noindent {\bf Proof}\qua  Proposition \ref{hrighgh}
and  Equation (\ref{defnofcL}) imply that  $K^-_0$ and $\cL_t$ are
transverse for $0 \leq t \leq 1$. Hence
$\Mas( K^-_0, \cL_{1-t}) = 0$. This
proves claim (iii).

Next consider claim (ii).  
Corollary \ref{trasnveree} implies that
 $\dim(L_Y \cap K^+_{t\pi/2}) =
     0$ for   $0<t\leq1 $. An exercise in linear algebra shows that, for
     small $s>0$, $\dim (e^{sJ}L_Y\cap K^{+}_{t\pi/2})=0$ unless
     $\tan (t\pi/2)=\tan (2s)$, and for this $t$
     (which is positive and close to 0) the
     intersection has dimension 2.  
 Apply Proposition \ref{princ} with $s(t) = t \pi /4$ to conclude
 that $\Mas(L_Y,K^+_{t\pi/2}) = -2.$ 
      
     Finally, consider claim (i).  It is easily verified that
 $$ \dim ( L_Y \cap \cL_t) =\begin{cases} 2 &
    \text{if $t=0,1$,} \\ 0 & \text{if $0 < t < 1$.}\end{cases}$$ 
    A direct calculation shows further that $e^{s J} L_Y
    \cap \cL_t \neq 0$ if and only if 
\begin{equation} \label{gimma}
0 = (1+\sin 2s ) t^2 +
    (1-\sin 2s )\, t + \sin 2s,
    \end{equation}
in which case  $ \dim( e^{s J} L_Y \cap \cL_t) = 2.$
We will apply Proposition
    \ref{princ} to the intersection of $L_Y $ and $\cL_t$ at $t=0$
    and the `reversed' result to the intersection at $t=1$
     (cf.~the remark immediately 
    following the proof of Proposition
    \ref{princ}).  
    The solutions $t =t(s)$  to (\ref{gimma}) are the two functions
     $$t_{\pm}(s) = \frac{1}{2}
    \pm \frac{1}{2} \sqrt{\frac{1+3\sin 2s}{1-\sin 2s}}.$$  
    Notice that
    $t_+(0)=1$ and $t'_+(0) > 0$ and $t_-(0)=0$ and $t_-'(0) < 0.$
    Apply Proposition \ref{princ} to $s_-(t)$ at $t=0$,
    and also apply its reversed result to $s_+(t)$ at $t=1$,
    where $s_{\pm}$ denote the inverse functions of $t_\pm$. It follows
    that
    $
    \Mas(L_Y, \cL_t; \ 0 \leq t \leq \de)=0 $ and $
 \Mas(L_Y, \cL_t; \ 1-\de \leq t \leq 1) = 0.$ \endproof

\section{Spectral flow on the solid torus} \label{sfsection}

In this section, we  carry out a detailed analysis of connections on the
    solid torus $Y$ and show how to compute the spectral flow between
    two nontrivial flat connections on $Y$.  We reduce the
    computation to an algebraic problem by explicitly constructing
    the Cayley graph associated to the gauge group using paths
    of connections.

\subsection{An SU(2) gauge group for connections on Y in normal form
on T}
\label{relative gauge group}
We begin by specifying certain groups of gauge transformations which
leave invariant the spaces of connections on $T$ and $Y$ which are in
normal form (on $T$ or along the collar).  We will
 identify $SU(2)$ with the 3--sphere $S^{3}$
of unit quaternions, and we identify the diagonal subgroup with
$S^{1}\subset S^{3}$.

Define $\tal, \tbe \co T\to S^{1}$ by the formulas
$$\tal (e^{ix}, e^{iy})=e^{ix}, \ \ \tbe (e^{ix},e^{iy})=e^{iy}.$$
Let $H$ be the abelian group generated by $\tal$ and $\tbe$, which
 act on $\cA_{\rm nf}(T)$ by
$$\tal \cdot a_{m,n}=a_{m+1,n},\ \ \tbe \cdot a_{m,n}=a_{m,n+1}.$$

Let $\cA_{\rm nf}(Y)$ denote the space of connections on $Y$ which are in normal
form on the collar (cf.~Definition \ref{Anf}),
$$\cA_{\rm nf}(Y)=\{ A\in \Om^{1}_{Y}\otimes su(2) \mid  A|_{[-1,0]\times T} \text{
is cylindrical and in normal form}\}.$$
  Let $r\co \cA_{\rm nf}(Y) \to \cA_{\rm nf}(T)$
denote the restriction map.  We define the gauge group
$$\cG_{\rm nf} =\{ \text{smooth maps } g\co Y\to S^{3} \mid g|_{[-1,0]\times T}=\pi^{*}h
\text{ for some } h\in H\},$$
where $\pi\co [-1,0]\times T\to T$ is projection.  It is clear that, for
$g\in \cG_{\rm nf} $ with $g|_{T}=h$, we have the commutative diagram
$$\begin{CD}
    \cA_{\rm nf}(Y) @>{g}>> \cA_{\rm nf}(Y) \\ @V{r}VV  @VV{r}V \\
    \cA_{\rm nf}(T) @>>{h}> \cA_{\rm nf}(T).\\
    \end{CD}$$
To clarify certain arguments about homotopy classes of paths, it is
convenient to replace the map $r\co \cA_{\rm nf}(Y) \to \cA_{\rm nf}(T)$ with the map
$Q\co \cA_{\rm nf}(Y) \to \RR^{2}$ defined by
$$
Q(A)=(m,n) \text{ where } A|_{T}=a_{m,n}.$$

The identity component $\cG_{\rm nf}^0 \subset \cG_{\rm nf} $ is a normal subgroup,
and we denote the quotient by $G=\cG_{\rm nf} /\cG_{\rm nf}^0$.

Recalling the orientation on $Y$ from Section \ref{torus decomp} and
using the orientation of $S^{3}$ given by the basis $\{ i,j,k \}$
for $T_{1}S^{3}$, we note that each $g\in \cG$ has a well-defined
degree, since $H_{3}(S^{3},S^{1}; \ZZ) = \ZZ$, and this degree
remains well-defined on $G$.

\begin{lemma}\label{hpic} {\sl Let $g,g'\in \cG_{\rm nf}$. Then $g$ is
homotopic to $ g'$ (ie, they represent the same element of $G$) if
and only if
    $(g|_T)=(g'|_T)$  and $\deg(g)=\deg(h)$.}
    \end{lemma}

\noindent {\bf Proof}\qua This is a simple application of obstruction theory that we
    leave to the reader.\endproof

It follows from Lemma \ref{hpic} that the restriction map descends to
a map $P\co G\to H$ which is onto, since
    $\pi_1(S^3)=\pi_2(S^3)=0$. Set $K=\ker P \cong\ZZ$, where the
    last isomorphism is given by the degree.

\begin{lemma} {\sl The kernel of $P\co G\to H$ is central.}\end{lemma}
\noindent {\bf Proof}\qua
    Suppose $k\in K$ and $g\in G$. After a homotopy, 
    we may assume
    that there is a 3--ball $B^3$ contained in the interior of $Y$
    such that $k|_{Y-B^3}=1$ and $g|_{B^3}=1$. It follows directly
    from this that $gk=kg$.
    \endproof

Using the cutoff function  $q(r)$ from Equation (\ref{eqnofbmn}),
we define $\al,\be,\ga \in G$ as follows (we make the definitions in $\cG_{\rm nf}$
    but they should be reduced mod $\cG^0_{\rm nf}$):
    \begin{enumerate}
    \item[(i)] $\al(re^{ix},w)=q(r)e^{ix}+\sqrt{1-(q(r))^2} j.$  
    \item[(ii)] $\be(re^{ix},w)=w$,  
     \item[(iii)] $\ga(z,w)=\hbox{ a generator of
    }K\hbox{ with }\deg(\ga)=1.$\end{enumerate}

It will be useful to denote by $\bar \alpha$ the map
$$\bar\al(re^{ix},w)=re^{ix}+\sqrt{1-r^2} j,$$
which is not in $\cG_{\rm nf}$ but is homotopic rel boundary to
$\al$ and has a simpler formula. Using  $\bar \alpha$ will simplify
the computation of
  degrees of maps involving $\al$.  Observe that
$$\begin{array}{lll} P(\al)=\tal&\hbox{and}&\deg(\al)=0\\
    P(\be)=\tbe&\hbox{and}&\deg(\be)=0\\
    P(\ga)=1&\hbox{and}&\deg(\ga)=1.\end{array}$$

Now $[\al,\ga]=[\be,\ga]=1$, hence $G$ is a central extension of $H$ by
    $K$: $$0 \lto K \lto G \lto H \lto 0.$$ Such extensions are
    classified by elements of $H^1(H;\ZZ),$ and to determine the cocycle
    corresponding to our extension, we just need to calculate which
    element of $K$ is represented by the map $[\al,\be]$. This amounts
    to calculating the degree of this map.

\begin{lemma}\label{commutatior} $[\al,\be]=\ga^{-2}$.\end{lemma}
    \noindent {\bf Proof}\qua Set $h = [\al,\be].$ Clearly, $h \in\ker(P)$, so we
    just need to calculate its degree. It is sufficient to compute the
    degree of $\bar h=[\bar \al, \be]$, since it is homotopic to $h$
    rel boundary.  Using the coordinates $(re^{ix},
    e^{iy})$ for $Y$ and writing quaternions as $A+jB$ for $A, B
    \in \CC,$ we compute that $$\bar h(re^{ix},e^{iy}) = r^2+(1-r^2)e^{-2iy}
    + j r\sqrt{1-r^2} e^{-ix}(1-e^{-2iy}).$$ To determine the degree of
    $\bar h$, consider the value $k\in S^3$ which we will prove is a regular
    value. Solving $\bar h(re^{ix},e^{iy})=k$ yields the two solutions
    $r={\frac{1}{\sqrt{2}}}$, $x={\frac{\pi}{2}}$,
$y={\frac{\pi}{2}}\hbox{
    or }{\frac{3\pi}{2}}$. Applying the differential $d\bar h$ to the
oriented
    basis $\{ \frac{\partial}{\partial r} ,{\frac{\partial}{\partial
    x}},{\frac{\partial}{\partial y}}\}$ for the tangent space of $Y$
    and then translating back to $T_1(S^3)$ by right multiplying by
    $-k$ gives the basis $\{-2\sqrt{2}k,i,i+j\}$ of $S^3$, which is
    negatively oriented compared to $\{i,j,k\}$. Since the computation
    gives this answer for both inverse images of $k$, it follows
    that $\deg(h)=-2$, which proves the claim.\endproof

We have now established the structure of $G$. Every element $g\in G$ can
    be expressed uniquely as $g=\al^a\be^b\ga^c$ where $a,b,c \in\ZZ$.
    Furthermore, with respect to this {\em normal form}, multiplication can
    be computed as follows: $$
    (\al^{a_1}\be^{b_1}\ga^{c_1})(\al^{a_2}\be^{b_2}\ga^{c_2}) =
    \al^{a_1+a_2}\be^{b_1+b_2}\ga^{2 b_1a_2+c_1+c_2}$$

The next result determines the degree of any element in normal form.

\begin{theorem}\label{oiewr} $\deg(\al^a\be^b \ga^c)=c-ab$.
    \end{theorem}
\noindent {\bf Proof}\qua We begin by computing the degree of
    $\al^a\be^b$. Let $f_{a}\co D^2\times S^1\to S^3$ be the map
    $$f_{a}(re^{ix}, e^{iy})= \al(r^{|a|} e^{iax}, e^{iby})= r^{|a|}
    e^{iax} +\sqrt{1- r ^{2|a|}}j .$$ Then $f_a$ is homotopic rel
    boundary to $\bar \al^a$ using Lemma \ref{hpic} since they agree on
    the boundary and both have degree $0$ (they factor through the
    projection to $D^2$).

The degree of $\al^a\cdot \be^b$
    equals the degree of $f_a\cdot \be^{b}$, since $\al^a$ is homotopic
to $f_a$.

But $ f_a\cdot \be^{b}$ factors as the composite of the map
    \begin{eqnarray*} &D^2\times S^1\lto D^2\times S^1 & \\ &(re^{ix},
    e^{iy})\mapsto (r^{|a |}e^{iax}, e^{iby})&\end{eqnarray*} and
    the map\begin{eqnarray*} &D^2\times S^1\lto S^3 & \\
    &(z,w)\mapsto\bar \al(z,w)\be(z,w).&\end{eqnarray*} The first map
    is a product of a branched cover of degree $a$ and a cover of
    degree $b$ and so has degree $ab$. The second restricts to a
    homeomorphism of the interior of the solid torus with $S^3-S^1$
    which can easily be computed to have degree $-1$.
    Thus $\al^a\be^b$ has degree $-ab$.

To finish proving the theorem, we need to calculate the effect of
    multiplying by $\ga$. For any $g\in \cG$, we can
    arrange by homotopy that $\ga$ is supported in
    a small 3-ball while $g$ is constant in the same 3-ball. It is
    then clear that for all $g\in G$, $\deg(g\ga)=\deg(g)+1$.
    \endproof

\subsection{The ${\CC^2}$--spectral flow on Y} \label{C2sfonW}

Suppose that $A_t \in \cA_{\rm nf}(Y)$ is a path between the flat connections
    $A_0$ and $A_1$ on $Y$. We will present a technique for computing
    $SF(A_t;Y; P^+),$ the spectral flow of the odd signature operator
    $$D_{A_t}\co \Om^{0+1}_Y\otimes\CC^2\lto \Om^{0+1}_Y\otimes\CC^2$$
    on $Y$ with $P^+$ boundary conditions. We assume that for all
    $t$, $Q(A_t) \in \RR^2 - \ZZ^2$. This implies that $P^+_{a_t}$
    varies continuously in $t,$ where $a_t$ denotes the restriction
    of $A_t$ to $T$ \cite{KKillinoisJ}. Moreover the exact sequence
    in
   Proposition \ref{exact} shows that the kernels of $D_{A_0}$ and $D_{A_1}$
    with $P^+$ boundary conditions are zero.

\begin{lemma}\label{splittingthm} {\sl  Let $Y_1$ and $Y_2$ be solid tori, and
    let $X=Y_1\cup Y_2$ be the  lens space obtained by gluing $\partial Y_1$
    to $\partial Y_2$ using an orientation reversing isometry
    $h\co \partial Y_1\to \partial Y_2$. Let $A_t$ be a path in
    $\cA_{\rm nf}(Y_1)$ and $B_t$ a path in $\cA_{\rm nf}(Y_2)$ so that
    $h^*(B_{t} |_{\partial Y_2})=A_{t}|_{\partial Y_1}$. Assume that
    $Q(A_t)\in\RR^2-\ZZ^2$ and that $A_0,A_1 ,B_0,B_1$ are flat.
Then $$SF(A_t\cup B_t;X)=SF( A_t;Y_1;P^+) + SF(B_t;Y_2; P^+).$$}
    \end{lemma}
    \noindent {\bf Proof}\qua   Write $T=\partial Y_1$ and let
    $a_t=A_t |_{T}$. The cohomology computation (\ref{cohomolofT}) shows
    that $H^{0+1+2}(T;\CC^2_{a_t})=0$ for all $t$. Hence $\ker
    S_{a_t}=0$ for all $t$. Also, the computation (\ref{cohomolofW})
    shows that $H^{0+1}(Y_1;\CC ^2_{A_i})=\ker D_{A_i}(P^+)=0$
    for $i=0,1$ and that $H^{0+1}(Y_1;\CC ^2_{B_i})=\ker D_{B_i}(P^+)=0$
    for $i=0,1.$
     
The lemma now follows from the splitting theorem for spectral flow of
    Bunke (Corollary 1.25 of \cite{bunke}). For a simple proof using the
    methods of this article see \cite{kirk-daniel}.\endproof

\begin{lemma}\label{homotopy} {\sl Suppose $A_t$ and $B_t$ are two paths in
    $\cA_{\rm nf}(Y)$ such that $A_i$ and $B_i$ are flat for $i=0,1$. Suppose
    further that the paths $Q(A_t)$ and $Q(B_t)$ miss
the integer
    lattice $\ZZ^2 \subset \RR^2$ for all $t \in [0,1].$ If $ A_i =
    g_i\cdot B_i $ for $i=0,1$ where $g_i\in \cG_{\rm nf}^0$ and if the
    paths $Q(A_t)$ and $Q(B_t)$ are homotopic rel
endpoints in
    $\RR^2 - \ZZ^2$, then $SF(A_t;Y; P^+)=SF(B_t;Y;P^+)$.}

\end{lemma}\noindent {\bf Proof}\qua First, note that a path of the form $g_tA$, where
    $g_t$ is a path in $\cG$, has spectral flow zero, because the
    eigenvalues are all constant. (This follows from the fact that
    the operators in the path are all conjugate.) Hence we may
    assume that $A_i=B_i$ for $i=0,1$ (if not, add a  path of
    the form $g_tA_i$ to each end of $B_t$ bringing the endpoints
    together). Now, using the fact that $\cA_{\rm nf}(Y)$ is a  bundle over $\cA_{\rm nf}(T)$
    with contractible fiber, it is easy to see that
     the homotopy between $Q(A_t)$ and $Q(B_t)$
     can be lifted to one between
    $A_t$ and $B_t$ which will, of course, avoid $Q^{-1}(\ZZ^2)$.
    Finally, homotopic paths of operators have the same spectral
    flow, proving the lemma.\endproof

Based on this lemma, we may now state precisely the question we wish to
    answer: Given a path of connections $A_t$ in $\cA_{\rm nf}(Y)$
    between two flat connections such that $Q(A_{t})$
    avoids $\ZZ^{2}\subset \RR^{2}$, how can one calculate
    $SF(A_{t};Y;P^{+})$ from $A_{0}$, $A_1$,
    and the image $Q(A_{t})$ in
    $\RR^{2}-\ZZ^{2}$?

The following lemmas serve as our basic computational tools in what follows.

\begin{lemma}\label{degree} {\sl Suppose $X$ is a closed oriented 3--manifold
    and $g\co X\to SU(2)$ is a gauge transformation. If $A_0$ is any
    $SU(2)$ connection on $X$, and $A_t$ is any path of connections
    from $A_0$ to $A_1=g\cdot A_0 = gA_0g^{-1}-dg~g^{-1}$,
    then $$SF({A_t};X)=-2\deg(g).$$}
    \end{lemma}
    \noindent {\bf Proof}\qua Recall that we are using the $(-\varepsilon,
    -\varepsilon)$ convention for computing spectral flows. The
    claim follows from a standard application of the Index Theorem.
    See for example the appendix to \cite{KKR}.\endproof

\begin{lemma}\label{internal} {\sl Let $A$ be any connection in $\cA_{\rm nf}(Y)$
    with $Q(A) \in \RR^2 - \ZZ^2$ and let $g\in\cG_{\rm nf} $ be a gauge
    transformation which is $1$ on the collar neighborhood
    of the boundary
    $T$. If $A_t$ is any path in $\cA_{\rm nf}(Y)$ from $A$ to $g\cdot A$ which
    is constant on $T$ (eg, the straight line from $A$ to $g\cdot A$),
    then $SF(A_t;Y; P^+)=-2\deg(g)$.}\end{lemma}

\noindent {\bf Proof}\qua Consider a path $B_t$ of connections on the double $D(Y)$ of
    $Y$ which is constant at $A$ on one side and is $A_t$ on the
    other side, and the gauge transformation $h$ which is $g$ on one
    side and the identity on the other. Then $B_1=hB_0$, and
    $\deg(h)=\deg(g)$. Lemma \ref{degree} shows that $SF(B_t;
    D(Y))=-2\deg(g)$. Now apply Lemma \ref{splittingthm}.\endproof

Since we are interested in paths between flat connections, we begin by
    analyzing the components of orbits of flat connections in  $\cA_{\rm nf}(Y)
    /\cG_{\rm nf}^0$.
    First, note that all the flat  connections in $\cA_{\rm nf}(Y)$
    project to $\ZZ\times\RR$ under $Q\co  \cA_{\rm nf}(Y) \to \RR^2$. Set $\tJ$
    equal to the open vertical  line segment $\tJ=\{(0,t) \mid 0 < t <
    1\}\subset \RR^{2}$.

A natural choice of gauge representatives for $Q^{-1}(\tJ)$ is the path
of connections  $J = \{-tidy \mid 0<t<1 \}\subset \cA_{\rm nf}(Y)$.
    The connection $-t idy$ is a flat connection on $Y$ whose holonomy
    sends $\mu$ to $1$ and $\la$ to $e^{2\pi i t}$.  Note that the spectral
    flow of any path $A_t$ whose image modulo $\cG_{\rm nf}^0$ lies
     in $J$ is 0, since $\ker
    D_{A_{t}}$ is constantly zero.

The set of all flat orbits  in $\cA_{\rm nf}(Y)
    /\cG_{\rm nf}^0$ not containing any gauge transformations of the trivial
    connection may be expressed as
    $\bigcup_{g\in G}(g\cdot J)$. For every nontrivial $g\in G$,
    $g\cdot J$
    is disjoint from $J$. This can be seen by considering the
    action of $P(g)$ on $J$, and using Lemma \ref{internal} above.
    The reader is encouraged to visualize the orbit of $J$ under
    $G$ as consisting of $\ZZ$ homeomorphic copies of $J$ sitting
    above each translate $(p,q)+\tJ$ in $\RR^2$, where $p$ and $q$ are
    integers.

We will now build a graph $\Ga$ with one vertex corresponding to each
    component of $G\cdot J$. Note that these vertices are also in one-to-one
    correspondence with $G $. Next, we will construct some directed
    edges with $J$ as their initial point. Actually, for specificity,
    we will think of their initial point as being $c_0=-{\tfrac{1}{2}}idy$ of
$J$.

Let $E_\al$ be the straight line path of connections from $c_0$ to
    $\al\cdot c_0$. We construct a corresponding (abstract) edge in
    $\Ga$ from $J$ to $\al J$, which we also denote by $E_\al$.
    Now for all $g\in G$, construct another edge $gE_\al$ from
    $g\cdot J$ to
    $g\al J$, which one should think of as corresponding to the
    path $gE_\al$ in $\cA_{\rm nf}(Y)$. Thus every vertex of $\Ga$ serves
    as the initial point of one $\al$--edge and the terminal point
    of another.

Next we construct a path $E_\be$ in $\cA_{\rm nf}(Y)$ from $c_0$ to $\be c_0$.
    We cannot use the straight line because its image in $\RR^{2}$
    hits the integer lattice, so
    instead we define $E_\be$ to be the path
in $\cA_{\rm nf}(Y)$ given by
  $$A_{t}= -\tfrac{1}{2}
    q(r) \cos t \ i dx - (1+ \tfrac{1}{2}\sin t) \ i dy,
    \ \ -\tfrac{\pi}{2} \leq t  \leq \tfrac{\pi}{2},$$
    where $q(r)$ is the radial bump function in Equation
    (\ref{eqnofbmn}).  Thus
    $Q(E_\be)$ is the semicircle
    $(\tfrac{1}{2}\cos t,1+\tfrac{1}{2}\sin t)$,
    $t \in[-\tfrac{\pi}{2}, \tfrac{\pi}{2}]$.   (As before, it
    is only the homotopy class of the path $E_\be$ in
    $Q^{-1}(\RR^2-\ZZ^2)$ rel endpoints that is important.)
For each $g\in G$, build an edge $gE_\be$ in $\Ga$ from $g\cdot  J$ to
    $g\be\cdot  J$ corresponding to the path of connections $gE_\be$.
 
Finally, construct a path of connections $E_\ga$ in $\cA_{\rm nf}(Y)$ from $c_0$ to
    $\ga c_0$ such that $Q(E_\ga)$ is the constant path in $\RR^2$
    at $(0,\frac{1}{2})$. A straight-line path would be acceptable
    in this case. Once again, for each $g\in G$, define an edge
    $gE_\ga$ in $\Ga$ from $g\tJ$ to $g\ga\tJ$. The resulting graph
    $\Ga$ is isomorphic to the Cayley graph of $G$, defined with respect
    to right multiplication by the generators $\{\al,\be,\ga\}$.

Notice that we have also constructed a 1--dimensional graph in $\cA_{\rm nf}(Y)$
the image of which in $\cA_{\rm nf}(Y)/\cG_{\rm nf}^0$ is  invariant under $G$;
this will provide us
    with a complete (up to homotopy and gauge transformation in  $\cG_{\rm nf}^0$)
    collection of paths of connections
    in $\cA_{\rm nf}(Y)$ connecting components of the flat connections in
    $Q^{-1}(\RR^{2}-\ZZ^{2})$.

The next step is to associate to each edge of $\Ga$ an integer, which
    will give the spectral flow of the odd signature operator $D_A$
    on the solid torus $Y$ with $P^+$ boundary conditions along the
    corresponding path of connections in $\cA_{\rm nf}(T)$. Of course, the
    integer associated to the path $gE_\al$ is independent of $g\in
    G$ since the gauge transformation $g$ induces a relation of
    conjugacy between $D_A$ and $D_{g\cdot A}$ for each $A$ in the path.
    An analogous fact holds for the edges $E_\be$ and $E_\ga$, as
    well. So we just need to find three    integers $k_\al$, $k_\be$,
    and $k_\ga$, one for each class of edges.

\begin{theorem}\label{constants} {\sl These constants have values $k_\al=2$,
    $k_\be=-2$, and $k_\ga=-2$.}\end{theorem}
    \noindent {\bf Proof}\qua The value
    of $k_\ga$ is calculated to be $-2$ in Lemma \ref{internal}, so
    we turn our attention to calculating $k_\al$ and $k_\be$.

Let $Y_1$ and $Y_2$ be two solid tori
with the same orientations.  
 For $i=1,2,$ set $T_i = \partial Y_i$ with coordinates $x_i,y_i$
 and
    let $\mu_i$, $\la_i$, $dx_i$, and $dy_i$  
    denote the loops and forms
    on $Y_i$. Glue $Y_1$ to $Y_2$
    by the homeomorphism of $T_1$ with $T_2$ which identifies
    $(e^{ix_1},e^{iy_1})$ with 
    $(e^{i(x_2+y_2)}, e^{i(2x_2+y_2)})$.
  Since this map is orientation reversing, we may give
    $Y_1\cup Y_2$ the orientation of both $Y_1$ and $Y_2$.
     Let $A_{1,t}$
    denote the path of connections on $Y_1$ corresponding to $E_\al$.
    When restricted to $T_1$, these connections are given by the straight
    line, ie, $A_{1,t} |_{T_1}=-tidx_1-\frac{1}{2}idy_1$.

We now need to construct a path of connections on $Y_2$ which is
    compatible along $T_2$ with $A_{1,t}|_T$. Pulling the connections $A_{1,t}|_{T_1}$
    back to $T_2$ by the above formula gives a path 
    $$a_t=-(t+1)idx_2
    -(t+\frac{1}{2})idy_2$$ of connections on $T_2$. 
    Under the identification
    $\cA_{\rm nf}(T_2)
    \cong \RR^2$, this is the straight line from
    $(1,\frac{1}{2})$ to $(2,\frac{3}{2})$. Now define the path $B_t$ in
    $\cA_{\rm nf}(Y_2)$ by first following the path $\al E_\be$ from $\al
    c_0$ to $\al\be c_0$, and then the path $\al\be E_\al$ from $\al\be
    c_0$ to $\al\be\al c_0$. Note that $B_t$ runs from $\al c_0$ to
    $\al\be\al c_0$, and that $Q(B_t)$ is a path which is homotopic
    rel endpoints in $\RR^2-\ZZ^2$ to the straight line from
    $(1,\frac{1}{2})$ to $(2,\frac{3}{2})$. Hence we may define a
    path $A_{2,t}$ in $\cA_{\rm nf} (Y_2)$ which is homotopic rel endpoints to
    $B_t$ in $Q^{-1}(\RR^2 - \ZZ^2) \subset \cA_{\rm nf}(Y_2)$ and has the
    property that $Q(A_{2,t})=(t+1,t+\frac{1}{2})$. 
    By Lemma \ref{homotopy}, $B_t$ and $A_{2,t}$ have
    the same spectral flow.

Consider the path of connections $A_t$ on $Y_1\cup Y_2$ 
defined to be $A_{1,t}$
    on $Y_1$ and $A_{2,t}$ on $Y_2$. Note that $A_1=gA_0$, where $g$ is
    a gauge transformation on $Y_1\cup Y_2$ equal to $\al$ on $Y_1$ and
    equal to $\al\be$ on $Y_2$. Since $\deg(\al)=0$ and $\deg(\al\be)=-1$
    by Theorem \ref{oiewr}, it follows that $\deg(g)=-1$. Hence $SF(A_t;
    Y_1\cup Y_2)=2$ by Lemma \ref{degree}. On the other hand, by Lemma
    \ref{splittingthm}, $$SF(A_t; Y_1\cup
    Y_2)=SF(A_{1,t};Y_1;P^+)+SF(A_{2,t};Y_2;P^+)=2k_\al+k_\be.$$ This gives
     the linear equation $2=2k_\al+k_\be$.

Repeating this process, we glue $Y_1$ to $Y_2$ by the homeomorphism of $T_1$
    with $T_2$ which identifies $(e^{ix_1},e^{iy_1})$ with
    $(e^{i(x_2+2y_2)}, e^{i(2x_2+3y_2)})$.  This gives another
    equation which can be used to solve for $k_\al$ and $k_\be$. Pulling
    back  the
    same path of connections $A_{1,t} |_{T_1}$ to $T_2$ using the new
    gluing map, we obtain $$a_t=-(t+1)idx_2-(2t+\tfrac{3}{2})idy_2,$$
    which under the identification $\cA_{\rm nf}(T_2)
    \cong \RR^2$ is the line segment  from $(1,\frac{3}{2})$ to
$(2,\frac{7}{2})$. We define the path $B_t$
    in $\cA_{\rm nf}(T_2)$ by first following the path $\al\be E_\be$ from
    $\al\be c_0$ to $\al\be^2 c_0$, then $\al\be^2 E_\be$ from $\al\be^2
    c_0$ to $\al\be^3 c_0$, then $\al\be^3 E_\al$ from $\al\be^3 c_0$ to
    $\al\be^3\al c_0=\al\be^2\ga^2\al\be c_0$ (the last equality is
    by the relation $[\al,\be]=\ga^{-2}$).

So, $B_t$ is a path in $\cA_{\rm nf}(Y_2)$ from $\al\be c_0$ to
    $\al\be^2\ga^2(\al\be c_0)$ with the property that $Q(B_t)$ is
    homotopic rel endpoints to the straight line from $(1,\frac{3}{2}) $
    to $(2, \frac{7}{2}) $ in $\RR^2-\ZZ^2$. Hence, as before, define a
    path $A_{2,t}$ in $\cA_{\rm nf}(Y_2)$ which is homotopic rel endpoints to
    $B_t$ in $ Q^{-1}(\RR^2 - \ZZ^2)$ and has  
    $Q(A_{2,t})= (1+t, \frac{3}{2}+2t)$. Note that
    $SF(A_{2,t};Y_2;P^+)=k_\al+2k_\be$ while, as before,
    $SF(A_{1,t};Y_1;P^+)=k_\al$. Gluing together $A_{1,t}$ on $Y_1$ and $A_{2,t}$
    on $Y_2$, we obtain a path $A_t$ of connections on $Y_1 \cup Y_2$.
    Note that $A_1=gA_0$, where $g$ is the union of $\al$ on $Y_1$
    and $\al\be^2\ga^2$ on $Y_2$. Since $\deg(g)= \deg(\al) +
    \deg(\al\be^2\ga^2)=0$ by Theorem \ref{oiewr},  it follows that
    $SF(A_t;Y_1\cup Y_2)=0$. On the other hand, Lemma \ref{splittingthm}
    says $$SF(A_t;Y_1 \cup Y_2)= SF(A_{1,t};Y_1;P^+) +
    SF(A_{2,t};Y_2;P^+)=2k_\al+2k_\be,$$  which yields the equation
    $0=2k_\al+2k_\be$.

Solving these two equations shows that $k_\al = 2$ and $k_\be = -2$ and
    completes the proof of the theorem.\endproof

\section{Dehn surgery techniques for computing gauge theoretic
    invariants} \label{sec-csandrho}   

In this section, we   apply the results from Sections 3 and 4
to develop formulas for
a variety of gauge theoretic invariants
of flat connections on Dehn surgeries $X=Y \cup_T Z.$
Given a path of flat connections on $Z$
whose initial point is the trivial connection and
whose terminal point extends  flatly  over $X$, we extend
this to a path $A_t$ of connections on $X$  such that
$A_0 = \Th$ and $A_1$ is flat.
We then apply Theorem  \ref{zorro} and the results of Subsection
\ref{C2sfonW} to derive a general formula for
 the $\CC^2$--spectral flow along this path.
We also give a formula for the
    Chern--Simons invariant of $A_1$ as an element in $\RR$ rather than
$\RR/\ZZ$.
    These formulas allow computations of the spectral flow
    and the Chern--Simons invariants
    in terms of easily computed
     homotopy invariant quantities associated to the path
     $(m_t,n_t) \subset \RR^2$
     introduced in Subsection \ref{defnofAt}.
    To illustrate how to use the formulas in practice,
    we present  detailed calculations  for $\pm 1$ surgery on the
    trefoil in
    Subsection \ref{trefoil}.

   Combining the formula  for the spectral flow with the one
   for the Chern--Simons invariant
    leads to a computation of the $SU(2)$
    rho invariants of Atiyah, Patodi, and Singer in  Subsection
    \ref{rhosubsec}.  Our ultimate aim  is to develop methods for computing
    the correction term for the $SU(3)$ Casson invariant \cite{bh}.
    Summing the rho invariants yields the correction term
    provided the $SU(2)$ representation variety is {\it regular} as a subspace
    of the $SU(3)$ representation variety
      (Theorem
    \ref{corrected}).   In Section 6,
    we will extend these computations to all surgeries on $(2,q)$ torus knots.

\subsection{Extending paths of connections to X} \label{constrpt}
Throughout this section, we denote by $A, B, C,$ and $a$ connections
on $X,Y,Z$ and $T$, respectively. With respect to the manifold splitting
$X= Y \cup_T Z,$ we have $A=B \cup_a C$.

Our starting point is the following. We are given a path $C_t$ of
    $SU(2)$ connections on $Z$ in normal form on the collar which are flat for
    $t$ near $0$ and at $t=1$ (in all the examples considered in
    this paper,  $C_t$ is flat for {\it all}
    $t$) such that
\begin{enumerate}   \item  $C_0= \Th,$ the trivial connection on $Z.$
 \item $C_1$  extends flatly over $X=Y\cup_T Z$.
  \item For all $t>0$  the restriction of
    $C_t$ to the boundary torus has nontrivial holonomy.   \item
    For all small positive $t$, $C_t$ is a nontrivial reducible
    connection.\end{enumerate}

Let $a_t = C_t |_T$ be the restriction of $C_t$ to the  boundary.
Then since $C_t$ is in normal form, we have   $$a_t=-m_t  i  dx - n_t  i
    dy$$
    for $(m_t,n_t) \in \RR^2.$ Conditions 1--4  imply that
    $(m_0,n_0)=(0,0)$, $m_1\in \ZZ$, $(m_t,n_t)\in \RR^2-\ZZ^2$ for
    $t>0$, and $n_t=0$ for small positive $t$.
Moreover, by reparameterizing $C_t$ and gauge transforming if necessary,
   we can assume
    that conditions 1--3
    of Subsection \ref{defnofAt} hold.

The path $(m_t,n_t)$ will usually be described starting with a path of
    representations. The following proposition is helpful.  This 
proposition follows from a relative version of the main theorem of
\cite{Fine-klassen-kirk}. Its proof, which follows the same outline as
\cite{Fine-klassen-kirk}, is omitted.

\begin{proposition} {\sl Suppose $\rho\co [0,1]\to {\rm Hom}(\pi_{1}Z,SU(2))$ is a
continuous
    path of representations with $\rho_t(\mu)=e^{2\pi i m_t}$ and
    $\rho_t(\la)=e^{2\pi i n_t}$.  Then there exists a path of flat
    connections $C_t$ on $Z$ in normal form such that 
    $\hol_{C_t} = \rho_t$ and
     $C_t|_T= -m_tidx
    - n_tidy$.  Moreover, if the initial point $C_0$ is specified, then  
     $C_t$ is uniquely determined  up to a gauge transformation homotopic
    to the identity for each $0<t\leq 1$.}\end{proposition}

Our next task is to construct a path $B_t$   of connections on
    $Y$ agreeing with $C_t$ along the boundary $T$.
    The resulting path $A_t =  B_t \cup_{a_t} C_t$ of connections
    on $X=Y \cup_T Z$ should satisfy conditions 1--3 of Subsection
    \ref{defnofAt}.

    We begin by defining three integers $a,b,c$ in terms of
    the path $(m_t,n_t)$.
    First, set $$a= m_1 \quad \text{and}  \quad b=[n_1],$$
    where $[x]$  denotes the
    greatest integer less than or equal to $x$. Since $C_1$ extends flatly over
    $X,$    $\hol_{C_1}(\mu) = 1$ hence $a \in \ZZ.$)

\begin{figure}[ht!]\small
\begin{center} 
\includegraphics[scale=.4]{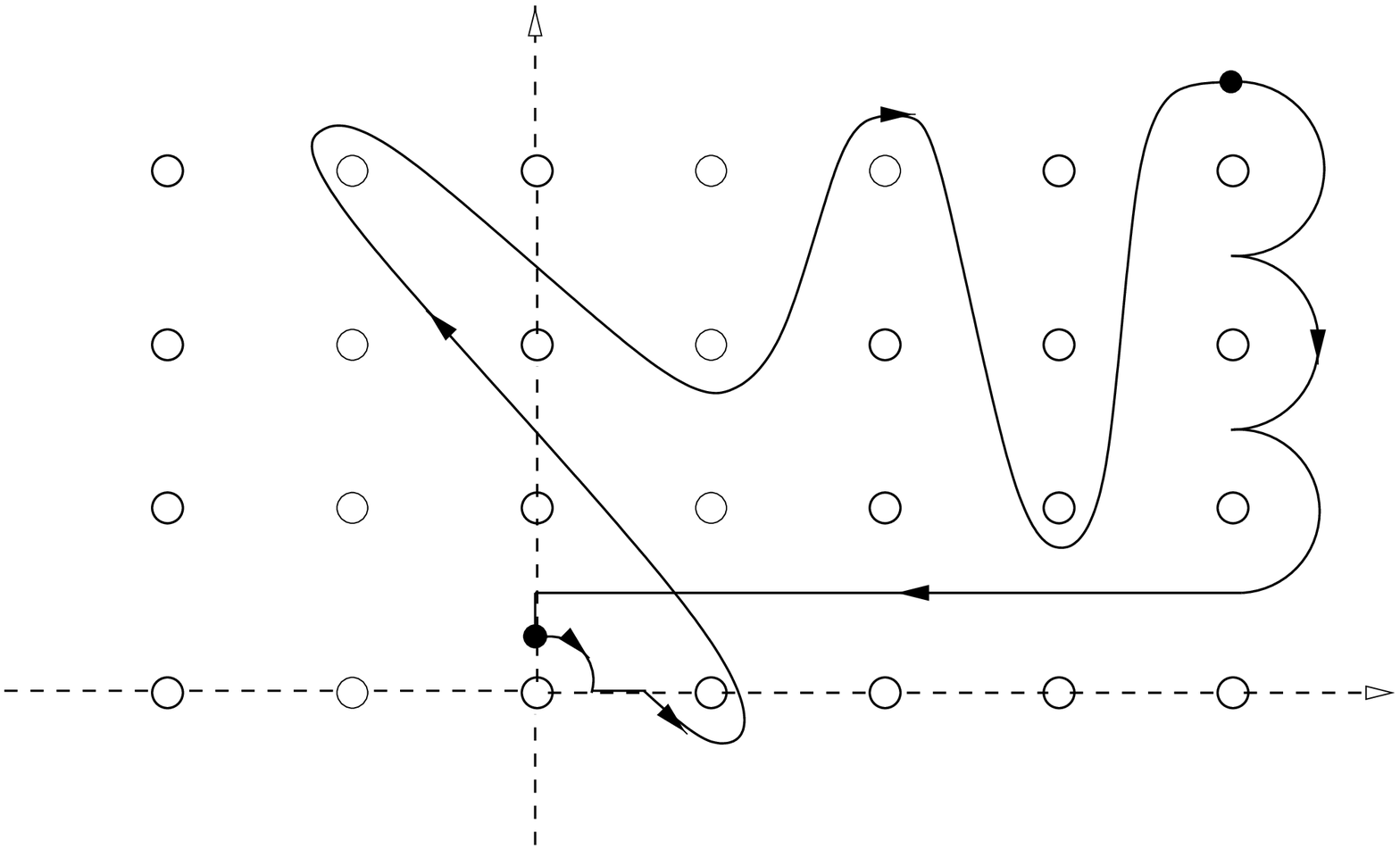}
\end{center} 
\vskip-5.1cm \hskip9.8cm $(m_1,n_1)$ 
\vskip-.2cm \hskip6.1cm $(m_t,n_t)$ 
 \vskip 1cm \hskip 10.3cm $p_3$ 
 \vskip 1.25cm \hskip 4.6cm $(0,\de)$ \hskip 2.1cm $p_4$ 
 \vskip 1cm
\caption{The loop $\ell$ and numbers $a=4$, $b=3$ and $c= -2(1-9)=16$}
\label{pathfig} 

\end{figure}

 Choose $\de >0$   as in condition 3 of Subsection \ref{defnofAt}.
Define the loop   $\ell =p_1 \cdot p_2 \cdot p_3 \cdot p_4\cdot p_5$ to
    be the composite of the following five paths in $\RR^2 - \ZZ^2$
    (see Figure \ref{pathfig}):
    \begin{enumerate}
\item[(i)] $p_1=\bar{\si}$ is the small quarter circle
    starting at $(0,\de )$ and ending at $(\de ,0)$, ie,
    $$p_1(t)=\left(\de  \cos(\tfrac{(1-t)\pi}{2}),\de 
\sin(\tfrac{(1-t)\pi}
    {2})\right).$$
\item[(ii)]
    $p_2=\eta$  is the path $ (m_t,n_t)$ for $\delta \leq t \leq 1$.
\item[(iii)]
    $p_3$ is the path from $(m_1,n_1)$ to $(m_1,n_1-b)$ which traverses
    the union of $|b|$ right hand semicircles of radius $\frac{1}{2}$.
    Setting $\varepsilon = \pm 1$ according to the sign of $b$,
    then $$p_3=\bigcup_{k=1}^{|b|} \{ (m_1+ \tfrac{1}{2}\cos t,
    n_1-\tfrac{\varepsilon}{2}(k+\sin t)) \mid -\tfrac{\pi}{2}
    \leq t \leq \tfrac{\pi}{2} \}.$$
\item[(iv)] $p_4$ is
    the horizontal line segment from $(m_1, n_1-b)$ to $(0,n_1-b)$.
\item[(v)] $p_5$ is the short vertical line segment from
    $(0,n_1-b)$ to $(0,\de )$.\end{enumerate}

We now define the integer $c$ in terms of the linking number of $\ell$
with the integer lattice $\ZZ^2 \subset \RR^2.$

\begin{definition}\label{defnofc}
Given any oriented closed loop  $L$ in $\RR^2-\ZZ^2$,
    define the linking number $ \lk(L,\ZZ^2)$ of $L$ and $\ZZ^2$
    to be the algebraic number of lattice
    points enclosed by $L$, normalized so that
    if $L(t) = (\de  \cos t,\de  \sin t)$
    for $t \in [0,2\pi],$ then $\lk(L,\ZZ^2)=1$.

Using the loop $\ell$ constructed above, we define an integer by
setting
    $$c=-2\ \lk (\ell,\ZZ^2).$$\end{definition}

Figure \ref{pathfig} shows how to compute the
integers $a,b$ and $c$ from the graph of $(m_t,n_t)$.

We can now define a path $B_t$ of connections on $Y$.
    \begin{definition}\label{ythgjg} Set $B_t = -q(r) t i dx$ for
    $0 \leq t \leq \de $, where $q(r)$ is defined in Equation
(\ref{eqnofbmn}).

  Also set   $$B_1=\al^{a} \be^b\ga^c (-(n_1-[n_1]) i dy)= \ga^c \al^{a}(-
    n_1  i dy).$$ Notice that $B_1|_T= a_1=-m_1 i  dx - n_1 i
    dy.$

Finally, define $B_t$ for $\de  < t \leq 1$ to be any path of connections
    in normal form interpolating from $B_\de $ to $B_1$ and
    satisfying $B_t|_{T}=a_t$. Such a path exists since the space of
    connections on $Y$ with a given normal form on the boundary is
    contractible.\end{definition}

Since the restrictions of $B_t$ and $C_t$ to the torus agree, and since
    they are in normal form, they can be glued together to form a
    path $$A_t=B_t\cup_{a_t}C_t$$ of connections on $X$. This path
    satisfies all the requirements of  Subsection \ref{defnofAt}, and
    hence we can apply Theorem \ref{zorro}.
Notice that the flat connection $B_1$ depends only on the
    homotopy class (rel endpoints) of the path $(m_t,n_t), \ t\in [\de ,1]$ in
    $\RR^2-\ZZ^2$.

\subsection{Computation of the spectral flow}\label{mandns}

\begin{theorem}\label{nicerformula} {\sl Let $A_t$   be the path of
    connections on $X$ constructed above and let $a$ and $b$ be the
    integers defined above for the path $(m_t,n_t)$ (so  $a=m_1$
    and $b=[n_1]$). Then} \begin{equation}\label{realniceformula}
    SF(A_t;X)=SF(A_\eta;Z;P^-) + 2(a-b)-2.\end{equation}\end{theorem}
\noindent {\bf Proof}\qua By Theorem \ref{zorro}, we only need to  show that
    $SF(A_{\bar{\si}\cdot \eta};Y;P^+)=2(a-b)$.   The path of
    connections $A_{\bar{\si}\cdot\eta} |_Y$ starts  at the flat
    connection $-\de   i   dy$ and ends at $B_1$. Its projects under  $Q$
    to
    the path $\bar{\si}\cdot \eta$ in $\RR^2-\ZZ^2$.

Recall that the path $\bar{\si}$ is the   small quarter circle  from
$(0,\de )$
    to $(\de , 0)$  and that $\eta$ is just $(m_t,n_t)$ starting at
    $t=\de $. Referring to the notation and results of Section
    \ref{sfsection}, we see that the homotopy class rel boundary of
    the path $\bar{\si}\cdot \eta$ in $\RR^2-\ZZ^2$  uniquely determines
    a word $w$   in $\al$ and $\be$ in the group $G$.
    This word  uniquely specifies a path $P$ in the
    Cayley graph, which we regard  as a path of connections
    on the solid torus.

For example, the word $w=\al^3\be\al^{-1}$ determines the path
     $$P=E_\al \cdot \al E_\al \cdot
    \al^2 E_\al \cdot \al^3 E_\be \cdot \al^3\be\al^{-1} E_\al^{-1},$$ where
    $E_\al^{-1}$ means $E_\al$ traversed backwards.
    By construction, the endpoint of this path  is
    $ \al^3\be \al^{-1} \cdot(-\tfrac{1}{2} i dy)$.

    Given any word $w$ in $\al$ and $\be,$ the associated path
    $P$ goes from $-\frac{1}{2} i dy$ to $w \cdot (
-\frac{1}{2} i dy)$.
    Thinking of $w$ as an element of $G$ and using Lemma
    \ref{commutatior} to put $w$ into normal form,
    it follows that $$w= \al^a \be^b \ga^c$$
    where $a=m_1, b=[n_1]$ and $c$ is defined relative to
   the path $\bar{\si} \cdot \eta$ as in Definition \ref{defnofc}.
 Thus, the terminal point of $P_w$
     is the flat connection $ \al^a \be^b \ga^c \cdot (-\frac{1}{2} idy).$

We now construct a path $\widetilde{P}$ by pre- and post-composing
   the given path $P$ so the initial and terminal
   points agree with those of the path   $A_{\bar{\si}\cdot\eta }|_Y$.
    This is done by adding short segments of
   nontrivial, flat connections. This will
   not affect the spectral flow since
    any nontrivial representation $\rho\co  \pi_1 Y \to SU(2)$
  has $H^{0+1}(Y;\CC^2_\rho)=0$.

Consider first the line segment from $(0, \de )$ to $(0,\frac{1}{2})$ in
$\RR^2$. Since it misses the integer lattice, it lifts to a straight line   
from $- \de  i dy$ to $- \frac{1}{2} i dy$. This lift is a path of 
nontrivial flat connections on $Y$.
Now consider the line segment
from $(m_1,[n_1] + \frac{1}{2})$ to $(m_1,n_1).$
It also misses the integer lattice,
hence it lifts to
a straight line from $ \al^a \be^b \ga^c \cdot (-\frac{1}{2} idy)$,
 the terminal point of $P$,  to $  \al^a \be^b \ga^c \cdot(-(n_1 - [n_1]) idy),$
the flat connection $B_1$. The second lift is also a path of nontrivial
flat connections on $Y.$

Precomposing $P$ by the first lift and post-composing
by the second defines a path $\widetilde{P}$ with the same
 $\CC^2$-spectral
flow as  $P$. Notice that the initial and terminal
points of $\widetilde{P}$ agree with those of $ A_{\bar{\si}\cdot
\eta}|_Y.$
  By Theorem \ref{constants}, if $g\in G$, then the
    spectral flow on $Y$ with $P^+$ boundary conditions along $g\cdot
    E_\al$ equals $2$ and along $g\cdot E_\be$ equals $-2$.
    Thus the spectral flow along the path $P$ is equal to
    $2(a-b)$. But since the spectral flow along
     $\widetilde{P}$ is the same as that along $P$,
     and since $ \widetilde{P}$ is homotopic rel endpoints
     to $  A_{\bar{\si}\cdot \eta}|_Y$,
     this shows that $SF(A_{\bar{\si}\cdot \eta};Y;P^+)=2(a-b)$
   and completes the
    proof.\endproof

\subsection{The Chern--Simons invariants}
The Chern--Simons function is  defined on the space $\cA_X$ of connection
    1--forms on a closed manifold $X$ by $$cs(A) = \frac{1}{8 \pi^2}
    \int_X {\rm tr}(A \wedge dA + \tfrac{2}{3} A \wedge A \wedge
    A).$$ With this choice of normalization, $cs\co \cA_X \to \RR$
    satisfies $cs(g\cdot A)= cs(A)- \deg g $ for gauge transformations
    $g$ (recall that $g\cdot A=gAg^{-1}-dg~ g^{-1}$).
    Since computing $cs$ modulo $\ZZ$ is not sufficient
    for the applications we have in mind, we work with connections
    rather than gauge orbits.

    Using the same path $A_t = B_t \cup_{a_t} C_t$ of
    connections on $X = Y \cup_T Z$
    as in Subsection \ref{constrpt}, we show how to compute
    $cs(A_1) \in \RR$.
This time, the initial data is a path of  {\it flat}
 connections $C_t$ on $Z$ in normal form  on the collar
    with $C_0$ trivial and $C_1$ extending flatly over $X$.

The restriction of $C_t$ to the boundary determines path
$(m_t,n_t)$
    (ie, $C_t|_T=-m_tidx -n_tidy$) which was used in Subsection
    \ref{constrpt} to construct a path $A_t$ of connections on $X$
    starting at the trivial connection and ending at a flat connection
    $A_1$. Reparameterize the path $(m_t,n_t)$ so that the coordinates
    are differentiable.  (It can always be made piecewise analytic
    by the  results of \cite{Fine-klassen-kirk},  and hence using
    cutoff functions  we can arrange that it is smooth.)

\begin{theorem}\label{csthm}  {\sl  The Chern--Simons invariant of $
    A_1$ is given by the formula $$cs( A_1)= -c+2\int_0^1n{\frac{dm}
    {dt}} \ dt. $$}

\end{theorem}

\noindent {\bf Proof}\qua We follow the  proof of Theorem 4.2 in \cite{chern-simons-1},
    being careful not to lose integer information.  Let $T(A)$ denote the
transgressed second Chern form $$T(A)= {\frac{1}
    {8\pi^2}}{\rm tr}(dA\wedge A+{\frac{2}{3}}A\wedge A\wedge A).$$
    Then $$cs(A_1)=\int_Y T(B_1) + \int_Z T(C_1)$$ since $A_1=B_1\cup
    C_1$ on $X=Y\cup_T Z$.  We compute these terms separately, using the
    following lemma.
\begin{lemma}\label{cslemma} {\sl Let $W$ be an oriented
    3--manifold with oriented boundary $T=S^1\times S^1$. Let $A_t$
    be a path of flat connections in  normal form on $W$. Assume
    that $A_t|_T= - m_t i dx - n_t i dy $, where $\{dx, dy\}$ is
    an oriented basis of $H_{1}(T;\ZZ)$. Then }$$\int_W
T(A_1)-\int_W
    T(A_0)=\int_0^1 (m\tfrac{dn}{dt}-\tfrac{dm}{dt}n) dt.$$\end{lemma}
    \noindent {\bf Proof}\qua Orienting $I\times W$   and $I\times \partial
    W$ as products and using the outward normal first convention,
    one sees that  the boundary
    $$\partial(I\times W)= (\{1\}\times W)-(\{0\}\times W) -
I\times\partial W.$$

The path of connections $A_t$ on $W$ can be viewed as a connection ${\bf
    A}$ on $I\times W$. Then the curvature form $F^{\bf A}$ equals
    $dt\wedge \omega$ for some $1$--form $\omega$. Hence $c_2({\bf
    A})=\tfrac{1}{4\pi^2}$tr$(F^{\bf A}\wedge F^{\bf A})=0$.

Using Stokes' theorem as in \cite{chern-simons-1}, one computes that
    \begin{equation}\label{goofy} 0=\int_{I\times W}c_2({\bf
    A}) = \int_{W} T(A_1) -T(A_0) -\int_{I\times \partial W}
    T({\bf a}),\end{equation}  where ${\bf a}$ denotes the connection
    $- m_tidx-n_tidy$ on $I\times \partial W$.
    Since $d {\bf a} = -\tfrac{dm}{dt}i dt dx-\tfrac{dn}{dt}i dt dy,$
    it follows that
$$d{\bf a}\wedge {\bf a} = \left( -\tfrac{dm}{dt}n+m\tfrac{dn}{dt}\right)
dxdy.$$
 Clearly ${\bf a}\wedge {\bf a}\wedge {\bf a}=0$, so  
$$T({\bf a})=\tfrac{1}{4\pi^2}(-\tfrac{dm}{dt}n+m\tfrac{dn}{dt})dtdxdy.$$
Hence\begin{eqnarray*}\int_{I\times \partial W} T({\bf a})&=&
    {\tfrac{1}{4\pi^2}}\int_{I\times \partial W}
    (-\tfrac{dm}{dt}n+m\tfrac{dn}{dt})dtdxdy\\
    &=&\int_0^1(-\tfrac{dm}{dt}n+m\tfrac{dn}{dt})dt.\end{eqnarray*}
    Substituting this into Equation (\ref{goofy}) finishes the proof
    of Lemma \ref{cslemma}.\endproof

    Returning to the proof of Theorem
    \ref{csthm},  we use Lemma \ref{cslemma} to compute   $\int_Z T(C_1)$.
Since
    $C_t$ is a path of flat connections on $Z$ starting at the trivial
    connection, Lemma \ref{cslemma} implies that $$ \int_Z T(C_1)=-\int_0^1
    (m\tfrac{dn}{dt}-\tfrac{dm}{dt}n) dt.$$ The sign change occurs because
    $\partial Z=-T$ as oriented manifolds.

Next we compute the  term $\int_Y T(B_1)$.  
Recall from Definition \ref{ythgjg}
    that $B_1=\ga^c \al^{a}(- n_1  i dy). $  Since $\ga$ is a degree 1
    gauge transformation supported in the interior of $Y$ and $cs(g\cdot
    A)=cs(A)-$deg$(g)$, $$\int_Y T(B_1)= -c +\int_Y T(\al^a(-n_1 i
    dy)).$$

Consider the path of flat connections on $Y$, $\tB_t=\al^a(-tn_1 idy), \
    t\in [0,1]$. Then $\tB_0=\al^a(\Theta)=-d(\al^a)(\al^a)^{-1}$
    and $\tB_1=\al^a(-n_1 i dy)$. The restriction of $\tB_t$ to the
    torus is $\tB_t|_{T}=\tal^a(-tn_1idy)=- a i dx - t n_1 i dy$.
    Recall that $a=m_1$. Applying Lemma \ref{cslemma} we conclude
    that $$\int_Y T(\al^a(-n_1 i dy))= m_1n_1+\int_Y
    T(-d(\al^a)(\al^a)^{-1}).$$  But $\int_Y T(-d(\al^a)(\al^a)^{-1})=0$
    since $-d(\al^a)(\al^a)^{-1}$ has no $dy$ component,
thus $$cs(A_1)= -c+m_1n_1 -\int_0^1 (m{\tfrac{dn}{dt}}-{\tfrac{dm}{dt}} n)dt
    =-c+2\int_0^1\tfrac{dm}{dt}n dt.\eqno{\qed}$$

 \subsection{Example: $\pm$ 1 Dehn surgery on the trefoil}\label{trefoil}

 In this section,
we use our
previous results to determine the $\CC^2$--spectral flow and
the Chern--Simons invariants
for flat connections on the homology spheres obtained
 by $\pm 1$ surgery on the right--hand trefoil  $K$.
More general results for surgeries on torus knots will be given
in Section \ref{torusknot}.

The key to all these  computations is a concrete description
of the $SU(2)$ representation variety
of the knot complement (see \cite{klassen-thesis}).
Let $K$ be the right hand  trefoil knot in $S^3$ and
 let $Z$ be the 3--manifold with boundary obtained by
    removing an open tubular neighborhood of $K$.
    Its fundamental group has presentation
$$\pi_1Z=\langle  x,y \mid  x^2=y^3\rangle.$$
    There are
    simple closed curves $\tmu$ and $\tla$ on $\partial Z = T$
     intersecting transversely in one point called the
   meridian and longitude of the knot complement.
   We use the right hand rule to orient the pair
   $\tla,\tmu$ (see Subsection \ref{cohtorusknot} for more details).
    In $\pi_1 Z$,
    $\tmu$ represents $xy^{-1}$ and $\tla$ represents
    $x^2(xy^{-1})^{-6}$
(cf.~Equation (\ref{natmerlong})).

\begin{figure}[ht!]\small
\begin{center}
\includegraphics[scale=0.8]{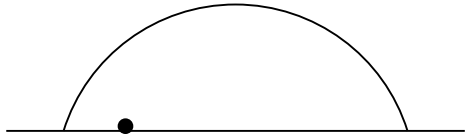}
\end{center}
\caption{$SU(2)$ representations of the trefoil}
\label{trefoilrep}
\end{figure}

The representation variety $\fR_{SU(2)}(Z)$
   can be described as the identification space
    of two closed intervals where the endpoints of the first interval are
    identified with two points in the interior of the second
    (see Figure \ref{trefoilrep}).
(In general, the representation variety of any
torus knot complement is a singular
    1--manifold with  `T' type intersections called
    {\it $SU(2)$ bifurcation points}, see \cite{klassen-thesis}.)

 Since
    $\tmu$ normally generates $\pi_1 Z,$
    any abelian representation $\be\co \pi_1 Z \to SU(2)$ is uniquely determined
    by the image $\be(\tmu)$.
    To each $t \in [0,\frac{1}{2}] $
    we associate the abelian representation
    $\be_t\co \pi_1 Z \to SU(2)$
    with
    $\be_t(\tmu)= e^{2\pi it} $.
    Thus, the interval $[0,\frac{1}{2}]$
     parameterizes the conjugacy classes of abelian or
    {\it reducible} representations.

The arc of
    nonabelian or {\it
    irreducible}  conjugacy classes of representations
    can be parameterized by the open interval
    $(0,1)$ as follows. For $t \in [0,1]$, let $\rho_t \co  \pi_1 Z \to SU(2)$
    be the representation  with $\rho_t(x) = i$
    and $$\rho_t(y) = \cos(\tfrac{  \pi}{3})+ \sin(\tfrac{ \pi}{3})
    (\cos( t \pi ) i + \sin(t \pi) j)   . $$ In \cite{klassen-thesis} it is
proved
    that every irreducible $SU(2)$ representation of $\pi_1 Z$
    is conjugate to one and only one $\rho_t$ for some $t \in (0,1).$
    The endpoints of $\rho_t$ coincide with the reducible
    representations  at ${1/12}$ and ${5/12}$.

    Restriction to the boundary defines a map
    $\fR_{SU(2)}(Z)\to \fR_{SU(2)}(T)$.
   To apply Theorems \ref{nicerformula} and \ref{csthm}
   to manifolds obtained by surgery on $K$,
   we need to lift the image
    $\fR^*_{SU(2)}(Z)\to \fR_{SU(2)}(T)$ under
    the branched cover $f\co \RR^2 \to \fR_{SU(2)}(T)$
    of Equation (\ref{pillowpara}). It is important
    to notice that $f$ depends on the surgery coefficients.
    Specifically, $f$ is defined in Equation (\ref{pillowpara})
    relative to the
    the meridian and longitude
    of the solid torus, as opposed to the meridian and
    longitude of the knot complement.
    We denote the former by $\mu$ and $\la$ and the latter by $\tmu$ and
$\tla$.
    For the manifold $X_k$ obtained by $\frac{1}{k}$ surgery on $K$, we have
    $\mu = \tmu \tla^k$ and $\la = \tla$.
    
 \begin{figure}[ht!]\small
\begin{center}
\includegraphics[scale=.30]{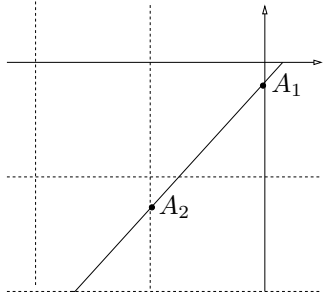}
\end{center}
\vskip-3.2cm \hskip8cm $A_1$
\vskip1.2cm   \hskip6.5cm $A_2$   
 \vskip 1.5cm
 \caption{ Two flat connections extending over $+1$ surgery on the right hand trefoil}
\label{surg1}
\end{figure}

For the Poincar\'e homology sphere (denoted here by $X_{+1}$),
 Proposition \ref{kDsurg}  implies
that one such lift is given by the
curve   (see Figure \ref{surg1})
 $$
    R_1(t)= (1-t)(\tfrac{1}{12},0) + t(-\tfrac{19}{12}, -2), \ \ 0 \leq t
\leq 1.
$$
All other lifts
are obtained by translating $R_1$ by integer pairs and/or
    reflecting it through the origin.
    
    \begin{table}[ht]\small  \renewcommand\arraystretch{1.5} \noindent\[
    \begin{array}{|c|  c c c|c|c|c|} \hline& \ a & b & c \ & 2\int
    m'n &cs(A) & SF(\Th,A)   \\ \hline A_1 & 0 & -1 & 0  &
    \frac{1}{120}&\frac{1}{120} & 0 \\ \hline A_2 & -1 & -2 & 2  &
    \frac{169}{120}&-\frac{71}{120} & 0  \\ \hline\end{array} \]
   \caption{$X_{+1} =$ +1 surgery on the right hand
    trefoil}\label{table1}\end{table}

    A  representation $\rho\co \pi_1 Z \to SU(2)$ extends over $X_{+1}$
    if and only if $\rho(\mu) =1$, hence it follows that the irreducible
    representations of $X_{+1}$ correspond to the points of $R_1$
    where the first coordinate is an integer.
    Figure   \ref{surg1} shows two such
points which represent
  two flat connections $A_1$ and $A_2$.
    Let $(m_t, n_t)$ be the path described as the composition of the
   horizontal line segment from $(0,0)$
    to $({\frac{1}{12}},0)$   with the path $R_1(s)$.
Then $A_1$ and $A_2$ are the flat connections
constructed as in Subsection \ref{constrpt} using the path
$(m_t, n_t)$, stopping on the $R_1$ portion at $R_1(1/20)$ for $A_1$
and at $R_1(13/20)$ for $A_2$.

Using the  path   $A_{\xi \cdot \eta}$ constructed from
$(m_t,n_t)$ as in Subsections
    \ref{defnofAt} and \ref{zorrostate}, we compute  the numbers
    $a_i, b_i, c_i$ associated to $A_i$ for $i=1,2.$ We get that
     $a_1=0,b_1= -1$, and
    $c_1=0$. Similarly  $a_2=-1, b_2=-2$, and
    $c_2=2$.

 \begin{figure}[ht!]\small
\begin{center}
\includegraphics[scale=.30]{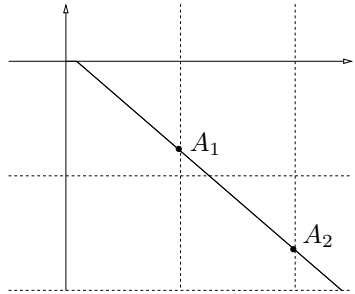}
\end{center}
\vskip-2.4cm \hskip6.7cm $A_1$
\vskip.8cm   \hskip8.2cm $A_2$   
 \vskip .5cm
 \caption{ Two flat connections extending over  $-1$ surgery on the right hand trefoil}
\label{surg2}
\end{figure}

It follows from Theorem \ref{2qsfpts} below that  $SF(A_\eta(t);
    Z;P^-)=0$ for both $A_1$ and $A_2$. By Equation
    (\ref{realniceformula}),
    $SF(\Theta, A_i; X_{+1})= 2(a_i-b_i) -2$, we conclude that
    $SF(\th,A_i;X_{+1})=0$ for $i=1,2.$
One can also compute the integral term  $ 2\int m'n$ arising
in Theorem \ref{csthm}, getting $2 \int m'n ={\frac{1}{120}}$ for
$A_1$ and $2 \int m'n ={\frac{169}{120}}$ for
$A_2$.
These results are summarized   Table \ref{table1}.

Similar computations for the manifold $X_{-1}$
are given in Table \ref{table2}. Here,
  Proposition \ref{kDsurg}  implies
that the lift of the image of $\fR^*_{SU(2)}(Z) \to \fR_{SU(2)}(T)$
under $f \co  \RR^2 \to \fR_{SU(2)}(T)$ is given by the
curve
$$
    R_1(t)= (1-t)(\tfrac{1}{12},0) + t(\tfrac{29}{12}, -2), \ \ 0 \leq t
\leq 1.
$$
As before, the numbers $a,b$ and $c$ and the integral term
$2\int m'n$ can be easily computed from   Figure \ref{surg2}.

 Using these results, we will determine the rho invariants of flat
    connections on $X_{+1}$ and $X_{-1}$ in Subsection
    \ref{rhosubsec}.

\begin{table}[ht]\small \renewcommand\arraystretch{1.5} \noindent\[
    \begin{array}{|c|  c c c|c|c|c|} \hline& \ a & b & c \ & 2\int
    m'n &cs(A) & SF(\Th,A) \\ \hline A_1 & 1 & -1 & -2  &
    -\frac{121}{168}& \frac{215}{168} & 2 \\ \hline A_2 & 2 & -2 &
    -6  & -\frac{529}{168}&\frac{479}{168} & 6 \\ \hline\end{array} \]
    \caption{$X_{-1}= -1$ surgery on the right hand
    trefoil}\label{table2} \end{table}

\subsection{The rho invariants} \label{rhosubsec} In this
    section, we present explicit formulas for the rho invariants
    based on our previous results. We first make clear {\it which}
    rho invariants we are computing.  Following \cite{APS}, an $SU(2)$
    connection $A$ on $X$ determines the self-adjoint odd signature
    operator  {\it with $\CC^2$ coefficients}:
    $$D_A\co \Om^{0+2}_X\otimes\CC^2\to \Om^{0+2}_X\otimes\CC^2.$$ The
    eta invariant of $D_A$, denoted here by $\eta_{A}(0)$, is the
    spectral invariant regularizing the signature;  it is the analytic
    continuation to $s=0$ of $$\eta_{A}(s)=\sum_{\la \neq 0}
    \frac{\hbox{sign}(\la)}{|\la|^{s}};$$ where the sum is over nonzero
    eigenvalues $\la$ of $D_A$.

If $A$ is a {\it flat} connection, then Atiyah, Patodi, and Singer show
    that the difference $$\varrho_{X}(A)=\eta_{A}(0)-\eta_{\Theta}(0)$$
    is a real number which is independent of the metric. Moreover,
    $\varrho_{X}(A)$ is  gauge invariant and hence defines a function $
    \varrho_X\co  \fM_{SU(2)}(X) \to \RR$  on the flat $SU(2)$ moduli space of
    $X$. Using the holonomy map to identify flat connections $A$
    and representations $\al \co \pi_1 X \to SU(2)$, the rho invariant
    can therefore also be viewed as a real-valued function on
    $\fR_{SU(2)}(X)$.

The rho invariants considered in this paper are those associated to the
    canonical representation of $SU(2)$ on $\CC^2$, {\it not} the
    adjoint representation on $su(2)$ which is more commonly studied in
    Donaldson and Floer theory. In the latter situation,  Fintushel
    and Stern developed a technique for computing the (adjoint) rho
    invariants of $SU(2)$ representations of Seifert--fibered spaces
    by extending them over the mapping cylinder of the Seifert
    fibration, viewed as a 4--dimensional orbifold \cite{FS}. This
    method   does not apply to our situation because generic fibers
    do not act trivially in the canonical representation  as they
    do in the adjoint representation.

\begin{theorem}\label{rhoinvts}  {\sl Suppose $C_t$ is a path of
    flat connections in normal form on $Z$ starting at the trivial
    connection and ending at a connection with trivial holonomy
    around $\mu$. Let $A_1$ be  any  flat connection $X=Y\cup_T Z$
    which extends $C_1$. Then the rho invariant of $A_1$ is given
    by the formula
    \begin{eqnarray}
    \varrho_{X}(A_1)&=&2SF(A_{\eta}(t);Z; P^-)+4(a-b+c) -2 \label{rhoformula} \\
    &&-\dim(\ker D_{A_1})-8\int m'n, \nonumber
    \end{eqnarray} where $C_t|T= -m_t i dx - n_t i dy$, and
    $a,b,c$ are the integer homotopy invariants of the path $(m_t,n_t)$
    defined in Subsection \ref{constrpt}.} \end{theorem}

\noindent {\bf Proof}\qua The rho invariant is  gauge invariant so every flat
    connection on $X$ gauge equivalent to $A_1$ has the same rho
    invariant.
Thus we are free to
use the path $A_t$ of connections constructed in Subsection \ref{constrpt}
from the path $C_t$ to compute $\varrho_{X}(A_1).$

A standard application of the Atiyah--Patodi--Singer index theorem shows
    that\begin{equation} \label{rhoAPS} SF( A_t;X) = 2 cs(A_1) +
    \tfrac{1}{2} \left(\varrho_X(A_1) - \dim(\ker D_\Theta) + \dim(\ker
    D_{A_1}) \right).\end{equation} This follows just as in the
    appendix to \cite{KKR}, keeping in mind that we are using the
    $(-\ep,-\ep)$-convention to compute spectral flow here whereas in
    that paper, the $(-\ep,\ep)$-convention is used (hence the sign
    change for the term $\dim(\ker D_{A_1})$).

Recall further that $\ker D_\Theta=H^{0+1}(X;\CC^2)\cong \CC^2$.  Using Theorem
    \ref{csthm} and Corollary \ref{nicerformula}  to substitute
    into Equation (\ref{rhoAPS}) and solving for $\varrho_X(A_1)$
    yields Equation (\ref{rhoformula}).\endproof

In general, from Equation (\ref{rhoAPS}), if $ \dim (\ker
    D_A)=0$, then
    \begin{equation} \label{rhorho} \varrho_{X}(A) = 2 SF(\Th,A;X)
    - 4 cs(A)+2.\end{equation}
    By Theorem \ref{vanish}, this holds for every nontrivial
    flat connection $A$ over a homology sphere $X$ obtained by surgery
    on a $(2,q)$ torus knot. \noindent

\begin{example} \label{rhorho1} Suppose $K$ is the right hand trefoil
    and consider the two sets of connections on $X_{\pm 1}$, the
    homology spheres obtained by $\pm 1$ surgery on $K$. Then, referring
    to Tables 1 and 2 and utilizing Equation (\ref{rhorho}), we conclude
    that:

{\bf Case 1}\qua For $+1$ surgery on $K$,   
    $\varrho_{X_{+1}}(A_1) = 59/30$  and $\varrho_{X_{+1}}(A_2) =
    131/30.$   

{\bf Case 2}\qua For $-1$ surgery on $K$,   
    $\varrho_{X_{-1}}(A_1) = 37/42$  and $\varrho_{X_{-1}}(A_2) =
    109/42.$
    \end{example}

Notice that while the quantities in Tables 1 and 2 depend on the choice
    of gauge representatives $A_1$ and $A_2$, the rho invariants
    do not. We shall
    extend these computations  considerably, first to all homology
    spheres obtained by Dehn surgery on the trefoil (Theorems
    \ref{23formulas} and \ref{-23formulas}) and later to all homology spheres obtained by surgery
    on a $(2,q)$ torus knot (Theorems \ref{+2qfull} and \ref{-2qfull}).

\subsection{The SU(3) Casson invariant} In \cite{bh}, an invariant of
    homology 3--spheres $X$ was defined by counting, with sign, the
    number of irreducible $SU(3)$ representations of $\pi_1(X)$ and
    subtracting a correction term. The correction term is given by
    a sum of $\CC^2$--spectral flows and Chern--Simons invariants applied
    to flat $SU(2)$  connections. One must typically incorporate
    the effect of perturbations on both of these sums, but in certain
    fortuitous cases the flat moduli space $\fM_{SU(3)}(X)$ is regular
    and no perturbations are needed. The aim of this subsection is to
    give a simple formula for the correction term  in this special
    case.

To begin, we recall the definition of the $SU(3)$ Casson
    invariant (cf.~Section 5 of \cite{bh}).

\begin{definition} \label{su3}  The $SU(3)$ Casson invariant for a
    homology sphere $X$ is given by the sum
    $$
    \la_{SU(3)}(X) = \la'_{SU(3)}(X) + \la''_{SU(3)}(X)
  $$ where
  \begin{eqnarray*}
    \la'_{SU(3)}(X) &=&   \sum_{\quad [A] \in {\fM}^*_{SU(3),h}(X)}
    (-1)^{SF_{su(3)}(\Th,A;X)}    \\ \la''_{SU(3)}(X) &=&   \sum_{ [A] \in
    {\fM}^*_{SU(2),h}(X)} (-1)^{SF_{su(2)}(\Th, A;X)} ( SF (\th,A;X) -2 cs(\hA)
    +1).\end{eqnarray*}  These are the first and second sums,
    respectively, of Definition 5.2 in \cite{bh}. All the spectral
    flows are taken with respect to the twisted odd signature operator
    $D_A$
    (this was denoted $K_A$ in \cite{bh}).
     The notation $SF(\Theta,A;X)$ refers to the $\CC^2$-spectral
    flow, ie, taking $SU(2)$ acting on $\CC^2$ (and counting complex
    eigenvectors) just as above. (The analogous term $SF_{{\mathfrak
    h}^\perp}(\th,A)$  in  Definition 5.2 of \cite{bh} counts {\it
    real} eigenvectors, which adds a factor
    of $\frac{1}{2}$ in front
    of the second sum in Definition 5.2 of \cite{bh}.).
     The notation $SF_{su(3)}$ and
    $SF_{su(2)}$ refers to the adjoint representations, ie, $SU(3)$
    acting on $su(3)$ and $SU(2)$ acting on $su(2)$ by the  adjoint
    representation  (and count  real eigenvectors). The function
    $h$ is a perturbation  function used to perturb the flatness
    equations. Then $\fM^*_{SU(3),h}(X)$ denotes the  moduli space of
    {\it irreducible} $h$-perturbed-flat $SU(3)$ connections on $X$, and
    similarly $\fM^*_{SU(2),h}(X)$ denotes the moduli space of {\it
    irreducible} $h$-perturbed-flat $SU(2)$ connections on $X$.

Notice that $\la_{SU(3)}$ is independent of the underlying orientation
    on the homology sphere. In fact this is   true for
$\la'_{SU(3)}$
    and $\la''_{SU(3)}$, namely $\la'_{SU(3)}(-X) = \la'_{SU(3)}(X)$ and
    $\la''_{SU(3)}(-X) = \la''_{SU(3)}(X)$.

Neither $ \la'_{SU(3)}(X)$ nor $ \la''_{SU(3)}(X)$ is generally
    independent of the choice of perturbation $h$, which must be
    small and chosen so that $\fM_{SU(3),h}(X)$ is regular as in Definition
    3.8 of \cite{bh}. To evaluate the correction term
    $\la''_{SU(3)}(X)$, one must also choose a representative $A$
    for each $[A] \in {\fM}^*_{SU(2),h}(X)$  along with a nearby flat,
    reducible connection $\hA.$\end{definition}

In certain cases, including surgeries on $(2,q)$ torus knots, the
    $SU(3)$ moduli space is regular. In this case the invariant
    $\la_{SU(3)}$ is calculable without perturbing.  In fact, whenever
    the $SU(2)$ moduli space is regular according to Definition 3.8
    of \cite{bh}, one can compute the correction term $\la''_{SU(3)}(X)$
    in terms of $SU(2)$ rho invariants.

\begin{theorem} \label{corrected}   {\sl Suppose $X$ is a homology
    sphere with $H^1(X; su(2)_A)=0$ and $H^1(X;
    \CC^2_A) =0$
    for every irreducible flat $SU(2)$
    connection $A$ on $X$.
  The first condition ensures that the moduli space
    $\fM_{SU(2)}^{*}(X)$ is a compact,
    0--dimensional manifold, and the second implies that the points
    in $\fM_{SU(2)}^{*}(X)$ are not limits of arcs of irreducible flat $SU(3)$
    connections. Then the correction term can be written as a
    sum of rho invariants, specifically } \begin{equation} \label{su3rho}
     \la''_{SU(3)}(X) = \sum_{[A] \in \fM^*_{SU(2)}(X)}
    (-1)^{SF_{su(2)}(\Th,A;X)}  \varrho_X(A)/2.\end{equation} 
    \end{theorem}
    \noindent {\bf Proof}\qua This follows by
    taking $\hA = A$ in Definition \ref{su3} (which is allowed
    since $\fM^*_{SU(2)}(X)$ is regular as a subspace of $\fM_{SU(3)}(X)$
    by hypothesis) and making a direct comparison with Equation
    (\ref{rhorho}).\endproof

In the next two examples, we present computations of the $SU(3)$ Casson
    invariant for $\pm 1$ surgery on the trefoil. In addition to
    the fact that the $SU(3)$ moduli space is regular, these cases
    avoid numerous other technical difficulties. For example, the
    sign of Equation (\ref{su3rho}) is constant for
    these manifolds. This goes back to a result of Fintushel and
    Stern which identifies the parity of the $su(2)$--spectral flow
    of irreducible flat $SU(2)$ connections on Brieskorn spheres
    (see \cite{FS} as well as the proof of Theorem \ref{2pcorrect}).  In
    the $SU(3)$ case, we know from \cite{boden} that the
    $su(3)$--spectral flow is {\it even}, which implies that $\la'_{SU(3)}(X_{\pm
    1})$ is given by simply counting the irreducible $SU(3)$
    representations of $\pi_1 (X_{\pm 1}).$ Moreover for $\pi_1(X_{\pm
    1})$,  all the important irreducible $SU(3)$ representations
    can be described in terms of representations of {\it finite}
    groups.

We compute the $SU(3)$ Casson invariant $\la_{SU(3)}$ for $\pm1$ surgery
    on the right hand  trefoil.   Recall that (at least as unoriented
    manifolds) $X_{+1}\cong\Si(2,3,5)$ and $X_{-1}\cong \Si(2,3,7)$.

\begin{example} Consider first $X_{+1}.$ Case 1 of Example  \ref{rhorho1}
    shows that the moduli space of flat $SU(2)$ connections on $X_{+1}$
    is $$ \fM_{SU(2)}(X_{+1}) = \{ [\Th], [A_1], [A_2] \},$$ where
    $\varrho_{X_{+1}}(A_1) = 59/30$ and $\varrho_{X_{+1}}(A_2)=
    131/30.$

 For any irreducible flat connection
 $A$ on $ X_{+1}$,  $SF_{su(2)}(\Theta,A;X_{+1}) $ is odd
 (see Theorem \ref{2pcorrect}).
Thus  Theorem \ref{corrected} implies $$
    \la''_{SU(3)}(X_{+1})  =
-\tfrac{1}{2}(\varrho_{X_{+1}}(A_1)+\varrho_{X_{+1}}(A_2))=-19/6.
$$
Since $\pi_1 (X_{+1})$ is the binary icosahedral group (which is finite),
    it is well-known that it has two irreducible rank 3 representations.
    (Setting $\al_i = \hol_{A_i}$ for $i=1,2$, these are the two $SU(3)$
    representations obtained by   composing  $\al_i \co  \pi_1 (X_{+1}) \to
    SU(2)$ with the sequence of maps $SU(2) \to SO(3) \hookrightarrow
    SU(3)$ given by the standard projection followed by the natural
    inclusion.)  Proposition 5.1 of \cite{boden} shows that the adjoint
    $su(3)$--spectral flow of $A_1$ and $A_2$ is  even, hence
    $\la'_{SU(3)}(X_{+1})=2$. Hence $$\la_{SU(3)}(X_{+1}) =
    \la'_{SU(3)}( X_{+1}) + \la''_{SU(3)}( X_{+1}) =2-19/6 =-7/6.$$
    \end{example}

\begin{example} Now consider $X_{-1}.$ Case 2 of Example  \ref{rhorho1}
    shows that the moduli space of $SU(2)$ connections on $X_{-1}$
    is $$\fM_{SU(2)}(X_{-1}) = \{ [\Th], [A_1], [A_2] \},$$ where
    $\varrho_{ X_{-1}}(A_1) =  37/42$ and $\varrho_{ X_{-1}}(A_2)=
    109/42.$ In this case, we know from Casson's invariant
    that the $su(2)$--spectral flow
    has the opposite parity from
    the previous case, namely $SF_{su(2)}(\Th,A_i;M_{-1})$
    is even for $i=1,2$. 
     Theorem \ref{corrected} implies
    \begin{equation} \label{M-1correct} \la''_{SU(3)}(X_{-1}) =
    \tfrac{1}{2}(\varrho_{ X_{-1}}(A_1)+\varrho_{ X_{-1}}(A_2)) =
     73/42.\end{equation}

In \cite{boden} it is shown that there are four irreducible $SU(3)$
    representations of $\pi_1 (X_{-1}).$  Two of these are obtained
    from the $SU(2)$ representations $\al_i = \hol_{A_i}$ as in the
    previous example and the other two are induced by representations of
    the  quotient $PSL(2, {\mathbb F}_7)$, which is a finite group
    of order $168$, as follows. Comparing the group presentations
    \begin{eqnarray*} \pi_1 (X_{-1})  &=& \langle x,y,z,h \mid h
    \hbox{ central}, x^2 h = y^3 h^{-1} = z^7 h^{-1} = xyz
=1\rangle \\
    PSL(2, {\mathbb F}_7) &=& \langle x,y,z \mid
    x^2=y^3=z^7=xyz=[y,x]^4=1 \rangle,\end{eqnarray*} it is clear
    that $ PSL(2, {\mathbb F}_7)$ is the quotient of $\pi_1 (X_{-1}) $
    by the normal subgroup $\langle h, [y,x]^4 \rangle.$ 

It is
    well-known that $ PSL(2, {\mathbb F}_7)$ has precisely two
    irreducible $SU(3)$ representations (see page 96 of 
    \cite{coxeter-moser}), thus the two  remaining $SU(3)$
    representations of $\pi_1 (X_{-1}) $ are obtained from
    $PSL(2,{\mathbb F}_7)$ by pullback. As before, Proposition 5.1
    of \cite{boden} implies that the adjoint $su(3)$--spectral flow
    of each of the four irreducible $SU(3)$ representations is
    even. Hence $\la'_{SU(3)}(X_{-1})=4$. Using this and Equation
    (\ref{M-1correct}), it follows that $$\la_{SU(3)}(X_{-1}) =
    4+73/42=241/42.$$\end{example}

\noindent Since  the $SU(3)$ Casson invariant is unchanged by a change
    of orientation,  we  conclude that $\la_{SU(3)}(\Si(2,3,5)) =
    -7/6$ and $ \la_{SU(3)}(\Si(2,3,7)) = 241/42.$

\section{Computations for torus knots} \label{torusknot}

Given a 3--manifold $X$ and a flat $SU(2)$ connection $A$ on it,
    Theorems \ref{nicerformula} and \ref{csthm}   determine the
    spectral flow and the Chern--Simons invariant of $A$ provided
    there exists a knot $K$ in $X$ so that $A$ is connected to the
    trivial connection $\Th$ by a path of flat connections on the
    knot complement $Z=X-N(K)$. In this section we apply
    our methods to perform explicit computations
    for homology spheres obtained
    by surgery on a torus knot.
    Computing the spectral flow on the complement of a torus knot
    is not hard, and it is especially straightforward for
     $(2,q)$ torus knot complements
    (see Theorem \ref{compofS}).
    In this way, we reduce the computation of
     all the gauge theoretic invariants,
    including the rho invariants,  to
    straightforward computations of the integers $a,b,c$,
   and the integral $2\int nm'$.

Our aim is to compute the $SU(3)$ Casson invariant
for surgeries on torus knots. In the general case,
 one needs to consider perturbed flat connections since
the $SU(2)$ representation variety may not be  cut out transversely
as a subspace of the $SU(3)$ representation variety.
For surgeries on $(2,q)$ torus knots, transversality holds
    so one can
   compute the $SU(3)$ Casson invariant
   without perturbing. As in
  Theorem \ref{corrected}, this has the happy consequence
    that the correction term $\la''_{SU(3)}$ can be expressed entirely
    in terms of the rho invariants. Using this approach,
    we compute $\la''_{SU(3)}$
    for homology spheres  obtained by surgery on  a $(2,q)$ torus
    knot. Coupling these results with the computations of
    $\la'_{SU(3)}$ in \cite{boden}, we   calculate $\la_{SU(3)}$
    for surgeries on $K(2,q)$ for various $q$ and use this data
    to conclude that
   $\la_{SU(3)}$ is not a finite type invariant of
    order $6$.

\subsection{Twisted cohomology of torus knot complements}
    \label{cohtorusknot}

We begin with a discussion of orientations and surgery
    conventions.   Any  knot $K$ in $S^3$ induces  a
    decomposition of $S^3$ into two pieces, a tubular neighborhood
    $N(K)$ and the knot exterior $Z= S^3-N(K)$.  This decomposition
    uniquely determines two isotopy classes of {\it unoriented}
    simple closed curves on  the separating torus: $\tmu$ is a
    curve which bounds a disc in $N(K)$ and $\tla$ is a  curve that
    bounds a surface in $Z$.  These can be represented by smooth
    curves intersecting transversely in one point.
Orient the pair $\{\tmu, \tla\}$ so that $\tmu \cdot \tla=1$ with the
outward normal first boundary orientation on $T=\partial (N(K))$. 

For any integer $k$,  consider the homology sphere
   $X=Y\cup_T Z$ obtained by performing $\tfrac{1}{k}$ surgery on $K$.
    This is the 3--manifold obtained by
gluing the solid torus $Y = D^2 \times S^1$ to $Z$ using a
diffeomorphism of their boundaries which takes
$\partial D^2 \times \{1\}$ to $\tmu  \tla^k$ and
$\{1\} \times S^1$ to $\tla.$  The curves

$$
 \mu = \tmu \tla^k, \ \la = \tla
$$
are called the meridian and longitude of the Dehn filling $Y$.
Notice that $\mu$ {\it does not} represent the usual
    meridian for the trefoil as a knot in $S^3$, but
    rather the meridian for the knot in $X$ which is the core of the
    Dehn filling.

  Now consider, for $p$ and $q$ relatively prime and positive,
 the $(p,q)$
    torus knot  $K(p,q)$.  This is  the knot
      $K\co [0,2\pi]\to
\RR^3$ parameterized by
$K(t)= ((2+ \cos qt) \cos pt, (2+ \cos qt) \sin pt, - \sin qt).$
The restriction that  $p$ and $q$ be positive is not serious; all the
    methods presented here are equally valid if either $p$ or $q$
    is negative. Notice however that  the result of $\frac{1}{k}$
    surgery on $K(p,q)$ is diffeomorphic to that of $-\frac{1}{k}$
    surgery on $K(p,-q)$ by an orientation-reversing diffeomorphism.
    Since the gauge theoretic invariants change in a predictable
    way under reversal of orientations, there is no loss in generality
    in assuming that $p$ and $q$ are positive.

The exterior $Z=S^3-N(K)$ of the $(p,q)$ torus knot $K$ has fundamental
    group $$\pi_1 Z  = \langle x,y \mid x^p=y^q \rangle.$$ Choose
    $r,s \in \ZZ$ with the property that $pr+qs = 1.$
Then the curves $\tmu$ and $\tla$ for $K(p,q)$ are represented in $\pi_1 Z $
    as\begin{equation} \label{natmerlong} \tmu = x^s
    y^r\quad \hbox{and} \quad \tla = x^p ( \tmu )^{-pq}.
    \end{equation}

For example, for the $(2,q)$ torus knot, one can take $s=1$ and
    $r=(1-q)/2$. Then $ \tmu = x y^{ (1-q)/2}$ and $ \tla = x^2(
    \tmu)^{-2q}.$

Theorem 1 of \cite{klassen-thesis} gives a general description of the
    variety of $SU(2)$ representations of torus knot groups. There
    it is shown that for a $(p,q)$ torus knot, $\fR_{SU(2)}(Z)$ is
    a connected, 1--dimensional singular manifold (smooth except for
    `T' type intersections, called $SU(2)$ bifurcation points, discussed below.)
     Figures \ref{trefoilrep},  \ref{repexamples}, 
     and \ref{repexamplespq} illustrate these representation varieties
     for several different torus knots.

Since $\tmu$ normally generates $\pi_1 Z,$
any reducible representation $\be\co \pi_1 Z \to SU(2)$
is uniquely determined by $\be(\tmu)$. Throughout this
section, we adopt the notation where $\be_s$ for $s \in [0,\frac{1}{2}]$
refers to
the reducible representation of $\pi_1 Z$ which is uniquely determined
up to conjugacy by the requirement that  $\be_s(\tmu) = e^{2 \pi i s}$. 
Since $\tla$ lies in the commutator subgroup of $\pi_1Z$ (it bounds a
Seifert surface) the reducible representation $\be_s$ sends $\la$ to $1$
so $\be_s(\tmu)=\be_s(\mu)$. 

The space $\fR^*_{SU(2)}(Z)$ of irreducible representations consists of
 $(p-1)(q-1)/2$ open arcs, the ends of which
    limit to  distinct  reducible representations.  Thus,   $\fR_{SU(2)}(Z)$
    is the space obtained by identifying the endpoints of a collection
    of closed arcs with distinct interior points of
     the interval $[0,\tfrac{1}{2}]$.  It follows that
     $\fR_{SU(2)}(Z)$  is path connected and
     any flat  connection $A$ on a homology sphere obtained from surgery
    on a torus knot  can be connected to  the trivial connection
    $\Theta$ by a path of  connections which are flat on $Z$
    and which satisfy conditions 1--3 of Subsection
    \ref{defnofAt}.

The next result is crucial to computations of $SF(A_\eta(t);Z;P^-)$ for
    torus knot complements.  It identifies the kernel of $D_A$ with
    $P^-$ boundary conditions at any flat connection on $Z$.

\begin{theorem} \label{pqthm}  {\sl Let $Z$ be the exterior of any
    $(p,q)$ torus knot and suppose  $\al\co \pi_1(Z)\to SU(2)$ is a
    representation,  defining a local coefficient system in $\CC^2$.
    \begin{enumerate}   \item[\rm(i)] If $\al$ is trivial, then
    $H^{0+1}(Z,\partial Z;\CC^2_\al)=0$.   \item[\rm(ii)] If $\al$
    is nontrivial,  then $H^0(Z,\partial Z;\CC^2_\al)=
    H^0(Z;\CC^2_\al)=0$ and $$H^1(Z,\partial Z;\CC^2_\al)=
    H^1(Z;\CC^2_\al)=\begin{cases} \CC^2& \text{if $\al(x^p)=1$
    and $\al(x) \neq 1 \neq \al(y)$,}\\ 0&\hbox{otherwise.}
    \end{cases} $$\end{enumerate} In particular if $A$ is a flat
    connection on $Z$ with nontrivial holonomy $\al=$hol$_A$ then
    the kernel of $D_A$ with $P^+$ boundary conditions  is isomorphic to
    $\CC^2$ if $\al(x^p)=1$, $\al(x)\ne 1$, and $\al(y)\ne 1$, and
    this  kernel is zero otherwise.}\end{theorem}
    \noindent {\bf Proof}\qua  The
    first statement is an easy exercise in cohomology since the
    coefficients are untwisted.   The chain complex for the universal
    cover of $Z$ is computed by the Fox Calculus to be (with
    $\pi=\pi_1(Z)$)  
    \begin{equation}\label{chcom}  0\to \ \ZZ[\pi]\
    \mapright{d_2} \ \ZZ[\pi] \oplus \ZZ[\pi] \ \mapright{d_1}\ZZ[\pi] \
    \to 0\end{equation} with $$d_1=\left[\begin{array}{c} x-1\\
    y-1\end{array}\right]$$ and $$d_2=\left[\begin{array}{cc} 1+x+
    \cdots + x^{p-1} \ &\  x^p y^{-q}(1+y+y^2+\cdots +y^{q-1})
    \end{array} \right].$$ \

The cohomology $H^*(Z;\CC^2_\al)$ is defined to be the homology of the
    chain complex obtained by applying Hom$_{\ZZ[\pi]}(-,\CC^2)$ to
    the complex (\ref{chcom}), where $\pi$ acts on $\CC^2$ via
    $\al$. The differentials are obtained by replacing $x$ and $y$
    in the matrices $d_2$ and $d_1$ by $\al(x)$ and $\al(y)$ and
    taking the transpose.   Denote by $d_1(\al)$ and
$d_0(\al)$ the
    resulting differentials.

Clearly $d_0(\al)=0$ if and only if $\al$ is trivial.  Thus  if $\al$ is
    nontrivial, then $H^0(Z;\CC^2_\al)=0$.  On the other hand, if
    $\al$ is nontrivial and $\al(x)=1$   then $\al(y)$ must be a
    nontrivial $q$--th root of unity and it follows that $d_1(\al) =
    \left[\begin{matrix}p\\q\end{matrix}\right] $ and $d_0(\al) = [   0 \quad
    {\al(y)-1}]$. Since $\al(y)\neq 1,$   we see that $\ker d_1(\al) =
    \im d_0(\al)$. This implies $H^1(Z;\CC^2_\al)=0$.

Similar arguments apply and give the same conclusion if $\al$ is
    nontrivial and $\al(y)=1.$

So assume $\al(x) \neq 1$ and $\al(y) \neq 1.$ This implies that $ \im
    d_1(\al)$ has complex dimension 2. If $\al(x^p) =1$, then $\al(x)$
    is a nontrivial $p$--th root of unity, which implies $1+\al(x)+\cdots
    + \al(x^{p-1}) =0.$ Since $x^p = y^q,$ it follows also that $\al(y)$
    is a nontrivial $q$--th root of unity. Thus $d_2(\al)$ is the
    zero map and so $H^1(Z;\CC^2_\al) = \ker d_1(\al)/ \im d_0(\al)$ is
    isomorphic to $\CC^2$ in this case.

On the other hand, if $\al(x^p) \neq 1,$ then $\al(x)$ is not a $p$--th
    root of unity. Hence $d_1(\al)$ is not the zero map and this
    forces $H^1(Z;\CC^2_\al)=0.$

We have seen (Equation (\ref{cohomolofT})) that if $\al$ restricts
    nontrivially to $T=\partial Z$, then the cohomology of $T$ vanishes.
    Of course,  since the meridian normally generates $\pi_1 Z,$
    any nontrivial representation $\al$ of $\pi_1 Z$ pulls back to
    a nontrivial representation of $\pi_1 T$. Now the long exact
    sequence in cohomology shows that $H^i(Z;\CC^2_\al)=H^i(Z,\partial
    Z;\CC^2_\al)$.
The last statement follows from Proposition \ref{exact}.\endproof

The characterization of the representation varieties $\fR_{SU(2)}(Z) $
    of torus knot groups in \cite{klassen-thesis} shows that
    $\hol_{A_t}(x)$ is   constant (up to conjugacy) along paths of
irreducible
    representations.  Let $A$ be a flat connection on $Z$ with
    nontrivial holonomy and set $\al = \hol_A$.  Use Proposition
    \ref{exact} to  identify $H^{1}(Z; \CC^2_\al)$
    with the kernel of $D_A$ on $Z$ with $P^-$
    boundary conditions. Then Theorem
    \ref{pqthm} implies the following result.

\begin{theorem} \label{vanish} {\sl If $A_t$ is a path of
    irreducible flat $SU(2)$ connections on the complement $Z$ of a
    $(p,q)$ torus knot, then the dimension of $H^1(Z;\CC^2_{A_t})$
    is independent of $t$. In particular $SF(A_t;Z;P^-)=0.$  Moreover,
    if $p=2$ then $\dim H^1(Z;\CC^2_{A_t})=0$ for all $t$.}
    \end{theorem}

\subsection{Jumping points and SU(2) bifurcation points}
As we have already seen, the representation variety $\fR_{SU(2)}(Z)$
of  the complement
of a torus knot can be described as the
space obtained by identifying endpoints of $(p-1)(q-1)/2$ closed
arcs with interior points in the line segment
$[0,\frac{1}{2}]$ parameterizing the reducibles.
The non-smooth points, which are precisely where
 the arcs are attached, are called {\it $SU(2)$ bifurcation points}.
It is a simple matter to
   characterize the  $SU(2)$ bifurcation points (for torus knots) in terms of the
Alexander
   polynomial, although which pairs of bifurcation points are endpoints of
the same arc of irreducibles is a more subtle problem.
    However, if $p=2,$ the answer is simple because
    the arcs are glued in a nested  way  (see
    Proposition    \ref{2qsfpts}).

    Closely related to the $SU(2)$ bifurcation points are the $\CC^2$ jumping  
    points. These play a central role in determining
    $SF(A_\eta;Z;P^-),$ the spectral
    flow along the knot complement.

\begin{definition} \label{defnsfbif}
Suppose $Z$ is the complement of a knot $K$ in a
    homology sphere $X$.\begin{enumerate}   \item[(i)]
    The {\it $\CC^2$ jumping points} are the gauge orbits of nontrivial
    reducible flat $SU(2)$
    connections $A$ on $Z$ where the kernel of $D_A$  with $P^-$
    boundary conditions jumps up in dimension.  This is the set of
reducible flat connections $A$ so that    
$ H^1(Z;\CC^2_A)
\neq 0$.
    \item[(ii)] The {\it $SU(2)$ bifurcation points} are gauge orbits of
    reducible flat connections $A$ on $Z$ which are limits of
    irreducible, flat connections.\end{enumerate}\end{definition}

Theorem \ref{vanish} implies that for torus knot complements, only certain
    reducible flat connections are $\CC^2$ jumping points. They
    are characterized in terms of the roots of the Alexander polynomial
    by the following theorem, which is reminiscent of the
    characterization of the $SU(2)$ bifurcation points in terms of square
    roots of the Alexander polynomial \cite{klassen-thesis}.

\begin{theorem}\label{corsing}  {\sl Suppose $K$ is the $(p,q)$
    torus knot and $Z$ is its exterior. Given a reducible, flat,
    nontrivial $SU(2)$
    connection $A$, the kernel of $D_{A}$ with $P^-$ boundary
    conditions  is nontrivial if and only if $\hol_A(\tmu)$ is a root of
    the  Alexander polynomial } $$\Delta_K(t) =
    \frac{(t^{pq}-1)(t-1)}{(t^p-1)(t^q-1)}.$$

\end{theorem}

\noindent {\bf Proof}\qua Let $\al=\hol_A$ be the reducible representation of $\pi_1 Z
    $ associated with $A.$ By Theorem \ref{pqthm}, $ H^{1}(Z, \partial
    Z; \CC^2_\al) \neq 0$ if and only if $\al(x)$ has order $p'$
    and $\al(y)$ has order $q'$ for $p' \neq 1 \neq q'$, where $p'$
    divides $p$ and $q'$ divides $q$.   Since $p'$ and $q'$ are
    relatively prime, this implies that $\al(x)$ and $\al(y)$ generate a
    cyclic group of order $p'q'$ equal to  $\langle \al(\tmu) \rangle$.
    Thus $\hol_A( \tmu)\neq 1$ is a $pq$--th root of unity, but it
    is neither a $p$--th root  nor a $q$--th root of unity.\endproof

For torus knots, the roots of
$\Delta_K(t)$ all lie on the unit circle. Thus $\CC^2$ jumping points
correspond to the reducible representations
$\be_t$ where $e^{2 \pi i s} $ is a root of $\Delta_K$.
Since we always use $[0,\frac{1}{2}]$ to parameterize reducible
$SU(2)$ representations,  we will  simply say  $s \in [0,\frac{1}{2}]$
is a $\CC^2$ jumping point
if $\be_s$ is. Similarly, we say that $s$ is an $SU(2)$ bifurcation point
if $\be_s$ is.

We give
two quick examples. First, if $K$ is the trefoil then
    $\Delta_K(t)= t^2-t+1$. Its roots are $e^{\pm \pi i/3}$.
    This gives one  $\CC^2$ jumping point at
    $s={1/6}$ (this is the black dot in Figure \ref{trefoilrep}).
    By  \cite{klassen-thesis},  the $SU(2)$ bifurcation points
    are the solutions to  $\Delta_K(t^2) =0$.
    So for the trefoil,  ${1/12}$ and
    ${5/12}$ are the $SU(2)$ bifurcation points. (These are just the square roots
    of the $\CC^2$ jumping point.)
    Next,
    if $K$ is the $(2,5)$ torus knot, then $\Delta_K(t) = t^4-t^3+t^2-t+1$.
   Its  roots are $e^{\pm \pi i/ 5}$ and $e^{\pm 3\pi i/ 5}$, yielding
    two $\CC^2$ jumping points at ${1/10}$ and ${3/10}$.
    There are four $SU(2)$ bifurcation points:
    $\{ \frac{1}{20}, \frac{3}{20}, \frac{7}{20}, \frac{9}{20}\}$.
Generalizing to $(2,q)$ torus knots, we obtain the following proposition.

\begin{proposition} \label{2qsfpts}  {\sl Suppose $Z$ is the
    exterior of a $(2,q)$ torus knot and consider its $SU(2)$
    representation variety $\fR_{SU(2)}(Z).$ 
    Parameterize the reducible representations by
    $[0,\frac{1}{2}]$ as above. Then there are $(q-1)/2$ arcs of 
irreducible
    representations, indexed  as $\hR_\ell$ for
    $\ell=1,\ldots, (q-1)/2$, such that  $\hR_\ell$ is attached to
    $[0,\frac{1}{2}]$ at the $SU(2)$ bifurcation points $\frac{2\ell-1}{4q}$ and
    $\frac{1}{2} - \frac{2\ell-1}{4q}.$ Thus the arcs of irreducible
    representations are nested.

Notice further that the  set of $SU(2)$ bifurcation points 
$$\{ \tfrac{2\ell-1}{4q},  \tfrac{1}{2} - \tfrac{2\ell-1}{4q}
\mid \ell=1,\ldots, (q-1)/2\}$$ is disjoint from the set 
$$\{  \tfrac{2\ell-1}{2q} \mid \ell=1,\ldots, (q-1)/2\}$$  of $\CC^2$ jumping points.}

\end{proposition}

\noindent {\bf Proof}\qua  
      The Alexander
    polynomial of $K(2,q)$ is $\Delta_K(t) = \frac{t^{q}+1}{t+1}$
    and its roots are the $q$--th roots of $-1$ different from $-1.$
    Theorem \ref{corsing} easily identifies the $\CC^2$ jumping points
    as the set $\{
    \tfrac{2\ell-1}{2q} \mid \ell=1,\ldots, (q-1)/2\}$.
The $SU(2)$ bifurcation points correspond to the
    reducible representations $\be$ with $\be(\tmu^2)$ a root of
    $\Delta_K(t)$  \cite{ klassen-thesis}.  Taking
square roots  gives the set
    $\{ \tfrac{2\ell-1}{4q},  \frac{1}{2} - \frac{2\ell-1}{4q}
\mid \ell=1,\ldots, (q-1)/2\}$ of $SU(2)$ bifurcation points.

The arcs $\hR_\ell$ can be described as follows. Since $\pi_1 Z =
    \langle x,y \mid x^2=y^q \rangle$  has infinite cyclic center
    generated by $x^2,$  if $\al$ is irreducible then $\al(x^2) = -1$
    (the centralizer of any nonabelian subgroup of $SU(2)$ is $\pm1$).
    Conjugating if necessary, we have $\al(x)=i$. Similarly, $\al(y)$
    can be conjugated to lie in the $ij$--plane and must be a $q$--th
    root of $-1.$  Drawing the great circles from $1$ to $\al(x)$
    and from $1$ to $\al(y)$,  one can use the angle between the
    great circles
    to parameterize the arcs $\hR_\ell$ as in the proof of Theorem
    1 in \cite{klassen-thesis}.

We will be even more specific. For  $\ell = 1, \ldots, (q-1)/2$
 define  a 1--parameter family of representations $\al_{\ell,t} \co  \pi_1 Z
    \to SU(2)$ for $t\in [0,1]$ by   setting $\al_{\ell,t}(x)=i$ and
    $$\al_{\ell,t}(y)=\cos(\tfrac{(2\ell-1) \pi}{q}) +
    \sin(\tfrac{(2\ell-1) \pi}{q})  (i \cos (\pi t) + j \sin(\pi
    t)).$$  These representations are irreducible except at the
    endpoints. The meridian of $K(2,q)$ is represented in $\pi_1 Z$
    by $\tmu = xy^{(1-q)/2}$,   and a simple computation shows that the
    endpoints of $\al_\ell$ are   the reducible representations $\be_s$ for
    $s\in \left\{\frac{q-2 \ell-2}{4q}, \frac{1}{2}- \frac{q-2
    \ell-2}{4q} \right\}$. As $\ell$ ranges from $1$ to $(q-1)/2$, the
    associated pairs of points in $[0,\frac{1}{2}]$ are nested
    ($\al_{\ell,t}$ parameterizes the arc $\hR_{(q-1)/2-\ell}$.)
    \endproof

\begin{figure}[ht!]\small \begin{center}
\includegraphics[scale=.45]{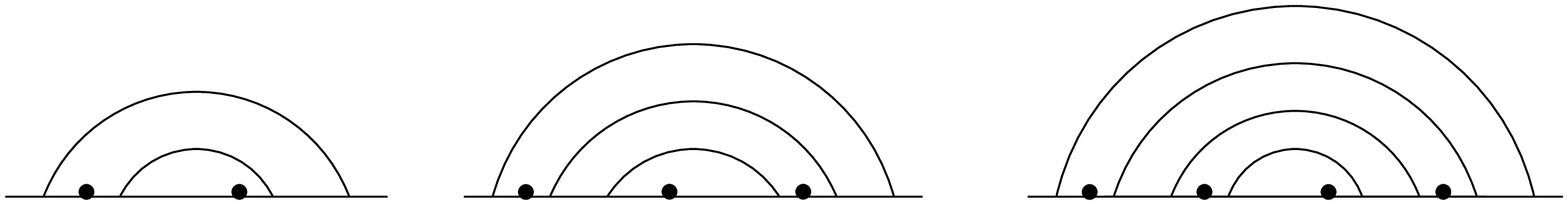} 
\end{center} 
\vskip-.6cm \hskip1.55cm $K(2,5)$\hskip2.45cm $K(2,7)$ \hskip3.3cm $K(2,9)$ 
 \vskip .1cm
 \caption{$SU(2)$ representation varieties of $(2,q)$ torus knot groups}
\label{repexamples} 
\end{figure}

\begin{figure}[ht!]\small 
\begin{center}
\includegraphics[scale=.70]{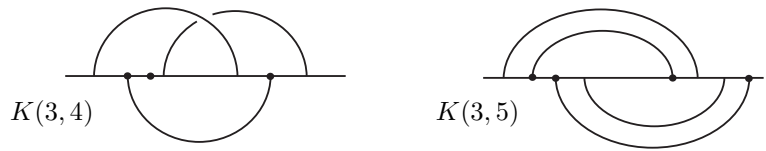}     
\end{center} 
\vskip-.9cm \hskip1.2cm $K(3,4)$\hskip4.5cm $ K(3,5)$
\vskip.3cm
\caption{$SU(2)$ representation varieties of $(p,q)$ torus knot groups}
\label{repexamplespq} 
\end{figure}

Figures \ref{repexamples} and
\ref{repexamplespq} show the $SU(2)$ representation varieties for
    several torus knots. The horizontal line segment denotes the
    reducibles, and the $\CC^2$ jumping points are dots. The   curved
    arcs are the irreducible components $\hR_\ell$.

Comparing Figures \ref{repexamples} and
\ref{repexamplespq}, one sees that most of the assertions of Proposition
    \ref{2qsfpts} fail for arbitrary $(p,q)$ torus knots. The sets
    of $SU(2)$ bifurcation points and $\CC^2$ jumping points are not in general
    disjoint. Nor are the irreducible components nested as for $K(2,q)$.
    For $K(p,q)$, one can still use Theorem \ref{corsing} to
    characterize the $\CC^2$ jumping points, they occur at the
   reducible representations $\be_s$ where $\be_s(\tmu)=e^{2 \pi i s}$ is a
    $pq$--th root  of unity which is neither a $p$--th root nor a
    $q$--th root  of unity.

\begin{corollary} \label{2qreg} {\sl Suppose $X$ is a homology sphere
    obtained by surgery on a $(2,q)$ torus knot and   $A$ is a nontrivial
    flat $SU(2)$ connection on $X.$  Then $H^i(X;\CC^2_A) = 0
    = H^i(X;su(2)_A)$ for all $i$. Hence, $\fM^*_{SU(2)}(X)$ is regular
    as a subset of $\fM_{SU(3)}(X)$ and the correction term $\la''_{SU(3)}(X)$
    can be computed
    in terms of $SU(2)$ rho invariants.}\end{corollary}
   \noindent {\bf Proof}\qua
    Since $X$ is a homology sphere, any nontrivial
    flat $SU(2)$ connection is irreducible. The restriction of such a
connection
    to $Z$ is  irreducible
    (since $\pi_1(Y)\cong\ZZ$), and its restriction to $T$ is   nontrivial
since the meridian $\tmu$ normally generates $\pi_1Z$.
   Theorem \ref{vanish}
    implies that $H^1(Z;\CC^2_{A})=0$, and the corollary
    follows from the computations of
    Subsection \ref{compsofco} using the Mayer--Vietoris sequence.
(The vanishing of first cohomology with both $su(2)$ and $\CC^2$ 
coefficients is the definition of regularity in \cite{bh}.)
    \endproof
 The corollary  reflects a rather special property of $(2,q)$ torus
 knots.  If $p,q>2$ and $X$ is a homology sphere obtained by surgery on $K(p,q)$
then there   exists an irreducible $SU(2)$ representation which is the
limit of an arc  of irreducible $SU(3)$ representations, and hence
$\fM^{*}_{SU(2)}(X)$ is not regular (when viewed as a subset of
$\fM_{SU(3)}(X)$), although it is a compact 0--dimensional manifold.

    Results of Fintushel and Stern shows that, for homology
    spheres $X$ obtained by surgery on a torus knot, $\fR^*_{SU(2)}(X)$
    is a finite set of points and the parity of
    $SF_{su(2)}(\Th,A;X)$ is independent of $[A] \in \fM^*_{SU(2)}(X)$.
    Corollary \ref{2qreg} and Theorem \ref{corrected} then imply
    that, for surgeries on a $(2,q)$ torus knot,
     $$\la''_{SU(3)}(X)
    = \pm {\tfrac{1}{2}} \sum_{[A]\in \fM^*_{SU(2)}(X)} \varrho_X(A).$$
     The next result determines the sign.  Recall if $\Delta_K(t)$
denotes the symmetrized Alexander polynomial of the knot $K$, then
$\Delta_K''(1)$ equals twice Casson's invariant of the knot (\cite{am}).  For
the torus knot $K=K(2,q)$, $\Delta''_K(1)=(1-q^2 )/4$.

\begin{theorem} \label{2pcorrect} 
{\sl Let $K$ be the $(2,q)$ torus knot.  Assume $k>0$ and denote by $X_{\pm k}$
    the result of $\pm \frac{1}{k}$ surgery on $K.$ By Corollary
    \ref{2qreg}, the $\fM_{SU(2)}(X_{\pm k})$ is regular as a subspace of
    $\fM_{SU(3)} (X_{\pm k})$.  If $b= (q^2-1)/4$ 
then the moduli spaces $\fM_{SU(2)}^*(X_{k})$
and $\fM_{SU(2)}^*(X_{-k})$ consist of $kb$
    points.

    \begin{enumerate}   \item[\rm(i)]
    Writing $\fM_{SU(2)} (X_{ k}) = \{ [\Th],[A_1],\ldots, [A_{kb}] \}$, 
then
    $SF_{su(2)}(\Th,A_i;X_{k})$ is odd and\begin{equation}
    \label{eq2pcorrect+} \la''_{SU(3)}(X_{ k}) = -{\tfrac{1}{2}}
    \sum_{i=1}^{kb}    \varrho_{X_k} (A_i).\end{equation}

  \item[\rm(ii)] Writing $\fM_{SU(2)} (X_{- k}) = \{ [\Th],[A_1],\ldots, [A_{kb}]
    \}$, then $SF_{su(2)}(\Th,A_i;X_{-k})$ is even and\begin{equation}
    \label{eq2pcorrect-} \la''_{SU(3)}(X_{-k}) = {\tfrac{1}{2}}
    \sum_{i=1}^{kb}    \varrho_{X_{-k}} (A_i).\end{equation}
    \end{enumerate} } \end{theorem}

    \noindent {\bf Proof}\qua
  Suppose  $X$ is a  homology sphere obtained
  by a positive surgery on $K(p,q).$   We claim that
  if $A$ is an irreducible flat
  $SU(2)$ connection on $X$, then $SF_{su(2)}(\Th,A;X)$ is odd. Theorem
\ref{corrected} then implies the first assertion.

By Taubes'
  theorem \cite{taubes} and the surgery formula for Casson's invariant,
  this implies that $SF_{su(2)}(\Th,A;X)$ is even in case
  $X$ is obtained by a negative surgery on $K(p,q)$, so the second
assertion follows from the first.  

Since Casson's knot invariant  
  is positive for all $(2,q)$ torus knots, it suffices to check one
example,
  which we take to be the Poincar\'e homology sphere.

   First, consider a path $A_t$ of $SU(2)$  connections on a
    homology sphere $X$ with $A_i$ flat for $i=0,1$.
    Let $\al_i = \hol_{A_i}$ for $i=0,1$ be the associated $SU(2)$
representations.
    Applying the Atiyah--Patodi--Singer
    index theorem as in the proof of Theorem \ref{rhoinvts} (cf.~Equation
    (\ref{rhoAPS})), we see that
\begin{eqnarray} \label{rhoadAPS}
SF_{su(2)}(A_t;X) &=& 8 (cs(A_1) -cs(A_0))+
\frac{1}{2}\left( \varrho_X({\rm ad}\, \al_1) -\varrho_X({\rm ad}\,\al_0)
\right) \nonumber \\ && 
 \hskip -.7in + \ \frac{1}{2}\left(\dim H^{0+1}(X; su(2)_{\al_1})-
\dim H^{0+1}(X; su(2)_{\al_0})\right), 
\end{eqnarray}
    where the rho invariants are defined relative to the odd signature
    operator acting on $su(2)$--forms via the adjoint action.
    (Comparing this to the formula at the end of Section 7 in
     \cite{KKR}, the  sign discrepancies are explained by the
    fact that in this paper, we are using the $(-\ep,-\ep)$ convention for
spectral
    flow.)

    Now consider the Poincar\'e homology sphere $X$, defined here to
    be $+1$ surgery on the right hand trefoil. Consider a path
    from $\Th$ to the flat connection $A_1$ constructed in Section
    \ref{trefoil}. Since $\pi_1 X$ is finite,
    one can compute that $\varrho_X({\rm ad}\,\al_1)  = \pm \frac{73}{15}$
    by standard methods \cite{FS, auckly}.
     The sign ambiguity
    comes about because the answer depends on how $X$ is oriented.
    This problem can be resolved  since we know
     that $cs(A_1) = \frac{1}{120}$ by the computations in Subsection
\ref{trefoil}.
     Applying Equation (\ref{rhoadAPS}) to the path $A_t$,
     we see that the only way the left hand side can be an integer
     is if $\varrho_X({\rm ad}\,\al_1) = \frac{73}{15}$, in which case
     $SF_{su(2)}(\Th,A_1;X)=1$, which is odd, as claimed.
  \endproof

\noindent{\bf Remark}\qua  It is well-known
     that, at least as unoriented manifolds, $\tfrac{1}{k}$-surgery on a
torus knot yields a Brieskorn homology sphere. These spaces admit a
natural orientation as the link of an algebraic singularity. Using
standard handlebody methods  (or alternatively, using the fact that the
$su(2)$--spectral flow from the trivial connection to a flat $SU(2)$
connection on a Brieskorn sphere is even \cite{FS}) it follows that, as
oriented manifolds,
$$X_k \cong - \Si(2,q,2qk-1) \quad \text{ and } \quad X_{-k} \cong
\Si(2,q,2qk+1).$$


 Theorems \ref{2pcorrect} and \ref{rhoinvts}   give a method for
 computing $\la''_{SU(3)}$
    for surgeries on $(2,q)$ torus knots.  Combining this with the computations
    of $\la'_{SU(3)}$ given in   \cite{boden}, we shall determine
    $\la_{SU(3)}$ for homology spheres obtained by surgery on
    $K(2,q)$.
    The analogous computation for $K(p,q)$
    is complicated by the fact that
    one must first apply a perturbation  to make the moduli space
    regular, so we postpone the calculations for surgeries on other
    torus knots to a future article.

Consider the complement $Z$  of $K(p,q)$, with $\partial Z=T$, the torus. 
Recall that
$\fR_{SU(2)}(T)$ denotes the variety of conjugacy classes of $SU(2)$
representations of $\pi_1T$.  Any choice of generators
$x,y\in\pi_1T=\ZZ\oplus\ZZ$ determine a branched cover $\RR^2\to
\fR_{SU(2)}(T)$ by assigning to the pair $(m,n)\in\RR^2$ the conjugacy
class of the homomorphism  taking $x$ to $e^{2\pi im}$ and $y$ to
$e^{2\pi i n}$.  (See Equation (\ref{pillowpara}).) We will need to consider
the cases
$x=\tmu\tla^k$ and
$y=\tla$ for various $k$ simultaneously. Thus we introduce the
notation
$$
f_k\co \RR^2\to \fR_{SU(2)}(T)
$$
for the map taking $(m,n)$ to the conjugacy class of the homomorphism
$\alpha\colon$\break 
$\pi_1T\to SU(2)$ satisfying $\alpha(\tmu\tla^k)=e^{2\pi i m}$ and
$\al(\tla)=e^{2\pi i n}$.  Letting $g_k\co \RR^2\to \RR^2$
denote the linear map
\begin{equation}\label{emsubk}
g_k(m,n)=(m+kn,n) 
\end{equation}
we see that  $$f_0= f_k\circ g_k.$$

Now consider the restriction map $\fR^*_{SU(2)}(Z) \to \fR_{SU(2)}(T)$. Each
component   of $\fR^*_{SU(2)}(Z)$ is an open arc.  Proposition
\ref{2qsfpts} enumerates these in the special case of $K(2,q)$; the path
components are denoted  $\hR_\ell$.   The image  in
$\fR_{SU(2)}(T)$ of each arc misses the branch points since  the branch
points correspond to central representations of $\pi_1T$, but $\tmu$ 
cannot be sent to the center $\pm 1$ by an irreducible (ie, nonabelian)
representation since it is
a normal generator of $\pi_1Z$.  Thus each path component of 
$\fR^*_{SU(2)}(Z)$ lifts to
$\RR^2$.

We claim that any such lift using the cover $f_0\co \RR^2\to \fR_{SU(2)}(T)$ 
takes the components of  $\fR^*_{SU(2)}(Z)$ to  arcs  of slope
$-pq$. One can see this as follows.
If $\al_t\co \pi_1Z\to SU(2)$ is any continuous path of
    irreducible representations, then it can be conjugated
    so that
    $\al_t(\tmu)$ and $\al_t(\tla)$ lie in the standard $U(1)$ subgroup
    of $SU(2)$.  If $\al_t$ is irreducible, then $\al_t(x^p) = \pm 1$.
Writing $\al_t(\tmu)=e^{2 \pi i  m_t}$
and $\al_t(\tla) = e^{2 \pi i  n_t}$,
then Equation (\ref{natmerlong}), namely that $\tla = x^p (\tmu)^{-pq}$, 
shows
 that $pq \,  m_t + n_t$ is constant.

To understand how the arcs in $\fR^*_{SU(2)}(Z)$ lift using 
 the
cover $f_k\co \RR^2\to  \fR_{SU(2)}(T)$, that is, with respect to $\mu$ and
$\la$ for the homology sphere obtained by $\tfrac{1}{k}$ surgery on
$K(p,q)$, one just applies the map $g_k$ of Equation (\ref{emsubk}). Thus 
each arc lifts using $f_k$ to arcs in $\RR^2$ of slope
$\tfrac{pq}{kpq-1}$.

The following proposition completes the identification of the lifts of each
arc $\hR_\ell$ for the $(2,q)$ torus knots.

 \begin{proposition}\label{kDsurg} {\sl Let $K$ be the $(2,q)$  torus knot
 and $Z$ its complement.
    For $\ell = 1,\ldots, (q-1)/2,$  consider the curve $R_\ell(t)$ defined for
    $ 0<t<1$ by setting
    $$R_\ell(t) = (1-t) \left(\tfrac{2 \ell -1}{4q},0\right) +  t
    \left(\tfrac{1}{2} -\tfrac{2 \ell -1}{4q} +k(2\ell-q-1), \,
    2\ell-q-1\right). $$
     Then $R_\ell$ is the lift under $f_k\co \RR^2 \to \fR_{SU(2)}(T)$ of
     $\hR_\ell\subset \fR^*_{SU(2)}(Z)$  (see Proposition \ref{2qsfpts}).
      All other lifts of $\hR_\ell$ are
    obtained by reflecting this lift through the origin and/or
    translating by an integer vector.}

\end{proposition}

\noindent {\bf Proof}\qua
We first determine the lift with respect to the  map
$$f_0\co \RR^2
\to 
\fR_{SU(2)}(T).$$
It was shown in the paragraph preceding this proposition 
that any lift of
$\hR_\ell$ has slope $-2q$.   Proposition \ref{2qsfpts} shows that the
endpoints of the arc
$\hR_\ell$ are reducible representations sending $\tmu$ to $e^{2\pi
i(2\ell-1)/4q}$  and
$e^{2\pi i(1/2-(2\ell-1)/4q)}$  and $\tla$ to $1$.  Thus there is a lift of
$\hR_\ell$ starting at $(\tfrac{2\ell-1}{4},0)$ and ending at
$ (e(\tfrac{1}{2}-\tfrac{2\ell-1}{4})+a,b)$ for some integers $(a,b)$ and
$e=\pm 1$.

We claim that $a=0$ and $e=1$. Assuming this for a moment,   the fact
that the slope is $-2q$ implies that $b=2\ell-q-1$. Thus the curve
$$
  (1-t) (\tfrac{2 \ell -1}{4q},0) +
t (\tfrac{1}{2} -\tfrac{2 \ell -1}{4q},2\ell-q-1), \ \ 0 \leq t \leq 1
$$
parameterizes the lift using $f_0$ of $\hR_\ell$ based at $(\frac{2 \ell
-1}{4q},0).$
Applying the map $g_k$ of Equation (\ref{emsubk}) finishes the proof.

It remains to show that $a=0$ and $e=1$. Suppose not. Then the lift of
$\hR_\ell$ to $\RR^2$ intersects one of the  vertical lines $x=0$ or
$x=\tfrac{1}{2}$. This means there exists a representation
$\alpha\in
\hR_\ell\subset
\fR^*_{SU(2)}(Z)$ so that
$\alpha(\tmu)=\pm 1$. But $\tmu$ normally generates $\pi_1Z$ and $\pm
1\in SU(2)$ is the center, so $\alpha$ is central, contradicting the fact
that $\al$ is irreducible.  
\endproof

\subsection{Dehn surgeries on the trefoil} \label{moretrefoil}
  In this subsection, we compute the gauge theoretic
    invariants  for flat connections on the manifolds $X_{\pm k}$ obtained
    by   $\pm \frac{1}{k}$
    surgery on the right hand trefoil $K$.

We consider the cases of positive and negative surgeries separately to make
counting arguments simpler in Theorems \ref{23formulas}, \ref{-23formulas},
\ref{+2qfull}, and \ref{-2qfull}. The reason for this is that the slopes of
the curves $R_\ell$ are positive for $k>0$ and negative for $k<0$, changing
the combinatorics of the numbers $a,b$ and $c$.  We combine the separate
results in the computations of the $SU(3)$ Casson invariant, so Theorems 
\ref{23su3} and Tables \ref{table3} and \ref{table4} are valid for all integers $k$ (including
$k=0$).

\begin{theorem}\label{23formulas} {\sl Suppose $k>0$ and denote
    by $X_{k}$ the result of $  \frac{1}{k}$ surgery on the right hand
trefoil $K$.
    Then $\pi_1 (X_{ k})$ admits $2k$ distinct conjugacy  classes
    of irreducible $SU(2)$ representations. In terms of  the moduli
    space of flat connections, this gives $$\fM_{SU(2)}(X_k) = \{
    [\Th],[A_1],\ldots, [A_{2k}]\}.$$ For $i=1,\ldots, 2k,$ we can
    choose $A_i$ a representative for the gauge orbit $[A_i]$ with
    \begin{eqnarray*} SF(\Th,A_i;X_k) &=&
   2-2i+2\left[\tfrac{i}{k+1}\right]\\
    cs(A_i) &=& 2-2i +(2k-2i+2)
    \left[\tfrac{i}{k+2}\right]  + \frac{(12i-11)^2}{24(6k-1)}\\
     \varrho_{X_k}(A_i)  &=&   4i-2   
    +4\left[\tfrac{i}{k+1}\right]  + 8(i-k-1)
    \left[\tfrac{i}{k+2}\right]  - \frac{ (12i-11)^2}{6(6k-1)}.
    \end{eqnarray*} Here, $[x]$ means the greatest integer less than or
    equal to $x$.}\end{theorem}

\noindent {\bf Proof}\qua By Proposition  \ref{kDsurg}, the lift of the one arc of
    irreducible representations from   $\fR_{SU(2)}(T)$ to $\RR^2$ is
    given by the curve $$R_t= (1-t)(\tfrac{1}{12},0) + t(
    \tfrac{5}{12}-2k,-2), \ \  0 \leq t \leq 1.$$

The flat connections which extend over $\frac{1}{k}$ surgery correspond
    to points along the path where the first coordinate $
    (1-t)\tfrac{1}{12}+ t(\tfrac{5}{12}-2k)$ is an integer. Let $A_i$ be
    the $i$--th such point along the arc $R_t$. Let
$t_i$ be the
    corresponding $t$ value.  Since $k>0$ we see that $t_i$ solves
    the equation $(1-t)\tfrac{1}{12}+ t (\tfrac{5}{12}-2k)=1-i$ for
    $i=1,\ldots,2k$ and so $$t_i = \frac{12i-11}{24k-4}, \quad
    i=1,\ldots, 2k.$$

Fix $i\in \{1,\ldots,2k\}$. Then there is a path of flat connections $C_t$ in
    normal form on $Z$ so that  the restriction to the torus $T$ is
    $a_{m,n}=-m_t i dx - n_t i dy$ with $(m_t,n_t)$ the composite of the
horizontal
    line segment   from $(0,0)$ to
    $(\tfrac{1}{12},0)$  with $R_t$, ending at $R_{t_i}$.

From this path we compute the integers $a_i,b_i$ and $c_i$ and construct
    the flat connection $A_i$ on $X_k$ and the path $A_t$ of connections
    on $X_k$ starting at the trivial connection and ending at $A_i$
    according to the method of Subsection \ref{constrpt}. (We hope
    the clash of notation $A_t|_{t=t_i}=A_i$ does not cause too much
    confusion. The integer $i$ is fixed throughout the rest of the
    argument.)

By definition, $a_i=1-i$ and
$$b_i=[-2t_i]=\left[-\tfrac{12i-11}{12k-2}\right]
=-1-\left[\tfrac{i}{k+1}\right].$$ Inspecting the graph  of the path
$(m_t,n_t)$ one can compute
   that 
$$ c_i =
    2i-2+(2i-2k-2)  \left[\tfrac{i}{k+2}\right].$$
To see this, 
observe that the 
loop constructed in Subsection
\ref{constrpt} encloses the lattice points $(1-j,0)$, $j=1,\cdots,i-1$. If
$i\ge k+2$, it also encloses the lattice points
$(1-j,-1)$, for $j=k+1,\ldots,i-1$. Since the loop
winds around all the lattice points clockwise, it follows that
$c_i=2(i-1) + 2(i-k-1) \left[\tfrac{i}{k+2}\right]$. 
   
 Now Theorem \ref{nicerformula}   
    implies that $SF(\Th,A_i) = 2-2i+2[\tfrac{i}{k+1}]$ because
    the spectral flow $SF(C_t;Z;P^-)$ along the knot complement vanishes.

To compute $cs(A_i)$, notice that the integral term $\int m'n$ in
    Theorem \ref{csthm} vanishes along the first part of the path
    (since $n_t=0$ along that part). On the second part one computes $$2
    \int_0^{t_i} m' n = ({\tfrac{1}{12}} -(1-i))(2 t_i)  = \frac{
    (12i-11)^2}{24(6k-1)},$$ and substituting this into the formula
    of Theorem \ref{csthm} gives $$cs(A_i) = 2-2i+(2k-2i+2)
   \left[\tfrac{i}{k+2}\right] + \frac{(12i-11)^2}{24(6k-1)}.$$

   Theorem \ref{rhoinvts} (or alternatively,
   Equation (\ref{rhorho})) gives the formula for the rho  invariants.
    \endproof

For negative surgeries, we get the following analogous result.

\begin{theorem} \label{-23formulas} {\sl Suppose $k>0$ and let
    $X_{-k}$ denote $-\frac{1}{k}$ surgery  on the right hand trefoil.
    Then $\pi_1 (X_{- k})$ admits $2k$ distinct conjugacy  classes
    of irreducible $SU(2)$ representations. In terms of  the moduli
    space of  flat connections, $$\fM_{SU(2)}(X_{-k})= \{ [\Th],[A_1],\ldots,
    [A_{2k}]\}.$$ For $i=1,\ldots, 2k,$ we can choose $A_i$ a
    representative for the gauge orbit $[A_i]$ with
    \begin{eqnarray*}
    SF(\Th,A_i;X_{-k}) &=&  2i+2\left[\tfrac{i}{k+1}\right]\\
    cs(A_i) &=&  2i+(2i-2k) \left[\tfrac{i}{k+1}\right] -
\frac{(12i-1)^2}{24(6k+1)} \\
    \varrho_{X_{-k}}(A_i) &=& 2-4i+ 4(2k-2i+1)
    \left[\tfrac{i}{k+1}\right]  + \frac{ (12i-1)^2 }{6(6k+1)}.
    \end{eqnarray*}} \end{theorem}
    \noindent {\bf Proof}\qua This theorem is proved using a similar argument as
    was used for positive Dehn surgery. The main difference is  that now
    the second part of path $(m_t,n_t)$ is given by $$
    (1-t)(\tfrac{1}{12},0) + t( \tfrac{5}{12}+2k,-2).$$  This path
    has first coordinate the integer $i\in\{1,\ldots,2k\}$ when $$t_i =
    \frac{12i-1}{24k+4}.$$

From the definitions, $a_i = i$ and $$ b_i =
[-2t_i]=\left[-\tfrac{12i-1
}{12k+2}\right]=-1-\left[\tfrac{i}{k+1}\right].$$

One
can compute from the graph of the path $(m_t,n_t)$ using a similar
analysis as in the proof of Theorem \ref{23formulas} that
    $$c_i = -2i+(2k-2i) \left[\tfrac{i}{k+1}\right].$$
    These determine as in Subsection \ref{constrpt}  a path of
    connections on $X_{-k}$ from the trivial connection to  a
    connection $A_i$ extending flatly
    over $X_{-k}$. 

As before, the spectral flow along $Z$
    vanishes.
The integral term is computed as
 $$2 \int_0^{t_1} m' n = -(i-{\tfrac{1}{12}})(2t_i) =
    -\frac{(12i-1)^2}{24(6k+1)}.$$
The proof is then completed by applying
Theorems \ref{nicerformula},
    \ref{csthm} and  \ref{rhoinvts}.\endproof

The following theorem   gives a general computation of the Casson
    $SU(3)$ invariant $\la_{SU(3)}$ for surgeries on the trefoil.

\begin{theorem} \label{23su3} {\sl  For any integer $k$ let $X_{k}$
    denote the homology sphere obtained by $\frac{1}{k}$ surgery on
    the trefoil. Then

\begin{eqnarray*}
    \la_{SU(3)}(X_{k}) &=&\frac{k(84k^2 -138k+19)}{6(6k- 1)}.
    \end{eqnarray*}}
    \end{theorem}

\noindent {\bf Proof}\qua Consider first the case $k>0.$  The results in Section
    5 of \cite{boden} show that $\la'_{SU(3)} (X_k) = 3k^2- k$.
    Using this and Equation (\ref{eq2pcorrect+}) and summing the
    rho invariants from Theorem  \ref{23formulas}, we see
    that\begin{eqnarray*} \la_{SU(3)}(X_k)  &=& \la'_{SU(3)}(X_k)
    + \la''_{SU(3)}(X_k) \\ &=&  3k^2-k - {\tfrac{1}{2}} \sum_{i=1}^{2k}
    \varrho_{X_k}(A_i) \\ &=& 3k^2-k \\
    &-& {\tfrac{1}{2}} \sum_{i=1}^{2k}
    \left(4i-2+4 \left[{\tfrac{i}{k+1}}\right]  + 8(i-k-1) \left[
    \tfrac{i}{k+2} \right] - \frac{(12i-11)^2}{6(6k-1)}
\right).\end{eqnarray*}

Using that
\begin{eqnarray*}
 \sum_{i=1}^{2k} 4 \left[{\tfrac{i-1}{k}}\right] +
8(i-k-1)
\left[\tfrac{i}{k+2} \right] 
&=& \sum_{i=k+1}^{2k}4 + \sum_{i=k+2}^{2k} 8(i-k-1) \\
&=& 4k + 4k^2 - 4k = 4k^2
\end{eqnarray*}
and standard summation formulas, we see that
$$\la_{SU(3)}(X_k) = \frac{k(84k^2 - 138k+19)}{6(6k - 1)}.$$

The proof for the case $k<0$ is similar, using Theorems \ref{-23formulas}
    and Theorem \ref{eq2pcorrect-}, and yields the same formula. 
    \endproof

Inspecting this proof one sees that the terms involving the greatest integer
function  in the sum defining $\la''_{SU(3)}(X_k)$ for $k>0$ contribute a
quadratic polynomial in $k$ to $\la_{SU(3)}(X_k)$, and the remaining terms 
contribute a rational function whose numerator is cubic in $k$ and whose
denominator is $6(6k-1).$  A
perfectly analogous computation in the case of
$K(2,q)$ treated below shows that the
$SU(3)$ Casson invariant of
$\tfrac{1}{k}$ surgery on $K(2,q)$ will always be a rational
function with cubic numerator and denominator
$2q(2qk-1)$ for
$k>0$. Similarly the
$SU(3)$ Casson invariant of
$-\tfrac{1}{k}$ surgery on $K(2,q)$ will always be a rational
function with cubic numerator and denominator
 $2q(2qk+1)$.

\subsection{Dehn surgeries on (2,q) torus knots} \label{2qcomp} 
     In this subsection, we  compute the
    spectral flow and the Chern--Simons invariants
    for flat connections
    on homology spheres obtained by surgery on a $(2,q)$ torus knot.
    We also determine the
    correction term $\la''_{SU(3)}$ by  summing the
    rho invariants and applying Theorem \ref{2pcorrect}.

The main  difference, which is illustrated in Figure \ref{Clfor(2,5)}
(see also Figure \ref{repexamples}),
 is that  the spectral flow
$SF(A_\eta; Z;P^-)$ along the knot complement need not vanish
as it did for the complement of the trefoil.
For example,  for $\tfrac{1}{k}$--surgery
    on $K(2,5)$, the two lifts $R_1$ and $R_2$ of the image
    of
    $\fR^*_{SU(2)}(Z) \to \fR_{SU(2)}(T)$  are separated
    by a $\CC^2$ jumping point.

\begin{figure}[ht!]\small
\begin{center}
\includegraphics[scale=0.8]{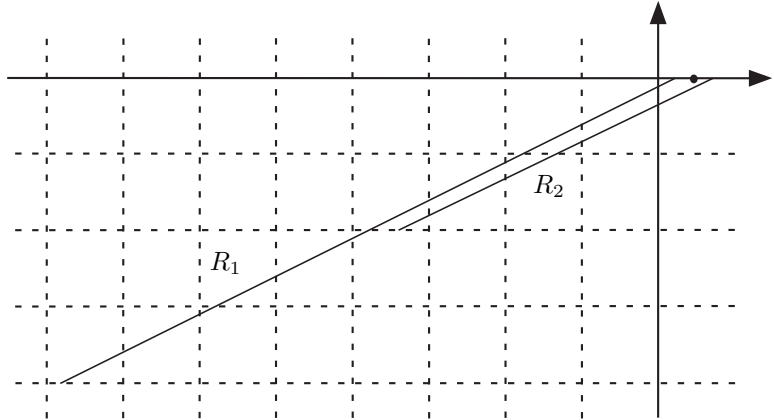} 
\end{center}
\vskip-3.5cm\hskip8.5cm $R_2$
\vskip.6cm \hskip4.2cm $R_1$
\vskip 1.8cm
\caption{$+\frac{1}{2}$ surgery on $K(2,5)$}
\label{Clfor(2,5)}
\end{figure}

  \begin{theorem}\label{compofS} {\sl Suppose $1 \leq \ell \leq (q-1)/2$
  and let $C_t$
  be a path of flat reducible connections
  on $Z$ in normal form such that
 $ {C_t}|_T = \frac{t+\ell-1}{q} i dx$ for $t \in [0,1]$. Notice that
 $C_t$ crosses one and only one $\CC^2$ jumping point
 (the one at $\frac{2\ell-1}{2q}$).
  Then
  $$SF(C_t;Z;P^-)=2.$$}
   \end{theorem}

    \noindent {\bf Proof}\qua
First notice that $-2\leq SF(C_t;Z;P^-)\leq 2$. This is because the
kernel of $D_{C_t}$ on $Z$ with $P^-$ boundary conditions is
$2$--dimensional at the jumping points, and $0$--dimensional  at  
non-trivial  reducible connections.  Thus two eigenvalues become zero at
the jumping point.

 We prove the theorem by comparing the rho invariant for gauge
    equivalent flat connections on
    the manifold $X_{+1}$ obtained by $+1$ surgery on $K(2,q)$.

    The path
    $$R_1(t) = (1-t)(\tfrac{1}{4q},0) +
t(\tfrac{1}{2}-\tfrac{1}{4q}+1-q,1-q)$$
(see Proposition \ref{kDsurg})    crosses the vertical axis at $t_0 =
\frac{1}{(q-1)(4q-2)}$.
    Let $(m_t,n_t)$
    be the composition of the short horizontal segment
    from $(0,0)$ to $(\frac{1}{4q},0)$  with the
    path $R_1(t)$ for $0\leq t \leq  t_0$ and let $A_t$ be the path
     of  connections on $X_{+1}$ which are flat along $Z$
     and correspond to the path $(m_t,n_t)$
     by the construction of Subsection \ref{constrpt} of $X_{+1}$.
    From this path, we compute  that $a =0, b = -1, c= 0$ and
   that  $2\int m'n = \frac{1}{4q(4q-2)}.$ Since $A_t$ misses
   all the $\CC^2$ jumping points, $SF(A_t;Z;P^-) = 0$
   and Theorem \ref{rhoinvts} implies
\begin{eqnarray*}
\varrho_{X_{+1}}(A_1) &=& 4(a-b+c)-2-8 \int m'n\\
    &=& 2-\frac{1}{q(4q-2)}.
    \end{eqnarray*}

   Now consider
    the path obtained by translating $R_1(t)$ by the vector $(q-1,q-1)$.
Proposition \ref{kDsurg} implies that this is another lift to $\RR^2$
of the arc $\hR_1\subset \fR^*_{SU(2)}(Z)$.
    Parameterized in the opposite direction (so that it starts on
    the horizontal axis), this is the curve
    $$\tR(t) = (1-t)(\tfrac{1}{2} - \tfrac{1}{4q},0) +
    t(\tfrac{1}{4q}+q-1, q-1).$$

    This crosses the vertical line   $x=q-1$ when
    $\tit_0 = 1 - t_0 =\frac{4q^2-6q+1}{4q^2-6q+2}$.
    Let $(\tm_t,\tn_t)$
    be the composition of the short horizontal segment
    from $(0,0)$ to $(\frac{1}{2} - \frac{1}{4q},0)$  with $\tR(t)$
    for $ 0 \leq t \leq \tit_0.$
    The corresponding path $\tA_t$ of connections
       crosses each of the
        $\CC^2$ jumping points exactly once and ends at
  $\tA_1$, which  is gauge equivalent to $A_1.$  
    Using the path, we compute as before   that
    $\tilde{a}=q-1$ and $ \tilde{b}=q-2$. To compute $\tilde{c}$,   observe
that $\tR(t)$ intersects the horizontal line $y=i$ in the point $(x_i,i)$
with $i<x_i<i+1$ if $ i=1,\cdots,q-1$.  Thus the loop constructed in 
Subsection
\ref{constrpt} encloses no lattice points of the form $(1,n)$, one lattice
point of the form $(2,n)$ (namely $(2,1)$), and in general encloses $j-1$
lattice points of the form $(j,n)$.  Hence in total, the loop encloses  
$1+2+\cdots +(q-2)=\tfrac{(q-2)(q-1)}{2}$ lattice points. These  are all
enclosed clockwise, so
$ \tilde{c}=\tfrac{(q-2)(q-1)}{2}=q^2-3q+2$.

The integral $2\int m'n $ is equal to $ \frac{(4q^2-6q+1)^2}{4q(4q-2)}.$
    Hence
    \begin{eqnarray*}
    \varrho_{X_{+1}}(\tA_1) &=& 4(\tilde{a}-\tilde{b}+\tilde{c})-2 - 8 \int
m'n +2 SF(\tA_t;Z;P^-) \\
    &=& 4-2q - \frac{1}{q(4q-2)} + 2 SF(\tA_t;Z;P^-).
    \end{eqnarray*}

    Since $A_1$ and $\tA_1$ are gauge equivalent, their rho
    invariants are equal.
    Setting $ \varrho_{X_{+1}}(A_1) = \varrho_{X_{+1}}(\tA_1)$
    and solving for $SF(\tA_t;Z;P^-)$ gives
    $$SF(\tA_t;Z;P^-) = q-1.$$
     Since there are
    exactly $\frac{q-1}{2}$ $\CC^2$--jumping points, the path $\tA_t$ passes
through all  of them, and each contributes at most $2$ to the spectral flow, 
the    spectral flow across each one
    is $  2.$ This proves the theorem.
    \endproof

The following lemma will be useful in simplifying formulas. 
\begin{lemma}\label{tricky} {\sl Let $q,k,\ell,i$ be positive integers with $q\ge 3$,
 $\ell\leq \frac{q-1}{2}$, and $i\leq k(q-2\ell+1)$. Let $[x]$
be the greatest integer less than or equal to $x.$  Then
$$\left[\frac{4q(1-i)-2\ell+1}{4qk-2}\right]=\left[-\frac{i}{k}\right] =
\left[\frac{2\ell-4qi -1}{4qk+2}\right].$$}
\end{lemma}
\noindent {\bf Proof}\qua
Letting $x=\tfrac{4q(1-i)-2\ell+1}{4qk-2}$,
one can easily check that
$0<x+\tfrac{i}{k}<\tfrac{1}{k}$.
  This implies that
$[x]=[-\tfrac{i}{k}]$.  Similarly, letting
$y=\tfrac{2\ell-4qi -1}{4qk+2}$, one checks
that $0<y+\tfrac{i}{k}<\tfrac{1}{k}$, which implies
$ [y]=[-\tfrac{i}{k}]$. 
\endproof

We can now turn our
attention to computing the gauge theoretic invariants.
Suppose $k>0$ and let $X_{\pm k}$ denote the manifold obtained by
$\pm \frac{1}{k}$ surgery
on $K(2,q).$
 By
Proposition \ref{kDsurg}, the curves
 $$R_\ell(t) =
    (1-t) \left(\tfrac{2 \ell -1}{4q},0\right) +  t \left(\tfrac{1}{2}
    -\tfrac{2 \ell -1}{4q} \pm k(2\ell-q-1), \, 2\ell-q-1\right),\ \
    0<t<1  $$
    for $\ell = 1, \ldots, (q-1)/2$ are lifts of the restrition map $\fR^*_{SU(2)}(Z)
\lto \fR_{SU(2)}(T)$ under the branched cover $f_{\pm k}\co  \RR^2 \lto
\fR_{SU(2)}(T)$.

Consider first the case of positive surgeries. We would like to
    determine the flat connections  which extend over
    $X_k$ for $k>0.$  Fixing $\ell,$ these correspond to  points on
$R_\ell(t)$ whose
    first coordinate is an integer.
This happens when $$t_i = \frac{4q(1-i) - 2 \ell +1}{(4qk-2)(2 \ell - q
    - 1)},$$ in which case the first coordinate of $R_\ell(t)$ is
    $1-i$, with $i\in \{1,\ldots,k(q-2\ell+1)\}$.

Fix $\ell$ and $i$
with $1 \leq \ell \leq \frac{q-1}{2}$ and
$1 \leq i \leq k(q-2\ell+1)$.
 Define the path $(m_t,n_t)$ to be the composition of
   the horizontal line
    segment  from  $(0,0)$  to $(\tfrac{2\ell-1}{4q},0)$ with the path
    $R_\ell(t)$ for   $t\in [0,t_i]$.
   Let $A_t$ be the  path of connections corresponding to $(m_t,n_t)$
    by the construction of Subsection \ref{constrpt}. The endpoint, which
we denote
      by $A_{\ell,i},$
      extends flatly over $X_{k}$.
Denote the integers $a,\ b,\ c$ associated to $A_{\ell,i}$ by
$a_{\ell,i},\  b_{\ell,i},\  c_{\ell,i}$.
Then, using Lemma \ref{tricky} one sees that
\begin{eqnarray*} a_{\ell,i} &=& m_{t_i} = 1-i\\ 
   b_{\ell,i} &=& [n_{t_i}] =[t_i(2\ell-q-1)]= \left[\frac{4q(1-i)-2
    \ell+1}{4qk-2}\right] = \left[ -\frac{i}{k}\right].  \end{eqnarray*}
Inspecting the graph of $(m_t,n_t)$ one sees  that
$c_{\ell,i}-c_{\ell,i-1}=2(-b_{\ell,i-1})$ and $c_{\ell,1}=0$, so
$$c_{\ell,i} =  -2
\sum_{j<i}
    b_{\ell,j} =  -2 \sum_{j=1}^{i-1} \left[-\frac{j}{k}\right].$$

To calculate the
    Chern--Simons invariant, we compute the integral:
\begin{eqnarray*}
   2\int m' n &=& 2\int_0^{t_i} [(k-\tfrac{1}{2q})(2\ell-q-1)]
    ( 2 \ell - q -1 ) t \  dt \\  &=& \frac{(4q(1-i) - 2\ell+1)^2}{4q(4qk-2)}.
    \end{eqnarray*}
Since the 
horizontal line segment from $(0,0)$ to $(\tfrac{2\ell-1}{4q},0)$
(ie, the first part of the path) passes
through the   $\CC^2$ jumping points at
$\tfrac{1}{2q},\tfrac{3}{2q},\cdots,\tfrac{2[\ell/2]-1}{2q}$, Theorem
\ref{compofS} implies that $SF(C_t;Z;P^-)=2[\tfrac{\ell}{2}]$.

    Applying Theorems \ref{nicerformula}, \ref{csthm} and \ref{rhoinvts},
    we compute the spectral flow, the Chern--Simons invariants and
    the rho invariants of $A_{\ell,i}$. The results are
    summarized in the following theorem.

    \begin{theorem} \label{+2qfull}
  {\sl Suppose $k>0$ and let
    $X_k$ be the result of $\frac{1}{k}$ surgery on the $(2,q)$
    torus knot. Let $A_{\ell,i}$ for $\ell=1,\ldots,  {(q-1)}/{2}$
    and $i = 1, \ldots, k(q+1-2 \ell)$ be the flat connections on $X_k$
    constructed above. Then
    \begin{eqnarray*}
 SF(\Th,A_{\ell,i};X_k) &=&  2  \left[\frac{\ell}{2}\right] -2i-
 2\left[-\frac{i}{k}\right]  \\
 cs(A_{\ell,i})  &=&  \frac{(4q(1-i) - 2\ell+1)^2}{4q(4qk-2)}
 + 2 \sum_{j=1}^{i-1} \left[-\frac{j}{k}\right] \\
 \varrho_{X_k}(A_{\ell,i}) &=&  4 \left[\frac{\ell}{2}\right] +  2-4i
   - \frac{(4q(1-i) - 2\ell+1)^2}{q(4qk-2)}  \\
   && -4\left[-\frac{i}{k}\right] - 8 \sum_{j=1}^{i-1} \left[-\frac{j}{k}\right].
 \end{eqnarray*}}\end{theorem}

Now consider the situation for negative surgeries on
    $K(2,q)$.    We would like to
    determine the flat connections  which extend over
    $X_{-k}$. (To make counting arguments simpler we still assume $k>0$).
    Fixing $\ell,$ these correspond to  points on $R_\ell(t)$ whose
    first coordinate is an integer.
This happens when $$t_i = \frac{4qi-2\ell + 1}{(q-2 \ell+1)(4qk+2)},$$
 in which case the first coordinate of $R_\ell(t)$ is
    $i$, with $i\in \{1,\ldots,k(q-2\ell+1)\}$.

    Fix $\ell$ and $i$
    with $1 \leq \ell
    \leq  (q-1)/2$ and  $1 \leq i \leq k(q-2\ell+1).$
    Define the path $(m_t,n_t)$ to be the composition
    of the horizontal line from $(0,0)$ to $(\frac{2\ell-1}{4q},0)$
    with $R_\ell(t)$ for $t \in [0,t_i].$
    Let $A_t$ be the  path of connections corresponding to $(m_t,n_t)$
    by the construction of Subsection \ref{constrpt}. The endpoint, which
we denote
      by $A_{\ell,i},$
      extends flatly over $X_{-k}$.
    We compute
     the numbers $a_{\ell,i}, \ b_{\ell,i}, \ c_{\ell,i}$ associated
     to  $A_{\ell,i}$. First,
 \begin{eqnarray*} a_{\ell,i} &=& m_{t_i} = i  \\
b_{\ell,i}
    &=& [n_{t_i}] =
[t_i(2\ell-q-1)]=\left[\frac{2\ell-4qi-1}{4qk+2}\right] =\left[-\frac{i}{k}\right]  
   \end{eqnarray*}
using Lemma \ref{tricky}.
Inspecting the graph of $(m_t,n_t)$ one sees  that
$$c_{\ell,i}=2 \sum_{j \leq i} b_{\ell,j} = 2 \sum_{j=1}^i
    \left[ -\frac{j}{k}\right].$$
Finally,  
    $$2\int_0^{t_i}  m' n = -\frac{(4qi-2\ell+1)^2}{4q(4qk+2)}.$$  

Just as in the case of positive surgery, the first part of the
path passes through the $\CC^2$ jumping points at
$\tfrac{1}{2q},\tfrac{3}{2q},\cdots,\tfrac{2[\ell/2]-1}{2q}$ and thus
$SF(C_t;Z;P^-)=2[\tfrac{\ell}{2}]$.
   
Theorems
\ref{nicerformula}, \ref{csthm} and \ref{rhoinvts} then give formulas for
the spectral flow, the Chern--Simons invariants, and
    the rho invariants for all connections on $X_{-k}$.

\begin{theorem} \label{-2qfull}
   {\sl Suppose $k>0$ and let
    $X_{-k}$ be the result of $-\frac{1}{k}$ surgery on the $(2,q)$
    torus knot. Let $A_{\ell,i}$ for $\ell=1,\ldots, \frac{q-1}{2}$
    and $i = 1, \ldots, k(q+1-2 \ell)$ be the flat connections
    constructed above. Then
    \begin{eqnarray*}
 SF(\Th,A_{\ell,i};X_{-k}) &=& 2  \left[\frac{\ell}{2}\right] + 2i-2 -
 2\left[-\frac{i}{k}\right]   \\
  cs(A_{\ell,i})  &=& -\frac{(4qi-2\ell+1)^2}{4q(4qk+2)}-2 \sum_{j=1}^i
    \left[ -\frac{j}{k}\right]  \\
  \varrho_{X_{-k}}(A_{\ell,i}) &=& 4\left[\frac{\ell}{2}\right]-2 + 4i 
  +\frac{(4qi-2\ell+1)^2}{q(4qk+2)} \\
  &&-4
\left[-\frac{i}{k}\right]+8 \sum_{j=1}^i
    \left[- \frac{j}{k}\right].
\end{eqnarray*} } \end{theorem}

Summing the rho invariants and applying Theorem
\ref{2pcorrect} yields the correction term
  $\la''_{SU(3)}$ for any homology sphere obtained by
  surgeries on a $(2,q)$ torus knot.  The results are summarized
    in Table \ref{table3}. (The computations of $\la'_{SU(3)}$
    can be found in \cite{boden}.)  
    
    \medskip

\begin{table}[ht]\small 
\renewcommand\arraystretch{2}
 \noindent\[
  \begin{array}{|c|cc|}     \hline 
& \la'_{SU(3)}(X_k) & \la''_{SU(3)}(X_k)  \\ \hline
 \   K(2,3) \ & 3k^2 - k
      & {\displaystyle \frac{-24k^3 - 84k^2+13k}{6(6k- 1)} }\\
\    K(2,5) \ & 33k^2 - 9k
      & {\displaystyle \frac{-200k^3 - 1620k^2+151k}{10(10k- 1)}} \\
 \   K(2,7) \ & 138k^2 - 26k
      & {\displaystyle \frac{-784k^3 - 9128k^2+606k}{14(14k- 1)}} \\
 \ K(2,9)\ & \ 390k^2 - 58k     \
      & \ {\displaystyle \frac{-2160k^3 - 33192k^2+1714k}{18(18k- 1)}}    \
      \\   \hline
      \end{array}  \]
    \caption{{$ \la'_{SU(3)}$ and $\la''_{SU(3)}$ for homology spheres $X_k$
    obtained by  $\tfrac{1}{k}$ surgery
    on  $K(2,q)$}}
    \label{table3}
     \end{table}

     \medskip
     
In completing this table we used the following fact.
Fix $q$ and let $X_k$ denote the manifold obtained by $\tfrac{1}{k}$
surgery on $K(2,q)$. As noted after the proof of Theorem \ref{23su3},
for positive $k$  the quantity
$2q(2q-1)\la''_{SU(3)}(X_k)$ is a cubic polynomial in $k$.
  Hence one can deduce $\la''_{SU(3)}(X_k)$ for all $k$ by computing it
in several examples and solving for the coefficients. Similar methods apply
if $k$ is negative.  

The entries in this table are valid
for any integer $k$, not just $k>0$. Despite the slight differences in the
statements and proofs of Theorems \ref{+2qfull} and \ref{-2qfull}, after
summing over all $A_i$, the resulting formulas give the same rational
function. For
$k=0$ the homology sphere is $S^3$ which has $SU(3)$ Casson invariant $0$
since it is simply connected. 

By summing $\la'_{SU(3)}$ and $\la''_{SU(3)},$
we  compute the $SU(3)$ Casson
    invariants for homology 3--spheres $X_k$ obtained
    by $\tfrac{1}{k}$
surgery on $K(2,q)$.

 \medskip
     
\begin{table}[ht]\small 
\renewcommand\arraystretch{2}
 \noindent\[
  \begin{array}{|c|c|}     \hline 
  & \la_{SU(3)}(X_k) \\ \hline
 \   K(2,3) \ & {\displaystyle \frac{84k^3 - 138k^2+19k}{6(6k- 1)}} \\
\    K(2,5) \ & {\displaystyle \frac{3100k^3  - 2850k^2+241k}{10(10k- 1)}} \\
 \   K(2,7) \ &  {\displaystyle \frac{26264k^3 - 16156k^2+970k}{14(14k- 1)}} \\
 \ K(2,9)\ &\  {\displaystyle \frac{124200k^3 - 59004k^2+2758k}{18(18k- 1)}
}  \
      \\   \hline
      \end{array}  \]
    \caption{{The $SU(3)$ Casson invariants for homology spheres $X_k$
    obtained by  $\tfrac{1}{k}$ surgery
    on  $K(2,q)$}}
    \label{table4}
     \end{table}
      \medskip

As remarked above, for $k>0$,
     $\pm \tfrac{1}{k}$ surgery on $K(2,q)$ is homeomorphic to
    the Brieskorn sphere $\Si(2,q,2qk\mp 1)$ up to a possible change of
    orientations. However,   $\la_{SU(3)}$ does not
 depend on
    the  choice of orientation.  Thus Table 
    \ref{table4} also gives
    the $SU(3)$ Casson
    invariants of $\Sigma(2,q,2qk\pm1)$ for $q=3$, $5$, $7$, and $9$.

From this data we conclude that $\la_{SU(3)}$ is not a finite type
invariant of low order.

 \begin{theorem} \label{notfinitetype}
    {\sl $\la_{SU(3)}$ is not a finite type invariant of order $\leq 6$.}
    \end{theorem}
    \noindent {\bf Proof}\qua
    We argue by contradiction.
    Suppose  $\la_{SU(3)}$ is a finite type invariant of order $\leq 6$. 
    Since $\la_{SU(3)}(S^3)=0$ and since it
    is invariant under
    change of orientation, it follows that  there exist constants
    $A$ and $B$ such that 
    \begin{equation} \label{nottrue}
    \la_{SU(3)} = A (\la_2 + 12 \la_{SU(2)}) + B
    \la_{SU(2)}^2,
    \end{equation}
    where $\la_2$ is the second Ohtsuki invariant
    \cite{lin-wang} and $\la_{SU(2)}$ is Casson's invariant \cite{am}.
    Both of these invariants satisfy surgery formulas
    (see   Theorem 4.3 in
    \cite{lin-wang} for $\la_2$) 
    and so can be computed in all the examples considered here.
   If Equation (\ref{nottrue}) were true,
   then each one of our computations would provide
   a linear constraint on $A$ and $B$. But just from surgeries on the trefoil,
    it follows that no such $A$ and $B$ exist.
    Thus $\la_{SU(3)}$ is not an invariant of finite type of order
    $\leq 6.$ 
\endproof

\noindent
{\bf Remark}\qua Stavros Garoufalidis has observed  
that our computations here prove that $\la_{SU(3)}$ is not a
finite type invariant of any order.

\subsection{Concluding remarks and open problems}

The methods we have developed apply more generally than
these computations suggest.
For example, although  we have restricted our attention
    to homology spheres obtained from surgeries on the
$(2,q)$ torus knots, the same methods apply to any Seifert
fibered homology sphere. For example, although
$\Si(2,5,7)$ is not obtained by surgery on a torus knot in $S^3,$ 
it can be described as surgery on a torus-like knot in a homology sphere
to which our main results apply. More generally,  one can compute  
the $\CC^2$--spectral flow, the  Chern--Simons invariants
and the rho invariants for Brieskorn homology spheres.
From this, one can deduce their $SU(3)$ Casson invariants
in case $p=2.$ On the other hand,
computing  $\la_{SU(3)}(\Sigma(p,q,r))$ when $p,q,r>2$ 
requires the use of perturbations and goes beyond the scope this article.
This problem will be
addressed in a later article. 

Our methods can also be used to compute $\la_{SU(3)}$ for Dehn surgeries on
knots other than
torus knots, eg, the figure eight knot.  The idea is to first notice
that one of
the surgeries on the figure eight knot gives  $\Si(2,3,7)$. 
This manifold can then  be used as a reference point
from which to calculate the invariants for other surgeries.  This is
especially interesting since most of the homology spheres obtained from
surgery on the figure eight knot are hyperbolic.  The crucial point in making
this idea work is   that our  formula for splitting the spectral flow
and the subsequent applications do
not assume the path $C_t$ of connections on $Z$ is flat.

In another direction, our technique for computing gauge
theoretic invariants can be generalized to 
groups other than $SU(2)$ and representations other than $\CC^2$.
For example, one can adapt our approach to compute
the ${\rm ad}\ su(2)$--spectral flow
which arises in Floer's instanton homology  and in  asymptotic
expansions of Witten's 3--manifold invariants (see \cite{JDGpaper}). 

One interesting and difficult   problem
is to determine   the extent to which
the  Atiyah--Patodi--Singer rho invariants  
fail to be invariant under homotopy equivalence. The results of \cite{fl} 
show that
the rho invariants of homotopy equivalent manifolds differ by a locally
constant function on the representation variety.  The cut--and--paste 
methods introduced here give a  technique to compute this difference on
the various path components of the representation variety. 
We plan to pursue this question in a future work.

In closing, we would like to mention one final
interesting problem raised by our results.
Although $\la_{SU(3)}$ is not a finite type invariant,  
it may still be possible to  express some of
the coefficients of the cubic polynomials in the numerators of
$\la_{SU(3)}(X_k)$ in Table \ref{table4} in terms of 
 the Alexander or Jones knot polynomials of the corresponding knot.  
For this problem, note that 
the denominators $2q(2qk \pm 1)$ appearing in Table \ref{table4} 
are just
the denominators of the Chern--Simons invariants of $\Si(2,q,2qk \pm 1)$.
We do not know if the Chern--Simons invariants are rational for general
homology spheres, or, alternatively, if the quantity $pq(pqk \mp 1)$ 
associated to $\pm \frac{1}{k}$ surgery on $K(p,q)$ extends naturally to define
an invariant for all homology spheres.

\end{document}